\documentclass[onefignum,onetabnum]{siamart190516}

\usepackage{lipsum}
\usepackage{amsfonts}
\usepackage{graphicx}
\usepackage{epstopdf}
\usepackage{amsmath}
\usepackage{subcaption}
\usepackage{pgfplots}

\usepackage[algo2e, ruled, noend, linesnumbered]{algorithm2e}
\usepackage{bm}
\ifpdf
  \DeclareGraphicsExtensions{.eps,.pdf,.png,.jpg}
\else
  \DeclareGraphicsExtensions{.eps}
\fi

\usepackage{placeins}
\usepackage{multirow}

\usepackage{booktabs}
\usepackage{tikz}
\usepackage[export]{adjustbox}

\newsiamremark{remark}{Remark}
\newsiamremark{hypothesis}{Hypothesis}
\crefname{hypothesis}{Hypothesis}{Hypotheses}
\newsiamthm{claim}{Claim}
\newsiamthm{prop}{Proposition}
\newsiamthm{defn}{Definition}
\newsiamthm{thm}{Theorem}
\newsiamthm{cor}{Corollary}
\newsiamthm{lem}{Lemma}

\headers{High-rank PSF Hessian approximation}{N. Alger, T. Hartland, N. Petra, and O. Ghattas}

\title{Point spread function approximation of high rank Hessians with locally supported non-negative integral kernels\thanks{To appear in SIAM Journal on Scientific Computing (SISC).\funding{This research was supported by the National Science Foundation under Grant No. DMS-1840265 and CAREER-1654311, DOD grant FA9550-21-1-0084, and DOE grants DE-SC0021239 and DE-SC0019303.}}}

\author{Nick Alger\thanks{Oden Institute for Computational Engineering and Sciences, The University of Texas at Austin, Austin, TX 
  (\email{nalger@oden.utexas.edu}).}
\and Tucker Hartland\thanks{Department of Applied Mathematics, University of California, Merced, Merced, CA. 
	(\email{thartland@ucmerced.edu}, \email{npetra@ucmerced.edu}).}
\and Noemi Petra\footnotemark[3]
\and Omar Ghattas\thanks{Oden Institute for Computational Engineering and Sciences and Walker Department of Mechanical Engineering, The University of Texas at Austin, Austin, TX 
	(\email{omar@oden.utexas.edu}).}}

\newcommand{\Aop}{\mathcal{A}}

\newcommand{\Aker}{\Phi}

\newcommand{\diraccomb}{\xi}
\newcommand{\combresponse}{\eta}

\newcommand{\impulseresponse}{\phi}

\newcommand{\febasis}{\psi}

\newcommand{\massmatrix}{\mathbf{M}}
\newcommand{\spatialvol}{V}
\newcommand{\spatialmean}{\mu}
\newcommand{\spatialcov}{\Sigma}
\newcommand{\genericdistribution}{\rho}
\newcommand{\pointbatch}{S}

\newcommand{\gdim}{d}
\newcommand{\fedim}{N}

\newcommand{\hrank}{k_h}
\newcommand{\nbatch}{n_b}
\newcommand{\numnbr}{k_{n}}

\newcommand{\nsamplepts}{m}

\newcommand{\icedomain}{\mathcal{D}}
\newcommand{\basalfriction}{{q}}
\newcommand{\normalvec}{\nu}
\newcommand{\searchdir}{{\widehat{\basalfriction}}}
\newcommand{\preconditioner}{\widetilde{H}}
\newcommand{\candidatepts}{X}

\newcommand{\candidatepoint}{x}
\newcommand{\velocity}{v}
\newcommand{\pressure}{p}
\newcommand{\rbfweight}{c}
\newcommand{\stress}{\sigma}
\newcommand{\tangentop}{T}
\newcommand{\strain}{\varepsilon}
\newcommand{\bodyforce}{f}

\newcommand{\rbf}{\varphi}
\newcommand{\fepoint}{p}
\newcommand{\stokesrobincoeff}{{s}}

\newcommand*{\horzbar}{\rule[.5ex]{2.5ex}{0.5pt}}

\newcommand{\concentration}{c} % NOTE: CHECK for collisions!!!
\newcommand{\advdomain}{\Omega}

\usepackage{amsopn}

\DeclareMathOperator{\Span}{span}

\DeclareMathOperator{\diam}{diam}

\DeclareFontFamily{U}{wncy}{}
\DeclareFontShape{U}{wncy}{m}{n}{<->wncyr10}{}
\DeclareSymbolFont{mcy}{U}{wncy}{m}{n}
\DeclareMathSymbol{\Sh}{\mathord}{mcy}{"58} 

\SetEndCharOfAlgoLine{}
\newcommand{\iFF}{\mathcal{F}}

\newcommand{\R}{\mathbb{R}}
\newcommand{\ipar}{\basalfriction} % make parameters consistent across problems

\makeatletter
\let\ftype@table\ftype@figure
\makeatother

\ifpdf
\hypersetup{
  pdftitle={Point spread function approximation of high rank Hessians with locally supported non-negative integral kernels},
  pdfauthor={N. Alger, N. Petra, T. Hartland, and O. Ghattas}
}
\fi

\begin{document}

\maketitle

% REQUIRED
\begin{abstract}
	We present an efficient matrix-free point spread function (PSF) method for approximating operators that have locally supported non-negative integral kernels. The PSF-based method computes impulse responses of the operator at scattered points, and interpolates these impulse responses to approximate entries of the integral kernel. To compute impulse responses efficiently, we apply the operator to Dirac combs associated with batches of point sources, which are chosen by solving an ellipsoid packing problem. The ability to rapidly evaluate kernel entries allows us to construct a hierarchical matrix (H-matrix) approximation of the operator. Further matrix computations are then performed with fast H-matrix methods. This end-to-end procedure is illustrated on a blur problem. We demonstrate the PSF-based method's effectiveness by using it to build preconditioners for the Hessian operator arising in two inverse problems governed by partial differential equations (PDEs): inversion for the basal friction coefficient in an ice sheet flow problem and for the initial condition in an advective-diffusive transport problem. While for many ill-posed inverse problems the Hessian of the data misfit term exhibits a low rank structure, and hence a low rank approximation is suitable, for many problems of practical interest the numerical rank of the Hessian is still large. The Hessian impulse responses on the other hand typically become more local as the numerical rank increases, which benefits the PSF-based method. Numerical results reveal that the preconditioner clusters the spectrum of the preconditioned Hessian near one, yielding roughly $5\times$--$10\times$ reductions in the required number of PDE solves, as compared to classical regularization-based preconditioning and no preconditioning. We also present a comprehensive numerical study for the influence of various parameters (that control the shape of the impulse responses and the rank of the Hessian) on the effectiveness of the advection-diffusion Hessian approximation. The results show that the PSF-based method is able to form good approximations of high-rank Hessians using only a small number of operator applications.
\end{abstract}

% REQUIRED

% REQUIRED

\section{Introduction}
\label{sec:intro}

We present an efficient \emph{matrix-free} point spread function (PSF) method for approximating operators $\Aop:L^2(\Omega) \rightarrow L^2(\Omega)'$ that have locally supported non-negative integral kernels. Here, $\Omega \subset \mathbb{R}^\gdim$ is a bounded domain, and $L^2(\Omega)'$ is the space of real-valued continuous linear functionals on $L^2(\Omega)$. By ``non-negative integral kernel,'' we mean that entries of $\Aop$'s integral kernel are non-negative numbers; this is not the same as positive semi-definiteness of $\Aop$. Such operators appear, for instance, as Hessians in optimization and inverse problems governed by partial differential equations (PDEs)~\cite{BorziSchulz11,Delosreyes15,HinzePinnauUlbrichUlbrich08}, Schur complements in Schur complement methods for solving partial differential equations and Poincare-Steklov operators in domain decomposition methods (e.g., Dirichlet-to-Neumann maps)~\cite{ChanMathew94,SmithBjorstadGropp96,ToselliWidlund04}, covariance operators in spatial statistics~\cite{ChenStein21,GentonKeyesTurkiyyah18,GeogaEtAl20,LindgrenRueLindstrom11}, and blurring operators in imaging~\cite{DenisEtAl11,NagyOleary98}. Here, ``matrix-free'' means that we may apply $\Aop$ and its transpose\footnote{Recall that $\Aop^T:L^2(\Omega)\rightarrow L^2(\Omega)'$ is the unique operator satisfying $\left(\Aop u\right)(w) = \left(\Aop^T w\right)(u)$ for all $u,w \in L^2(\Omega)$, where $\Aop u \in L^2(\Omega)'$ is the result of applying $\Aop$ to $u\in L^2(\Omega)$, and $\left(\Aop u\right)(w)$ is the result of applying that linear functional to $w \in L^2(\Omega)$, and similar for operations with $\Aop^T$. }, $\Aop^T$, to functions,
\begin{equation}
\label{eq:eval_maps}
u \mapsto\Aop u \qquad \text{and} \qquad w \mapsto \Aop^T w,
\end{equation} 
via a black box computational procedure, 
but cannot easily access entries of $\Aop$'s integral kernel. Evaluating the maps in~\eqref{eq:eval_maps} may require solving a subproblem that involves PDEs, or performing other costly computations.

\begin{figure}
	\begin{subfigure}[b]{0.49\textwidth}
		\centering
		\includegraphics[scale=0.3]{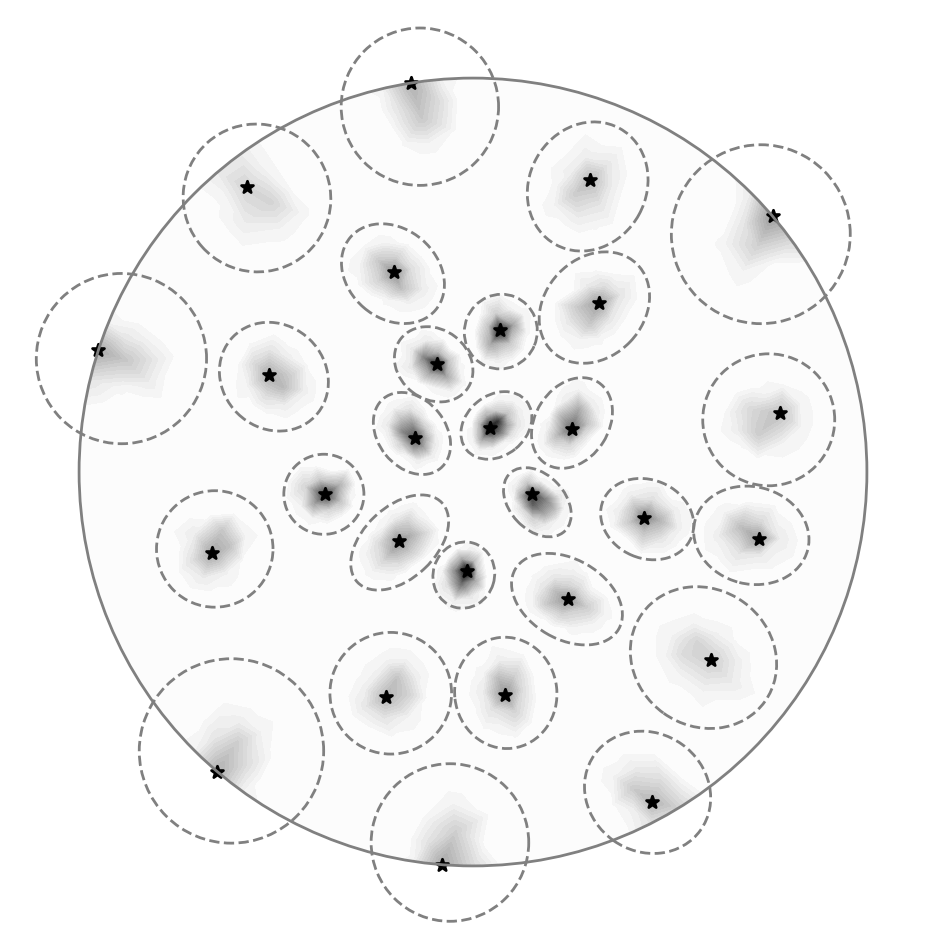}
		\caption{One impulse response batch}
		\label{fig:batches_intro}
	\end{subfigure}
	\begin{subfigure}[b]{0.49\textwidth}
		\centering
		\includegraphics[scale=0.7]{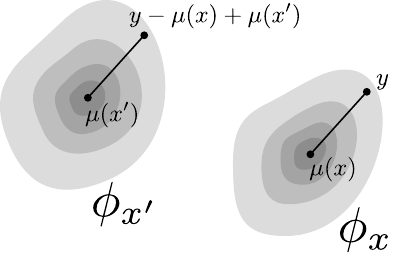}
		\caption{Local mean displacement invariance}
		\label{fig:mean_displacement_invariance}
	\end{subfigure}
	\caption{(\ref{fig:batches_intro}) One batch, $\eta_b$, of normalized impulse responses, $\impulseresponse_{x}$, that arise from applying $\Aop$ to a weighted sum of scattered point sources (see Section \ref{sec:get_impulse_response}). Here, $\Aop$ is the ice sheet inverse problem data misfit Gauss-Newton Hessian described in Section~\ref{sec:numerical_results}. 
	Black stars are point source locations. Shading shows the magnitude of the normalized impulse responses (darker means larger function values). Dashed gray ellipses are estimated impulse response support ellipsoids based on the moment method in Section \ref{sec:intromoments}. The large circle is $\partial \Omega$. (\ref{fig:mean_displacement_invariance}) Illustration of impulse responses, $\impulseresponse_x$ and $\impulseresponse_{x'}$, corresponding to points $x$ and $x'$. The operator $\Aop$ is locally mean displacement invariant (Section \ref{sec:local_mean_displacement_invariance}) if $\impulseresponse_{x}(y) \approx \impulseresponse_{x'}\left(y - \spatialmean(x) + \spatialmean(x')\right)$ when $x$ is close to $x'$. Here, $\spatialmean(z)$ denotes the mean (center of mass) of $\phi_z$.
	}
	\label{fig:batches_and_mdi}
\end{figure}

\begin{figure}
	\begin{center}
		\includegraphics[scale=0.3]{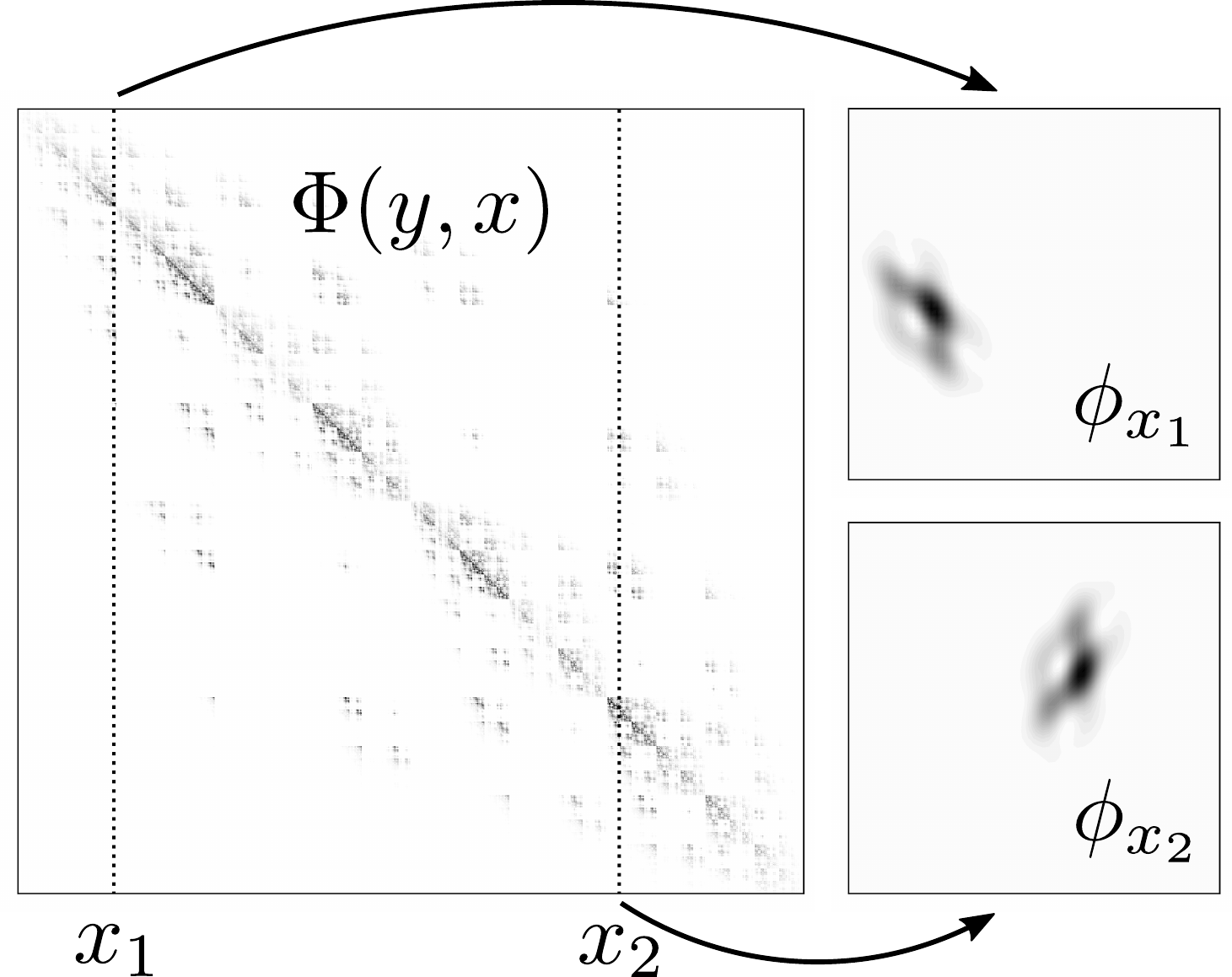}
	\end{center}
	
	\caption{Left: Matrix created by evaluating the integral kernel $\Phi$ for $\Aop$ (Equation~\ref{eq:kernel_representation}) at all pairs of mesh vertices. This illustration is for the integral kernel in Equation~\ref{eq:frog_kernel}. Dark colors indicate large entries and light colors indicate small entries. Rows and columns are ordered according to a kd-tree hierarchical clustering. Right: Impulse responses associated with points $x_1,x_2\in \Omega$, shown by the two dotted vertical lines. Intuitively, one may think of impulse responses as ``columns'' of the integral kernel.
	}
	\label{fig:frog_kernel_impulse_responses}
\end{figure}

\begin{figure}
	\begin{subfigure}[b]{0.65\textwidth}
		\centering
		\includegraphics[scale=0.22]{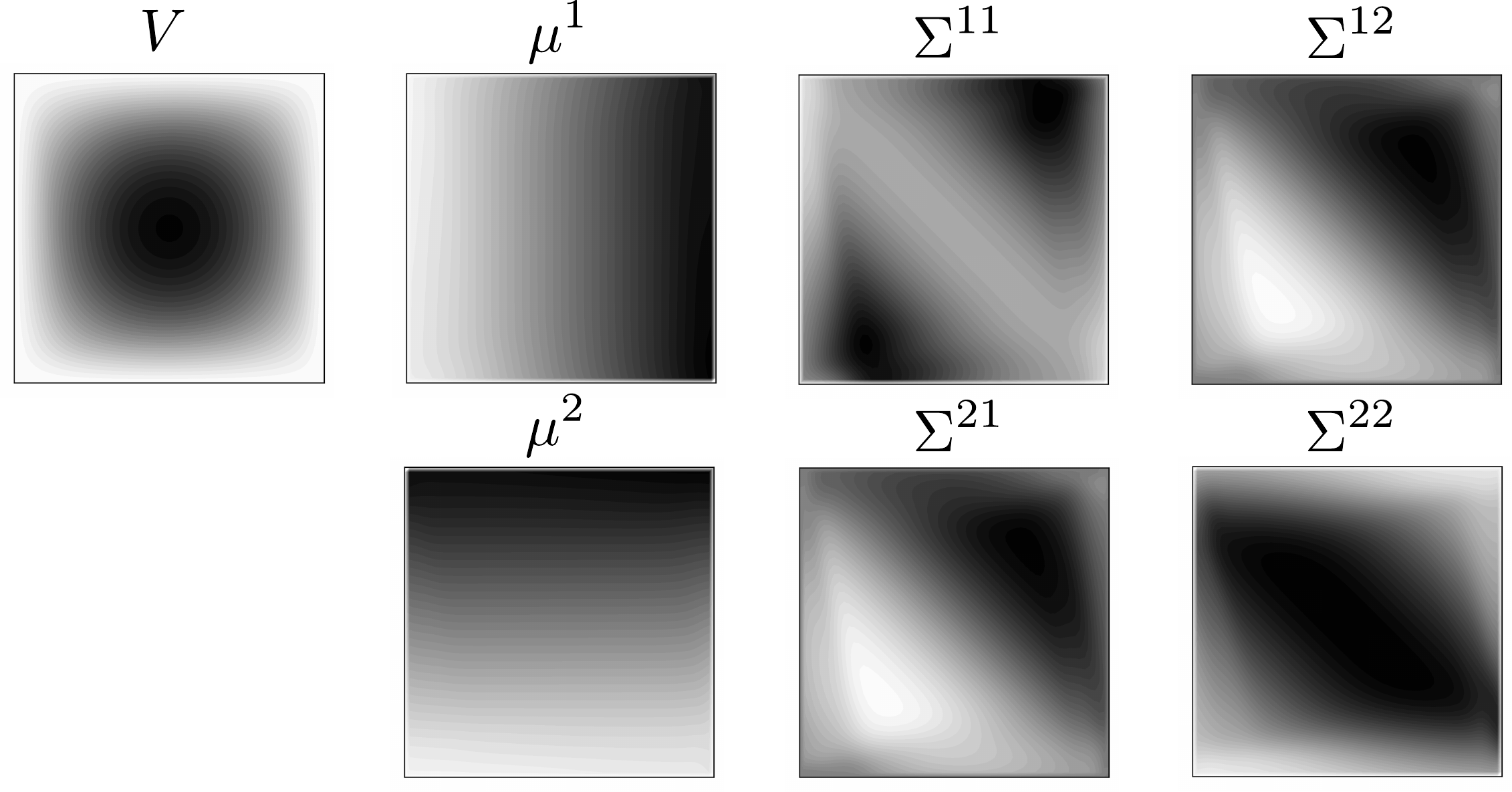}
	\end{subfigure}
	\begin{subfigure}[b]{0.34\textwidth}
		\centering
		\includegraphics[scale=0.44]{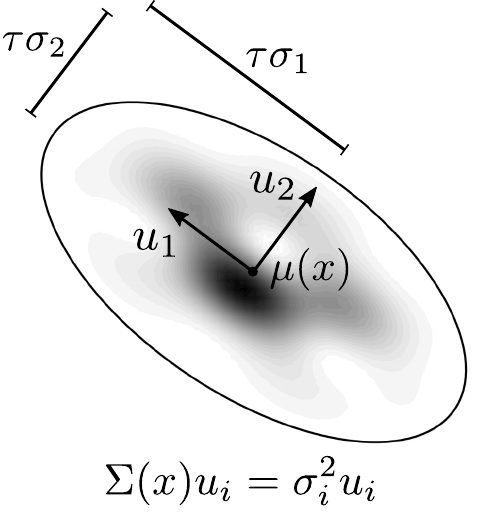}
	\end{subfigure}
	\caption{Left: Impulse response moments. Scaling factor ($V$), mean ($\mu$), and covariance ($\Sigma$). For each point $x \in \Omega$, the quantity $V(x)$ is the integral of $\phi_x$ over $\Omega$, $\mu(x)$ is the location that $\phi_x$ is centered at, and $\Sigma(x)$ is a matrix with eigenvectors and eigenvalues that characterize the width of the support of $\phi_x$ about $\mu(x)$ (see Section~\ref{sec:intromoments}). Right: Ellipsoid support for an impulse response.  This ellipsoid is the set of points within $\tau$ standard deviations of the mean of the Gaussian distribution with mean $\spatialmean(x)$ and covariance $\spatialcov(x)$. The scaling factor $V(x)$ characterizes the magnitude of $\phi_x$.
	}
	\label{fig:frog_moments_ellipsoid}
\end{figure}

\begin{figure}
	\centering
	\includegraphics[scale=0.3]{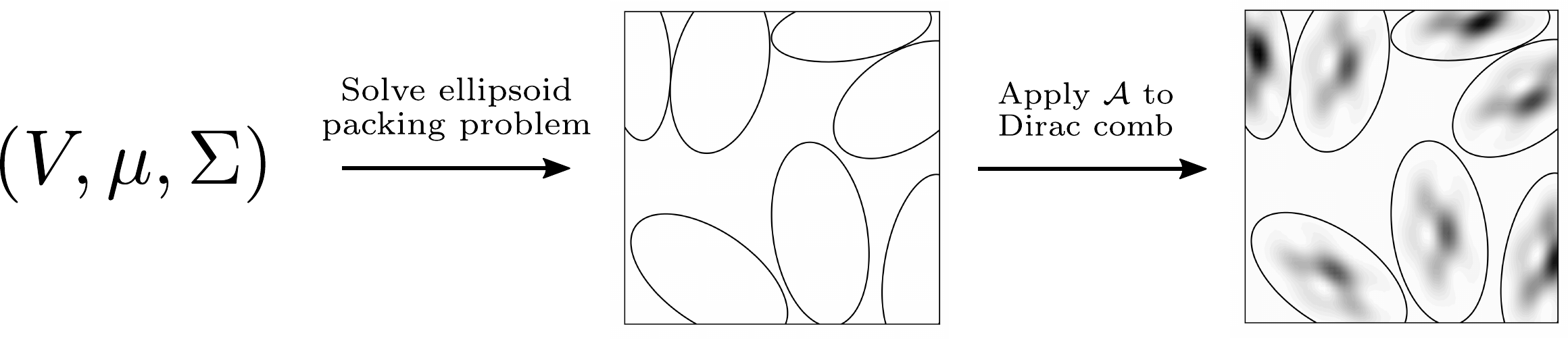}
	\caption{Illustration of the process to compute one impulse response batch. Impulse response moments are first used to form ellipsoid shaped estimates of the supports of impulse responses (Equation~\ref{eq:support_ellipsoid}). Then, an ellipsoid packing problem is solved to choose batches of non-overlapping support ellipsoids (Section~\ref{sec:sample_point_selection}). Finally, $\Aop$ is applied to a Dirac comb associated with the points $x_i$, which correspond to the ellipsoids (Section~\ref{sec:get_impulse_response}). The process is repeated to form more batches.
	}
	\label{fig:frog_batches}
\end{figure}

\begin{figure}
	\centering
	\includegraphics[scale=0.18]{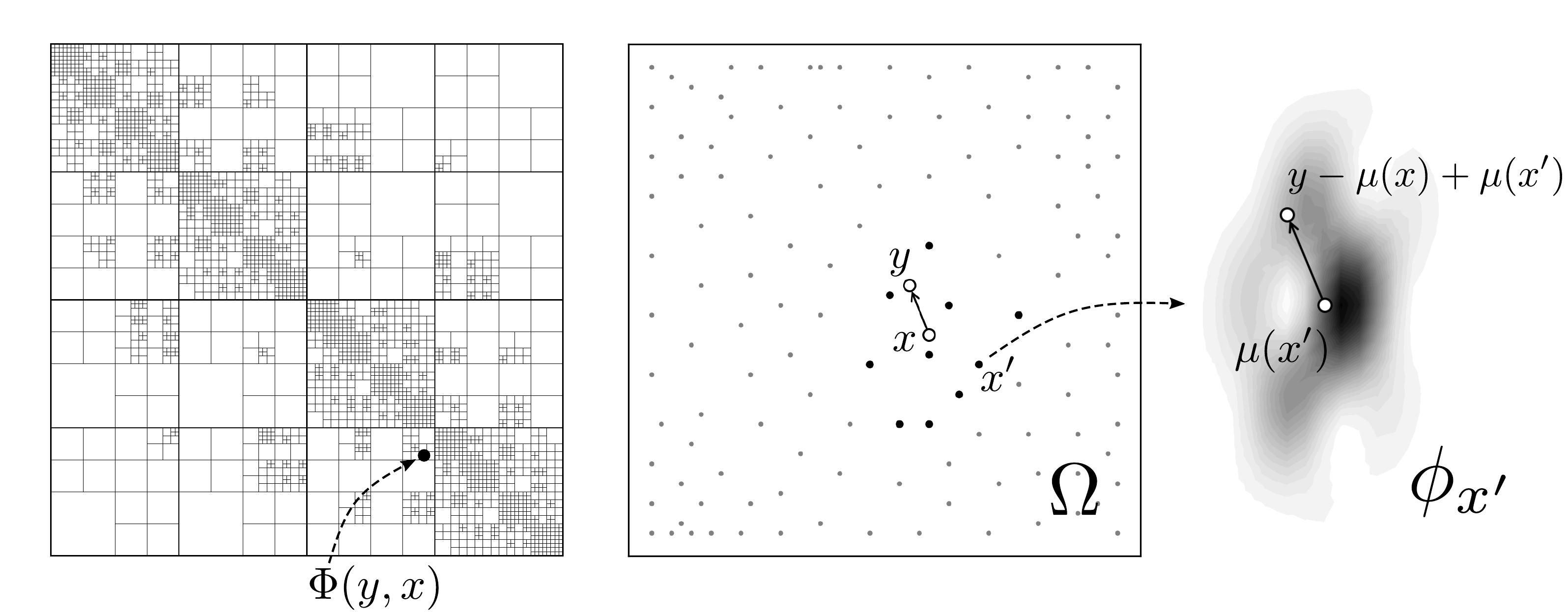}	
	\caption{Left: H-matrix structure for $\Aker$. Computing an entry of this matrix requires evaluating the integral kernel, $\Phi(y,x)$, at a pair of points $(y,x) \in \Omega \times \Omega$.  Center: Kernel evaluation points $x$ and $y$ (black circles),  sample points for the approximation (light gray and black dots), and the $k_n$ sample points, $x'$, that are nearest to $x$ (black dots). Right: Known impulse response at $x'$. Using radial basis function interpolation, the desired kernel entry is approximated as a weighted linear combination of translated and scaled versions of impulse responses at the points $x'$ (Section~\ref{sec:approximate_kernel_entries}).
	}
	\label{fig:hmatrix_neighbors_frog}
\end{figure}

\begin{figure}
	\begin{subfigure}[b]{0.30\textwidth}
		\centering
		\includegraphics[scale=0.32]{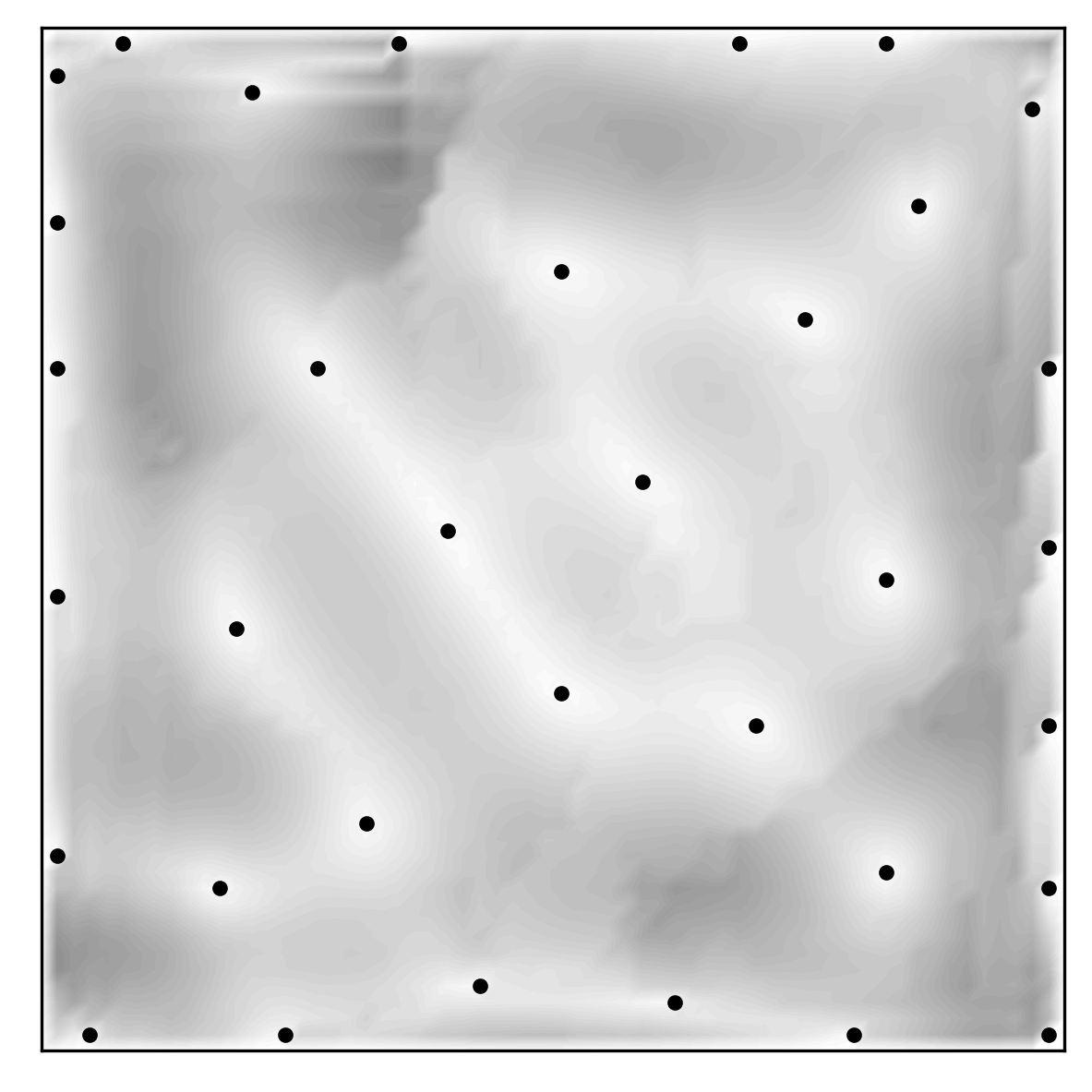}
		\label{fig:blur_error_psf5}
	\end{subfigure}
	\begin{subfigure}[b]{0.30\textwidth}
		\centering
		\includegraphics[scale=0.32]{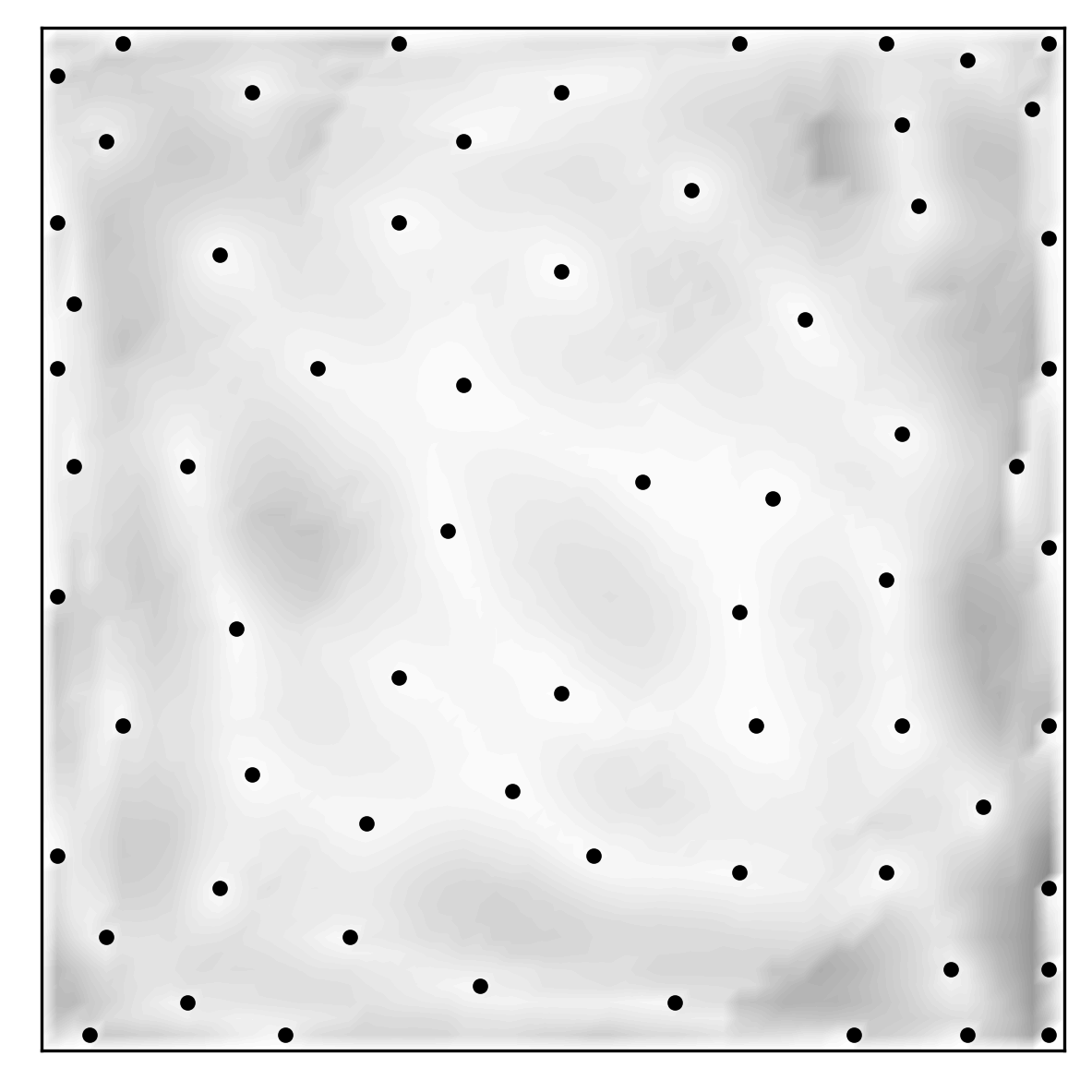}
		\label{fig:blur_error_psf10}
	\end{subfigure}
	\begin{subfigure}[b]{0.30\textwidth}
		\centering
		\includegraphics[scale=0.32]{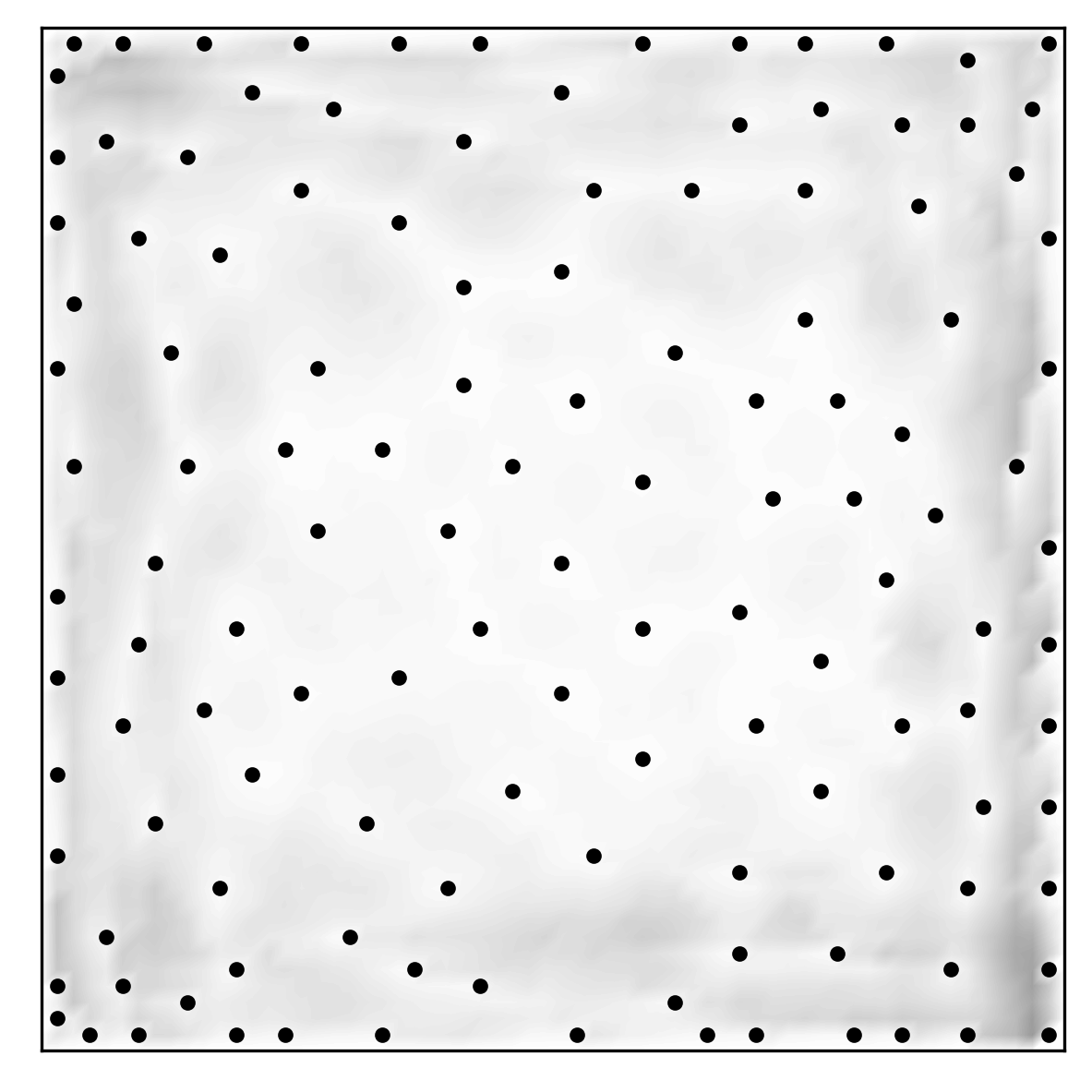}
		\label{fig:blur_error_psf20}
	\end{subfigure}
	\begin{subfigure}[b]{0.08\textwidth}
		\centering
		\includegraphics[scale=0.32]{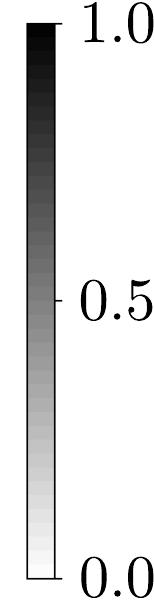}
	\end{subfigure}
	
	\caption{Relative error, $||\Phi(\cdot,x) - \widetilde{\Phi}(\cdot,x)||/||\Phi(\cdot,x)||$, in the approximation of the ``column'' of the integral kernel associated with $x$, using 5 (left), 10 (center) and 20 (right) impulse response batches. Sample points are indicated by black dots. The error associated with the point $x$ is the shade of the image at location $x$, with white indicating zero error and black indicating $100\%$ error. At the sample points, the error is zero. The further the point $x$ is from the sample points, the larger the error. Adding more batches yields a more accurate approximation.
	}
	\label{fig:frog_column_errors}
\end{figure}

The idea of the proposed method, which we refer to throughout the paper as the ``PSF-based method,'' is to use \emph{impulse response interpolation} to form a high rank approximation of $\Aop$ using a small number of operator applications. The impulse response, $\phi_x$, associated with a point, $x$, is the Riesz representation\footnote{Recall that the Riesz representative of a functional $\rho \in L^2(\Omega)'$ with respect to the $L^2$ inner product is the unique function $\rho^* \in L^2(\Omega)$ such that $\rho(w) = \left(\rho^*,w\right)_{L^2(\Omega)}$ for all $w \in L^2(\Omega)$.} of the linear functional that results from applying $\Aop$ to a delta distribution (i.e., point source, impulse) centered at $x$. We compute batches of impulse responses by applying $\Aop$ to weighted sums of delta distributions associated with batches of points scattered throughout the domain (see Figure~\ref{fig:batches_intro}). Batches of impulse responses may be thought of intuitively as sets of ``columns'' of the kernel (Figure \ref{fig:frog_kernel_impulse_responses}). To choose the batches, we form ellipsoid estimates for the supports of all $\phi_x$ via a moment method (Figure~\ref{fig:frog_moments_ellipsoid}) that involves applying $\Aop^T$ to a small number of polynomials (see Section~\ref{sec:intromoments}). We then use a greedy ellipsoid packing algorithm (Figure~\ref{fig:frog_batches}) to maximize the number of impulse responses per batch. Then we interpolate translated and scaled versions of these impulse responses to approximate entries of the operator's integral kernel (Figure~\ref{fig:hmatrix_neighbors_frog}). Adding more batches yields impulse responses at more points, increasing the approximation accuracy at the cost of one operator application per batch (Figure~\ref{fig:frog_column_errors}).

The PSF-based method we propose is loosely based on ``product convolution'' (PC) approximations, which are approximations of an operator by weighted sums of convolution operators with spatially varying weights. PC and PSF methods have a long history dating back several decades. We note the following papers (among many others) in which the convolution kernels are constructed from sampling impulse responses of the operator to scattered point sources:~\cite{Adorf94,AlgerEtAl19,BigotEscandeWeiss19,EscandeWeiss15,EscandeWeiss22,EscandeWeiss12,FishEtAl96,NagyOleary98,ZhuLiFomelEtAl16}. For background on PC and PSF methods, we recommend the following papers: \cite{DenisEtAl15,EscandeWeiss17,GentileCourbinMeylan13}. The proposed PSF-based method improves upon existing PC and PSF methods in the following ways: (1) While PC and PSF approximations are typically based on an assumption of local translation invariance, the method we propose is based on a more general assumption we call ``local mean displacement invariance'' (Section~\ref{sec:local_mean_displacement_invariance} and Figure~\ref{fig:mean_displacement_invariance}), which improves the interpolation of the impulse responses. (2) In our previous work~\cite{AlgerEtAl19}, we chose point sources in an adaptive grid via a sequential procedure; the refinements to the adaptive grid were chosen to maximally reduce the error at each step. However, in that work each point source required a separate operator application, making the previous method expensive when a large number of impulse responses is desired. In this paper, we use a new moment method (Section \ref{sec:intromoments}) which permits computation of many impulse responses (e.g., 50) per operator application. We are inspired by resolution analysis in seismic imaging, in which
$\Aop^T$ is applied to a random noise function, and the width of the support of $\phi_x$ is estimated to be the autocorrelation length of the resultant function near $x$~\cite{FichtnerLeeuwen15,TrampertFichtnerRitsema13}. The moment method that we use estimates the support of $\phi_x$ more accurately than random noise probing in resolution analysis, at the cost of the additional constraint that $\Aop$ has a non-negative integral kernel. (3) The PSF-based method we propose never evaluates computed impulse responses outside of their domain of definition. This eliminates ``boundary-artifact'' errors (see \cite[Section 1.1]{AlgerEtAl19}) that plague conventional PC and PSF methods.

The ability to rapidly approximate entries of $\Aop$'s integral kernel allows one to approximate discretized versions of $\Aop$ using the full arsenal of tools for matrix approximation that rely on fast access to matrix entries.
In this work, we form a hierarchical matrix~\cite{BormGrasedyckHackbusch03,Hackbusch99} approximation of a discretized version of $\Aop$. H-matrices are a compressed matrix format in which the rows and columns of the matrix are re-ordered, then the matrix is recursively subdivided into blocks in such a way that many off-diagonal blocks are low rank, even though the matrix as a whole may be high rank. 
H-matrix methods permit us to perform matrix-vector products cheaply, and perform other useful linear algebra operations that cannot be done easily using the original operator. These operations include matrix-matrix addition, matrix-matrix multiplication, matrix factorization, and matrix inversion. 
The work and memory required to perform these operations for an $\fedim \times \fedim$ H-matrix with rank $\hrank$ blocks scales as $O\left(\hrank^a \fedim \log(\fedim)^b\right)$ where $a,b \in \{0,1,2,3\}$ are constants which depend on the type of H-matrix used and the operation being performed~\cite{GrasedyckHackbusch03}\cite[Section 2.1]{Kriemann13}.

\section{Why we need more efficient approximations of high rank Hessians}
\label{sec:hessian}

While the PSF-based method proposed in this paper may be used to approximate any operator that has a locally supported non-negative integral kernel, we are primarily motivated by approximation of high-rank Hessians in distributed parameter inverse problems governed by PDEs. In this section, we provide a brief background on this topic, and explain why existing Hessian approximation methods are not satisfactory.

In distributed parameter inverse problems governed by PDEs, one seeks to infer an unknown spatially varying parameter field from limited observations of a state variable that depends on the parameter implicitly through the solution of a PDE. 
Conventionally, the inverse problem is formulated using either a deterministic framework~\cite{BanksKunisch89,Vogel02},
or a Bayesian probabilistic framework~\cite{KaipioSomersalo05,Stuart10a,Tarantola05}. In the deterministic framework, one solves an optimization problem to find the parameter that best fits the observations, subject to appropriate regularization~\cite{EnglHankeNeubauer96,Vogel02}. In the probabilistic framework, Bayes' theorem combines the observations with prior information to form a posterior distribution over the space of all possible parameter fields, and computations are performed to extract statistical information about the parameter from this posterior. The Hessian of the objective function with respect to the parameter in the determinstic optimization problem and the Hessian of the negative log posterior in the Bayesian setting are equal or approximately equal under typical noise, regularization, and prior models, so we refer to both of these Hessians as ``the Hessian.'' The Hessian consists of a data misfit term (the \emph{data misfit Hessian}), which depends on a discrepancy between the observations and the associated model predictions, and a regularization or prior term (the \emph{regularization Hessian}) which does not depend on the observations. For more details on the Hessian, see~\cite{Alger19,GhattasWillcox21,VillaPetraGhattas21}. 

Hessian approximations and preconditioners are highly desirable because the Hessian is central to efficient solution of inverse problems in both the deterministic and Bayesian settings. When solving the deterministic optimization problem with Newton-type methods, the Hessian is the coefficient operator for the linear system that must be solved or approximately solved at every Newton iteration. Good Hessian preconditioners reduce the number of iterations required to solve these Newton linear systems with the conjugate gradient method~\cite{Saad03}. In the Bayesian setting, the inverse of the Hessian is the covariance of a local Gaussian approximation of the posterior. This Gaussian distribution can be used directly as an approximation of the posterior, or it can be used as a proposal for Markov chain Monte-Carlo methods for drawing samples from the posterior. For instance, see~\cite{KimEtAl23,PetraEtAl14} and the references therein. 

Owing to the implicit dependence of predicted observations on the parameter, entries of the Hessian are not easily accessible. Rather, the Hessian may be applied to a vector via a computational process that involves solving a pair of forward and adjoint PDEs which are linearizations of the original PDE~\cite{GhattasWillcox21,PetraStadler11}.
The most popular matrix-free Hessian approximation methods are based on low rank approximation of either the data misfit Hessian, or the data misfit Hessian preconditioned by the regularization Hessian, e.g.,~\cite{BuiEtAl13,CuiEtAl14,FlathEtAl11,PetraEtAl14,SpantiniEtAl15}. Krylov methods such as Lanczos or randomized methods~\cite{Cheng05,HalkoMartinssonTropp11} are typically used to construct these low rank approximations by applying the Hessian to vectors. Using these methods, the required number of Hessian applications (and hence the required number of PDE solves) is proportional to the rank of the low rank approximation. Low rank approximation methods are justified by arguing that the numerical rank of the data misfit Hessian is insensitive to the dimension of the discretized parameter. This means that the required number of PDE solves remains the same as the mesh used to discretize the parameter is refined.
However, in many inverse problems of practical interest the numerical rank of the data misfit Hessian, while mesh independent, is still large, which makes it costly to approximate the Hessian using low rank approximation methods~\cite{AmbartsumyanEtAl20,BuiEtAl13,IsaacEtAl15}. 

Examples of inverse problems with high rank data misfit Hessians include large-scale ice sheet inverse problems \cite{Hartland_2023,IsaacEtAl15}, advection dominated advection-diffusion inverse problems~\cite{AkcelikBirosDraganescuEtAl05}\cite[Chapter 5]{Flath13}, high frequency wave propagation inverse problems~\cite{BuiEtAl13}, inverse problems governed by high Reynolds number flows, and more generally, all inverse problems in which the observations highly inform the parameter. The eigenvalues of the data misfit Hessian characterize how informative the data are about components of the parameter in the corresponding eigenvector directions,
hence more informative data leads to larger eigenvalues and a larger numerical rank~\cite{AlexanderianGloorGhattas16}\cite[Section 1.4 and Chapter 4]{Alger19}. Roughly speaking, the numerical rank of the data misfit Hessian is the dimension of the subspace of parameter space that is informed by the data. The numerical rank of the regularization preconditioned data misfit Hessian may be reduced by increasing the strength of the regularization, but this throws away useful information: components of the parameter that could be learned from the observations
would instead be reconstructed based on the regularization~\cite[Section 4]{AlgerEtAl17}\cite[Chapters 1 and 7]{Vogel02}. Hence, low rank approximation methods suffer from a predicament: if the data highly inform the parameter and the regularization is chosen appropriately, then a large number of operator applications are required to form an accurate approximation of the Hessian using low rank approximation methods. 
High rank Hessian approximation methods are thus needed.

Recently there have been improvements in matrix-free H-matrix construction methods in which an operator is applied to structured random vectors, and the response of the operator to those random vectors is processed to construct an H-matrix approximation~\cite{LevittMartinsson22,LinYing11,Martinsson11,Martinsson16,MartinssonTropp20}. 
These methods (which we do not use here) have been used to approximate Hessians in PDE constrained inverse problems~\cite{AmbartsumyanEtAl20,Hartland_2023}. 
Although these methods are promising,
the required number of operator applications is still large (e.g., hundreds to thousands). For example, using the method in \cite{LinYing11}, the required number of operator applies to construct an $H^1$ matrix with hierarchical rank $r$ for problems in a 2D domain discretized with a regular grid is $\#\text{levels} \cdot 64 \cdot (r + c)$, where $\#\text{levels}$ is the depth of the hierarchical partitioning, $r$ is the rank of the blocks (hierarchical rank), and $c$ is an oversampling parameter (see \cite[Section 2.4]{LinYing11}). On a $64 \times 64$ grid with depth $4$, hierarchical rank $10$, and oversampling parameter $c=5$, this works out to $4 \cdot 64 \cdot (10 + 5) = 3840$ operator applies. In Section \ref{sec:frog}, we see numerically that the randomized hierarchical off-diagonal low rank (HODLR) method in \cite{Martinsson16} requires hundreds to thousands of matrix-vector products to construct approximations of the integral kernel for a blur problem example with modest (e.g., $10\%$) relative error.
Matrix-free H-matrix construction is currently an active area of research, hence these costs may decrease as new algorithms are developed.
In this paper, we also form an H-matrix approximation. However, to reduce the required number of operator applications, we first form a PSF approximation of the data misfit Hessian by exploiting locality and non-negative integral kernel properties, then form the H-matrix using classical techniques. Using this two stage approach, we reduce the number of operator applications to a few dozen at most. 

Not all data misfit Hessians satisfy the local non-negative integral kernel properties. We note, in particular, that the wave inverse problem data misfit Hessian and Gauss-Newton Hessian have a substantial proportion of negative entries in their integral kernels. In this case more specialized techniques have been developed using, eg., pseudodifferential operator theory \cite{DemanetEtAl12,Symes08}, and sparsity in the wavelet domain \cite{HerrmannMoghaddamStolk08}. However, many data misfit Hessians of practical interest do satisfy the local non-negative integral kernel properties (either exactly or approximately), and the PSF-based method we propose is targeted at approximating these Hessians.

\section{Preliminaries}
\label{sec:prelims}

Let $\Omega \subset \mathbb{R}^\gdim$ be a bounded domain (typically $\gdim=1$, $2$, or $3$). We seek to approximate integral operators $\Aop:L^2(\Omega)\rightarrow L^2(\Omega)'$ of the form
\begin{equation}
\label{eq:kernel_representation}
(\Aop u)(w) := \int_\Omega \int_\Omega w(y) \Aker(y,x) u(x) \mathrm{d}x \mathrm{d}y.
\end{equation}
The linear functional $\Aop u \in L^2(\Omega)'$ is the result of applying $\Aop$ to $u\in L^2(\Omega)$, and the scalar $\left(\Aop u\right)(w)$ is the result of applying that linear functional to $w \in L^2(\Omega)$.
The integral kernel, $\Aker:\Omega \times \Omega \rightarrow \mathbb{R}$, exists but is not easily accessible. In this section we describe how to extend the domain of $\Aop$ to distributions, which allows us to define impulse responses (Section~\ref{sec:impulse_response_2}), we then state the conditions on $\Aop$ that the PSF-based method requires (Section~\ref{sec:conditions_2}), and detail finite element discretization (Section~\ref{sec:finite_element_kernel}).

\subsection{Distributions and impulse responses}
\label{sec:impulse_response_2}

The operator $\Aop$ may be applied to distributions\footnote{I.e., generalized functions such as the Dirac delta distribution. See, for example, \cite[Chapter 5]{ArbogastBona08}.} if $\Aker$ is sufficiently regular. 
Given $\rho \in L^2(\Omega)'$, let $\rho^* \in L^2(\Omega)$ denote the Riesz representative of $\rho$ with respect to the $L^2(\Omega)$ inner product. We have
\begin{subequations}
	\label{eq:extension_to_distributions}
	\begin{align}
		\left(\Aop \genericdistribution^*\right)(w) &= \int_\Omega \int_\Omega w(y) \Aker(y,x) \genericdistribution^*(x) \mathrm{d}x ~\mathrm{d}y \\
		&= \int_\Omega w(y) \int_\Omega \Aker(y,x) \genericdistribution^*(x) \mathrm{d}x ~\mathrm{d}y 
		= \int_\Omega w(y) \genericdistribution\left(\Aker(y,\cdot)\right) \mathrm{d}y, \label{eq:action_on_distribution}
	\end{align}
\end{subequations}
where $\Aker(y,\cdot)$ denotes the function $x \mapsto \Aker(y,x)$.
Now let $\mathcal{D}(\Omega) \subset L^2(\Omega)$ be a suitable space of test functions and let $\rho:\mathcal{D}(\Omega) \rightarrow \mathbb{R}$ be a distribution. In this case, 
$\rho^*$ may not exist, so the derivation in~\eqref{eq:extension_to_distributions} is not valid. However, if $\Aker$ is sufficiently regular such that the function $y \mapsto \rho\left(\Aker(y,~\cdot~)\right)$ is well-defined for almost all $y\in\Omega$, and if this function is in $L^2(\Omega)$, then the right hand side of~\eqref{eq:action_on_distribution} is well-defined. Hence, we \emph{define} the application of $\Aop$ to the distribution $\rho$ to be the right hand side of~\eqref{eq:action_on_distribution}. We denote this operator application by ``$\Aop \genericdistribution^*$,'' even if $\rho^*$ does not exist.

Let $\delta_x$ denote the delta distribution\footnote{Recall that the delta distribution $\delta_x:\mathcal{D}(\Omega)\rightarrow \mathbb{R}$ is defined by $\delta_x(w) = w(x)$ for all $w\in \mathcal{D}(\Omega)$.} (i.e., point source, impulse) centered at the point $x \in \Omega$.
The \emph{impulse response} of $\Aop$ associated with $x$ is the function $\impulseresponse_x:\Omega \rightarrow \mathbb{R}$, 
\begin{equation}
\label{eq:impulse_response_delta_action}
\impulseresponse_x := \left( \Aop \delta_x^* \right)^*,
\end{equation}
that is formed by applying $\Aop$ to $\delta_x$ (per the generalized notion of operator ``application'' defined above), then taking the Riesz representation of the resulting linear functional. Using~\eqref{eq:action_on_distribution} and the definition of the delta distribution, we see that $\phi_x$ may also be written as the function $\impulseresponse_x(y) = \Aker(y, x)$.

\subsection{Required conditions}
\label{sec:conditions_2}

We focus on approximating operators that satisfy the following conditions: 
\begin{enumerate}
	\item The kernel $\Aker$ is sufficiently regular so that $\phi_x$ is well-defined for all $x\in \Omega$.
	\item The supports of the impulse responses $\impulseresponse_x$ are contained in localized regions.
	\item The integral kernel is non-negative\footnote{Note that having a non-negative integral kernel is different from positive semi-definiteness. The operator $\Aop$ need not be positive semi-definite to use the PSF-based method, and positive semi-definite operators need not have a non-negative integral kernel.} in the sense that $$\Aker(y,x) \ge 0 \quad \text{for all} \quad (y,x) \in \Omega \times \Omega.$$
\end{enumerate}
The PSF-based method may still perform well if these conditions are relaxed slightly. 

It is acceptable if the support of $\impulseresponse_x$ is not perfectly contained in a localized region (violating Assumption 2), so long as the bulk of the ``mass'' of $\impulseresponse_x$ is contained in a localized region. In principle, the PSF-based method can be applied even if the impulse responses are widely dispersed. However, in this case only a small number of impulse responses can be computed per batch, which means more batches, and hence more operator applies, are needed to form an accurate approximation.

\begin{figure}
	\centering
	{
		\begin{tabular}{cccc}
			\includegraphics[scale=0.25, valign=c]{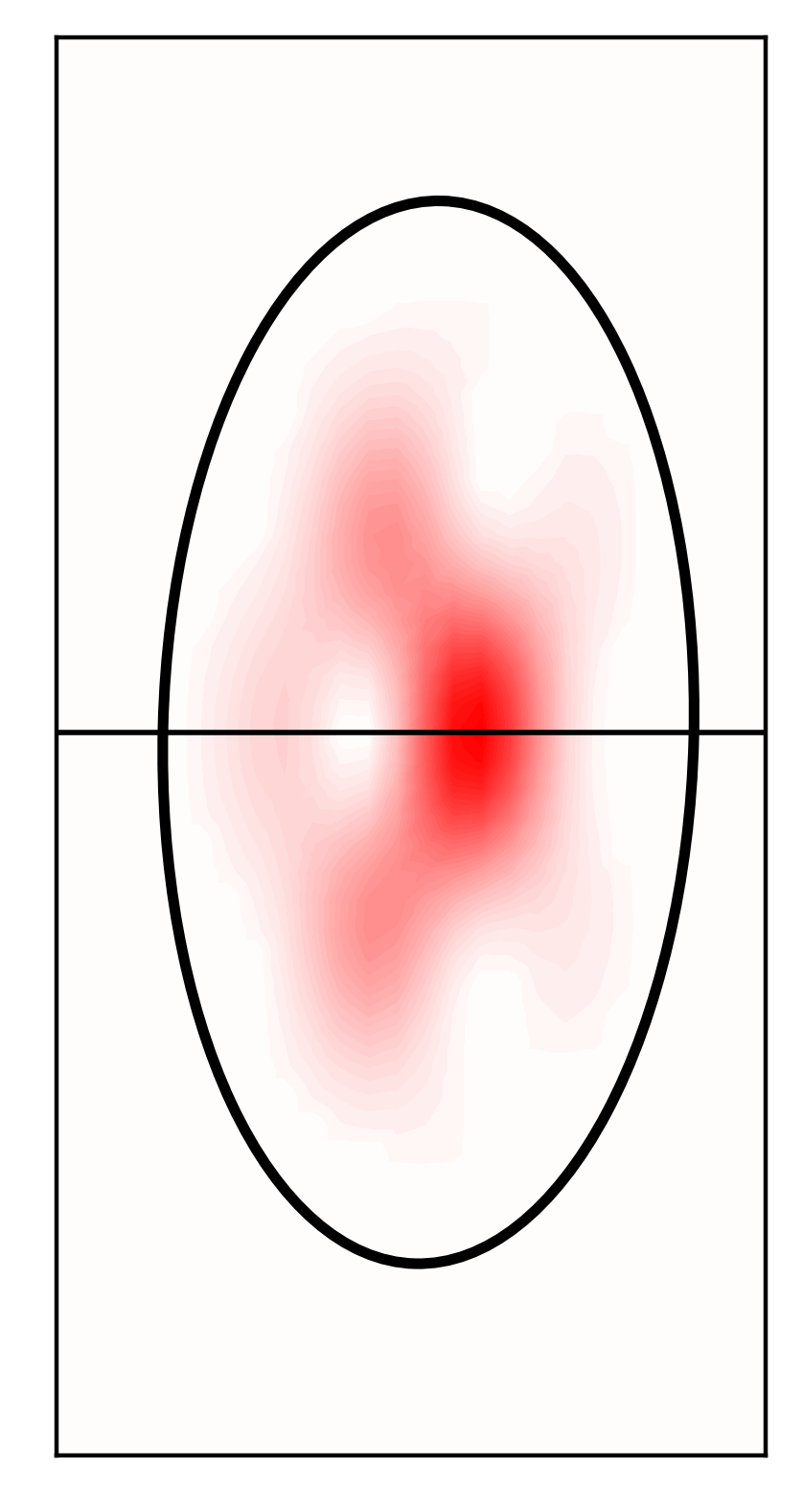} & 
			\includegraphics[scale=0.25, valign=c]{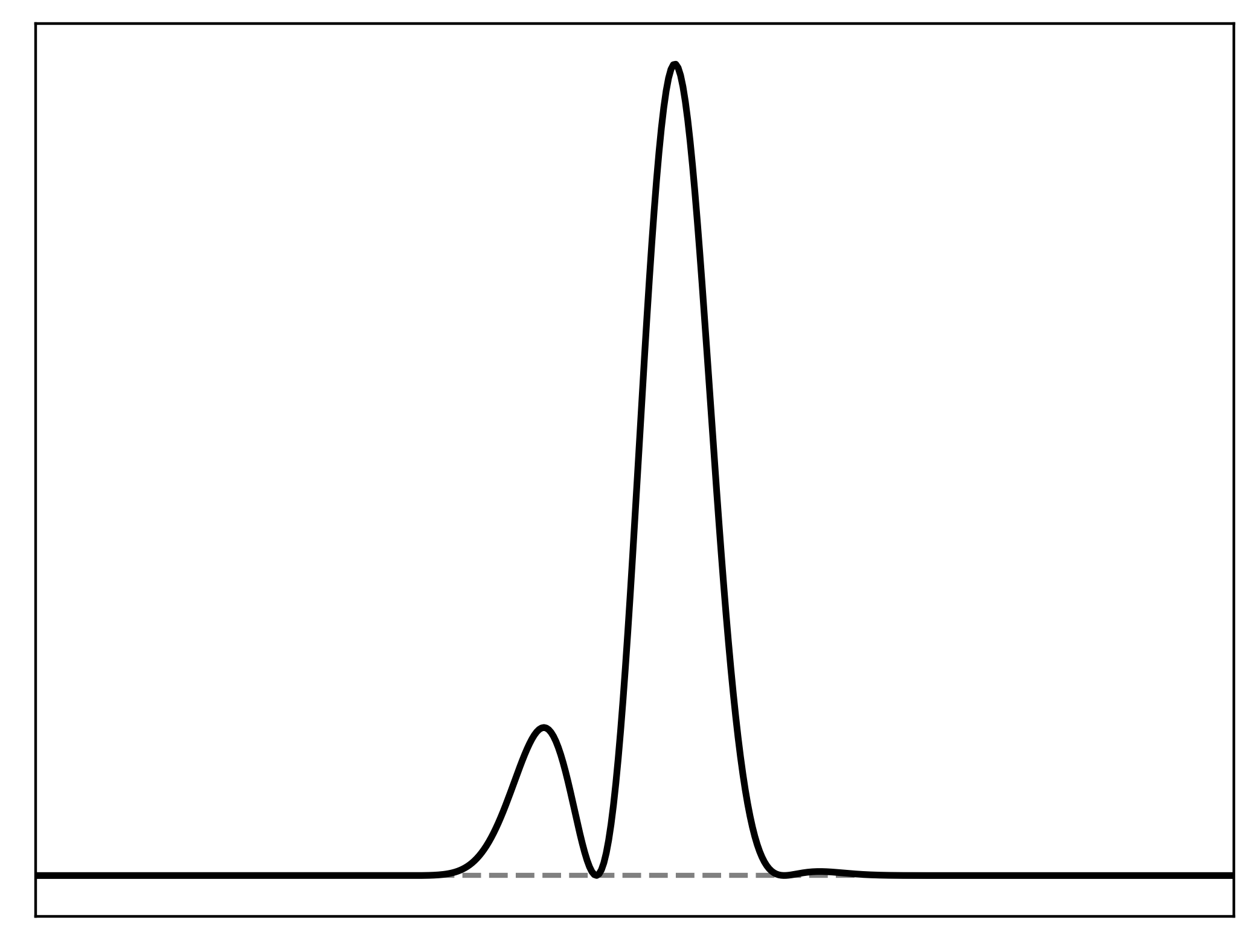} & 
			\includegraphics[scale=0.25, valign=c]{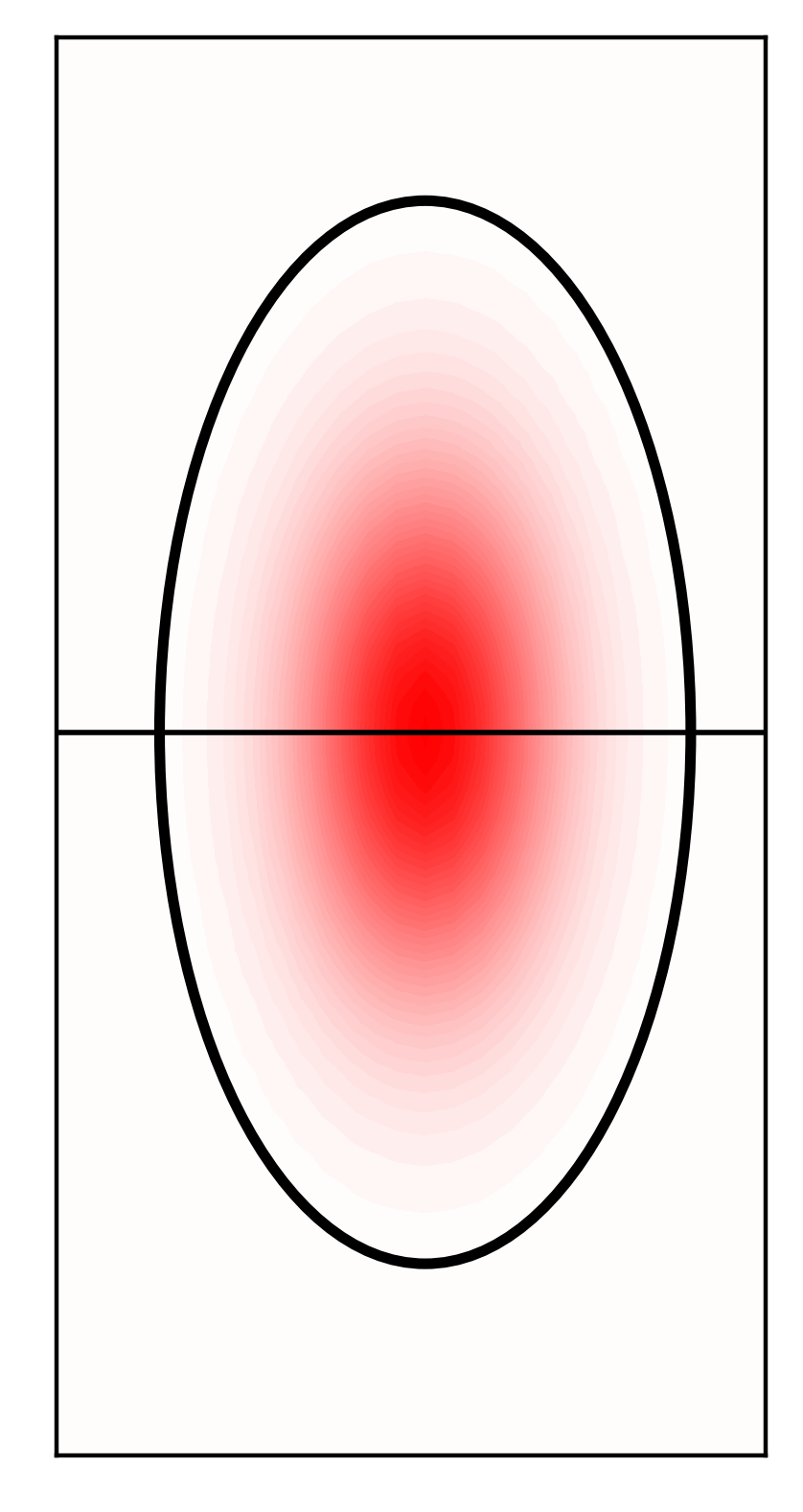} & 
			\includegraphics[scale=0.25, valign=c]{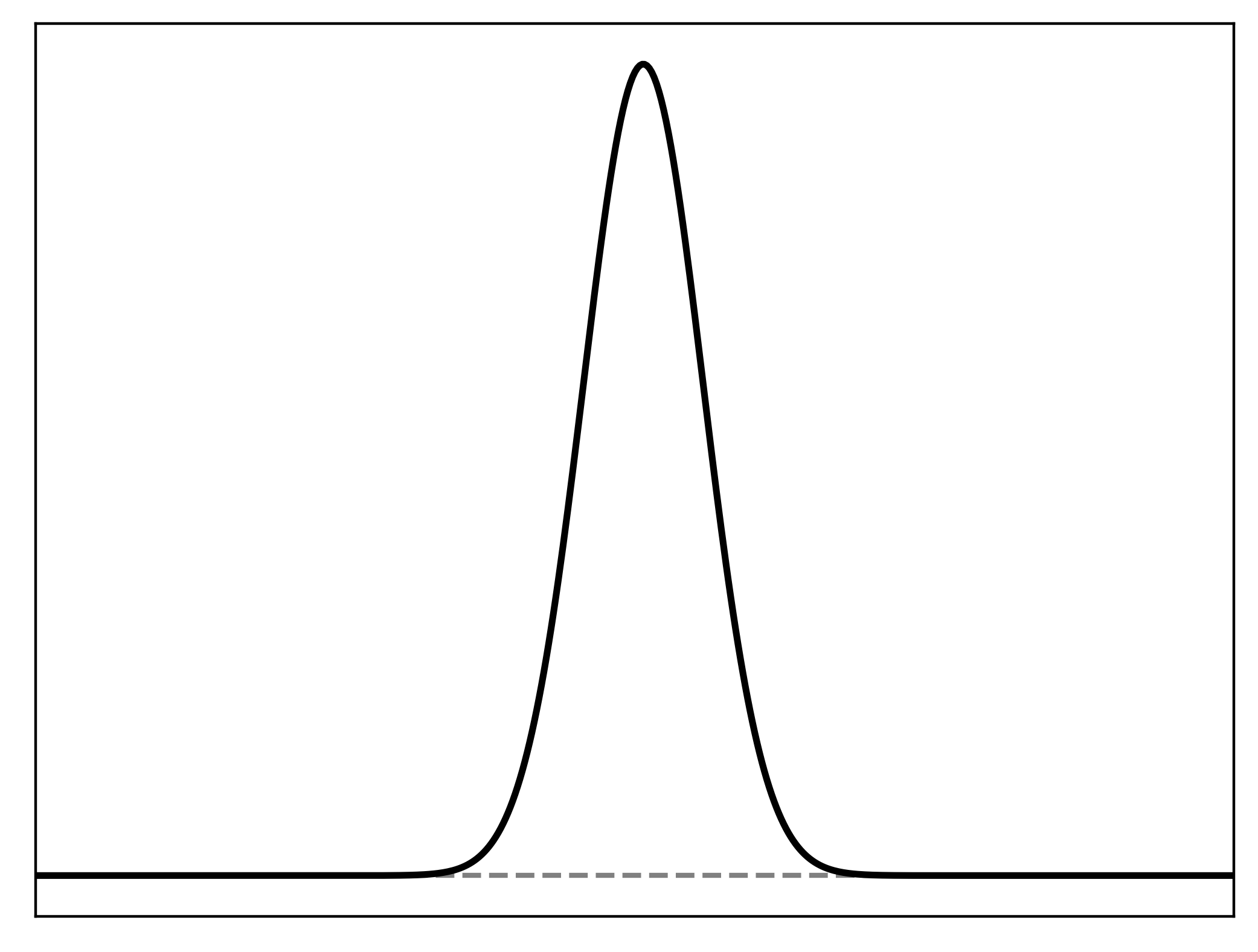} \\
			\includegraphics[scale=0.25, valign=c]{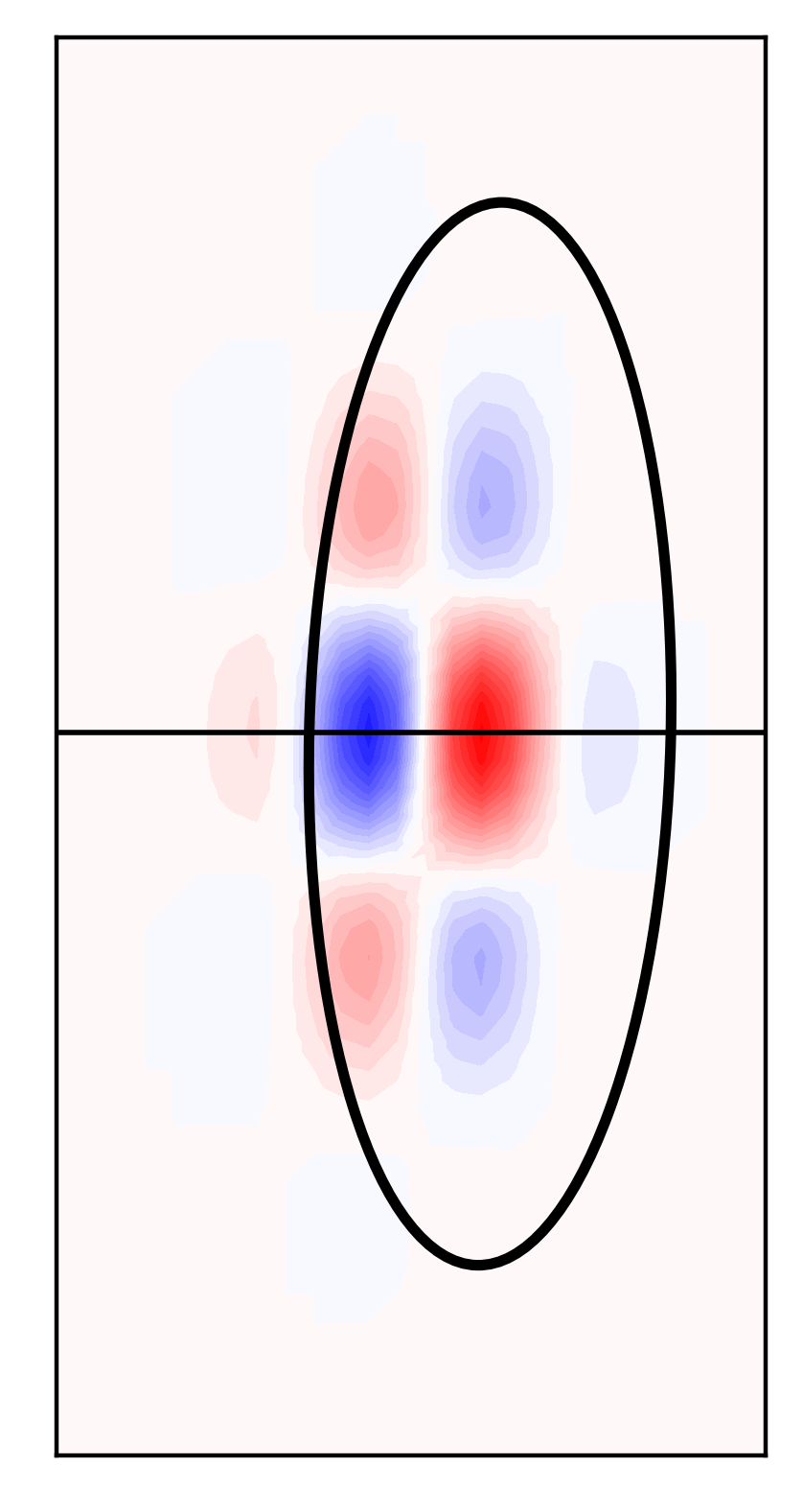} & 
			\includegraphics[scale=0.25, valign=c]{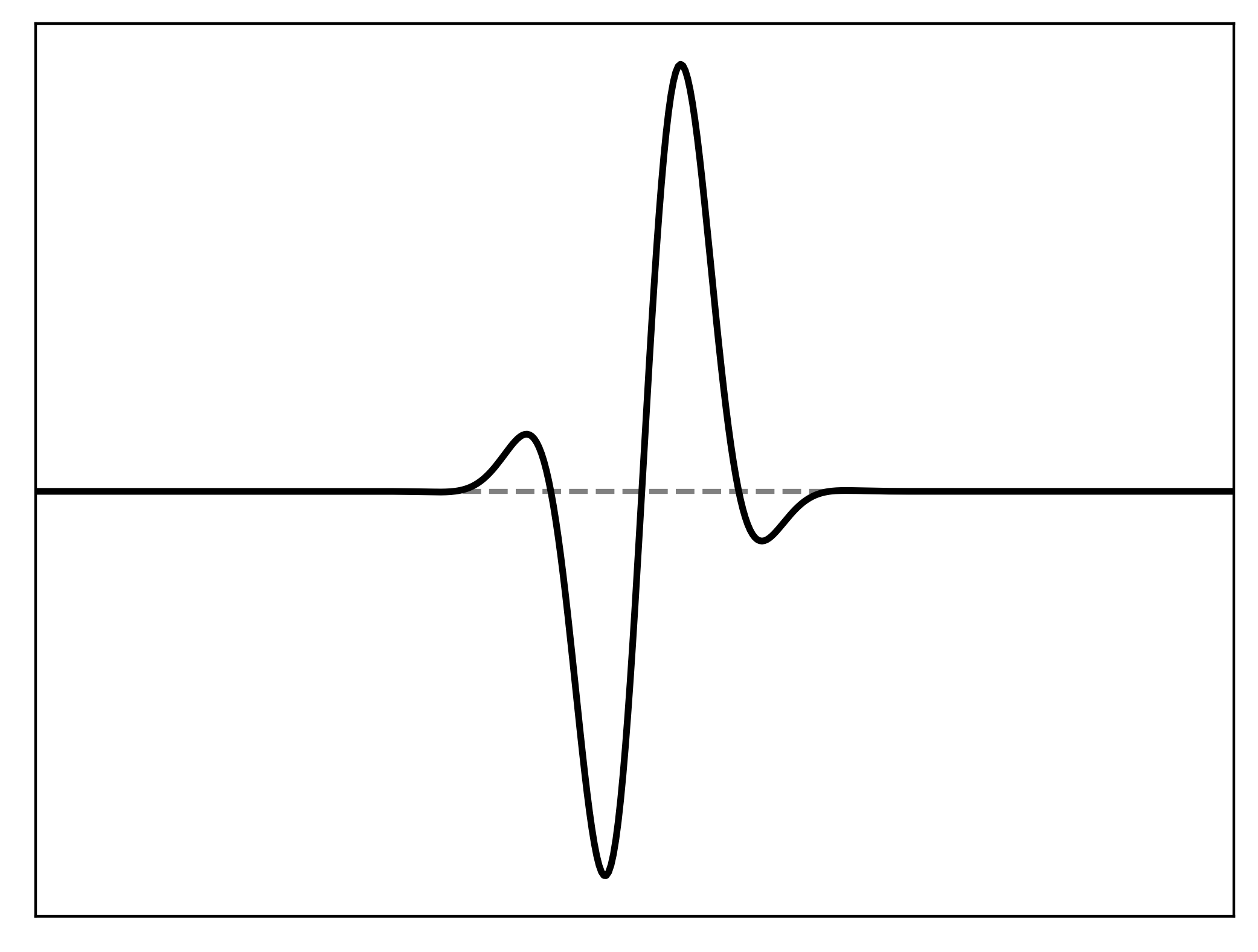} & 
			\includegraphics[scale=0.25, valign=c]{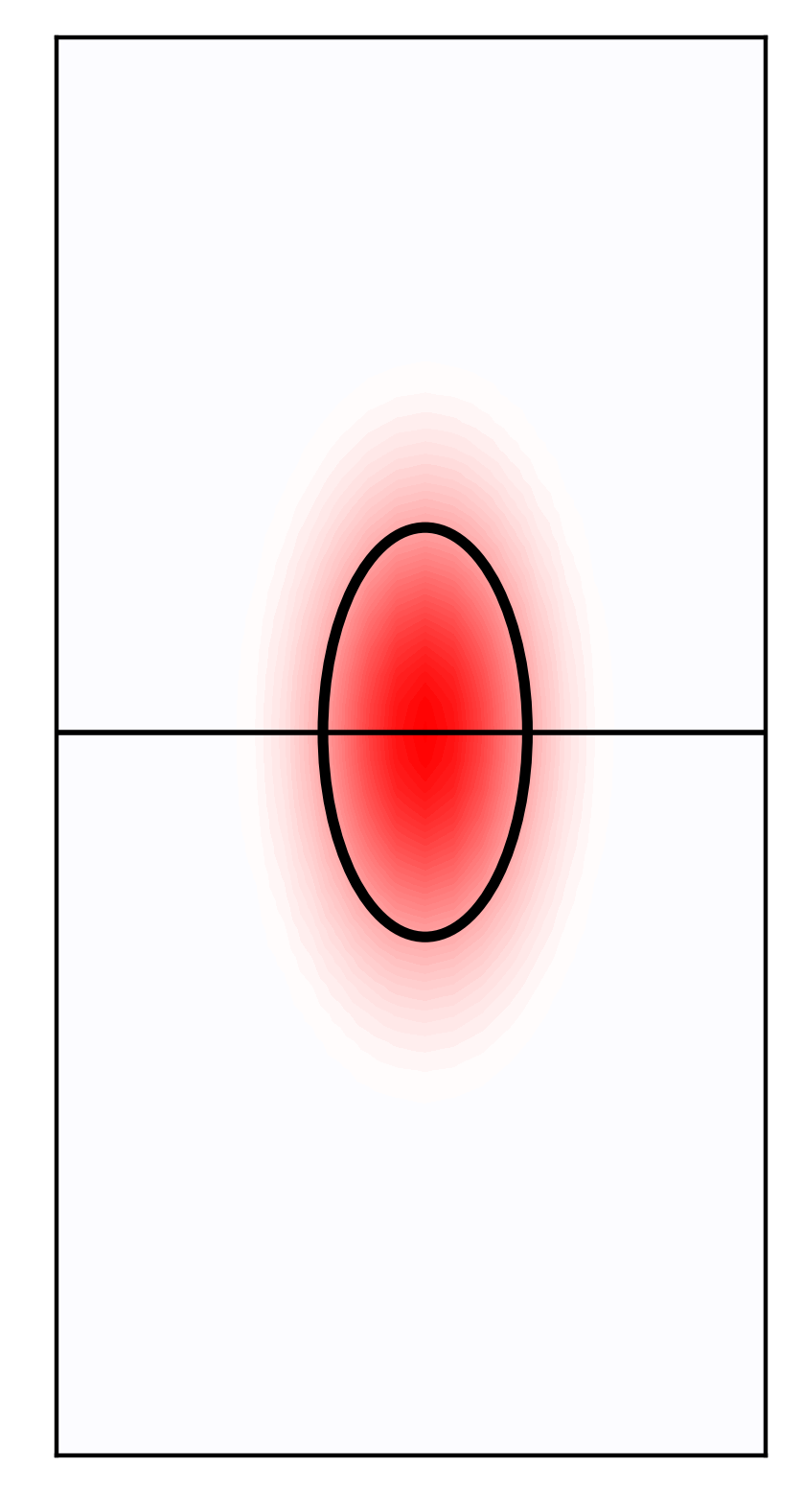} & 
			\includegraphics[scale=0.25, valign=c]{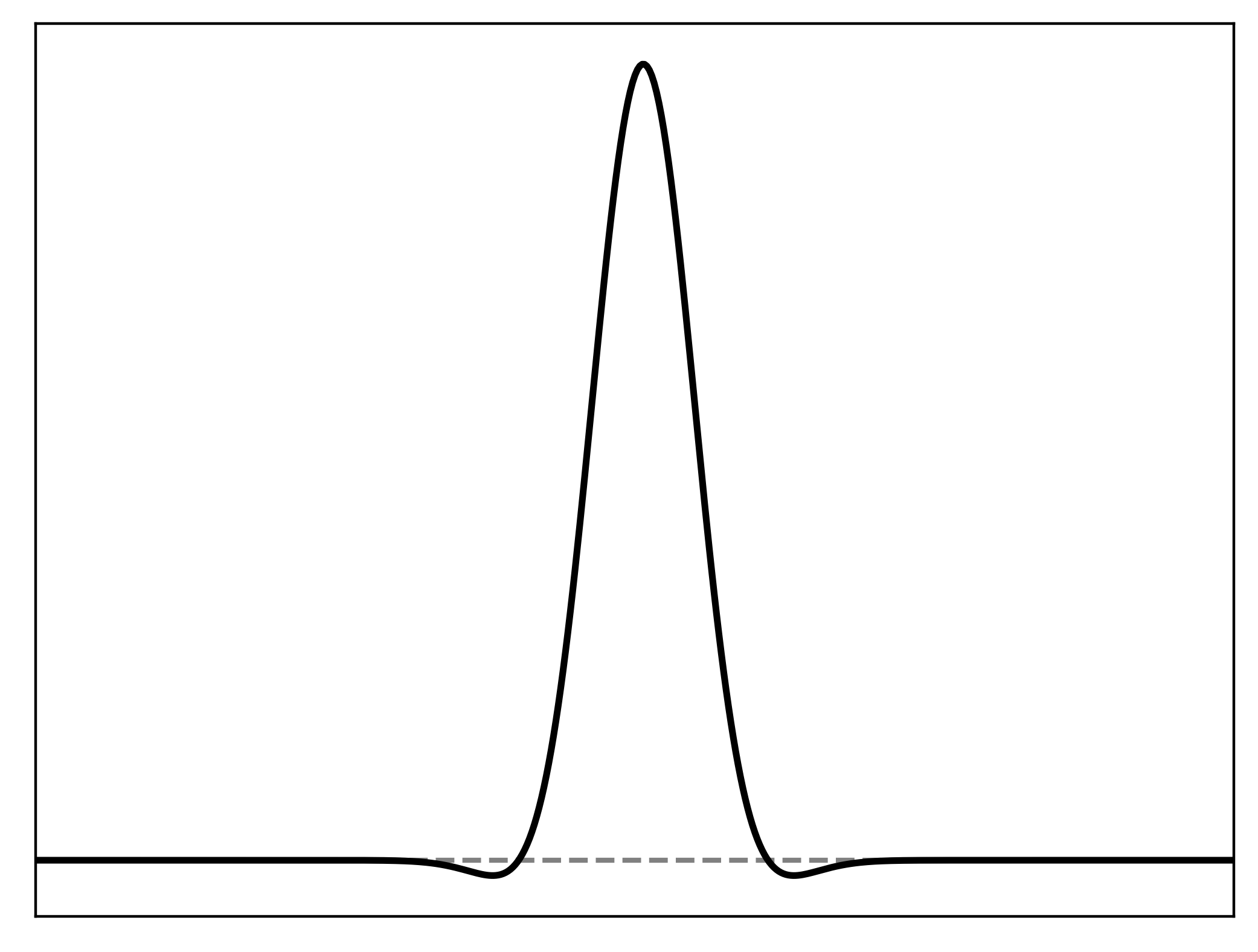} \\
			\includegraphics[scale=0.25, valign=c]{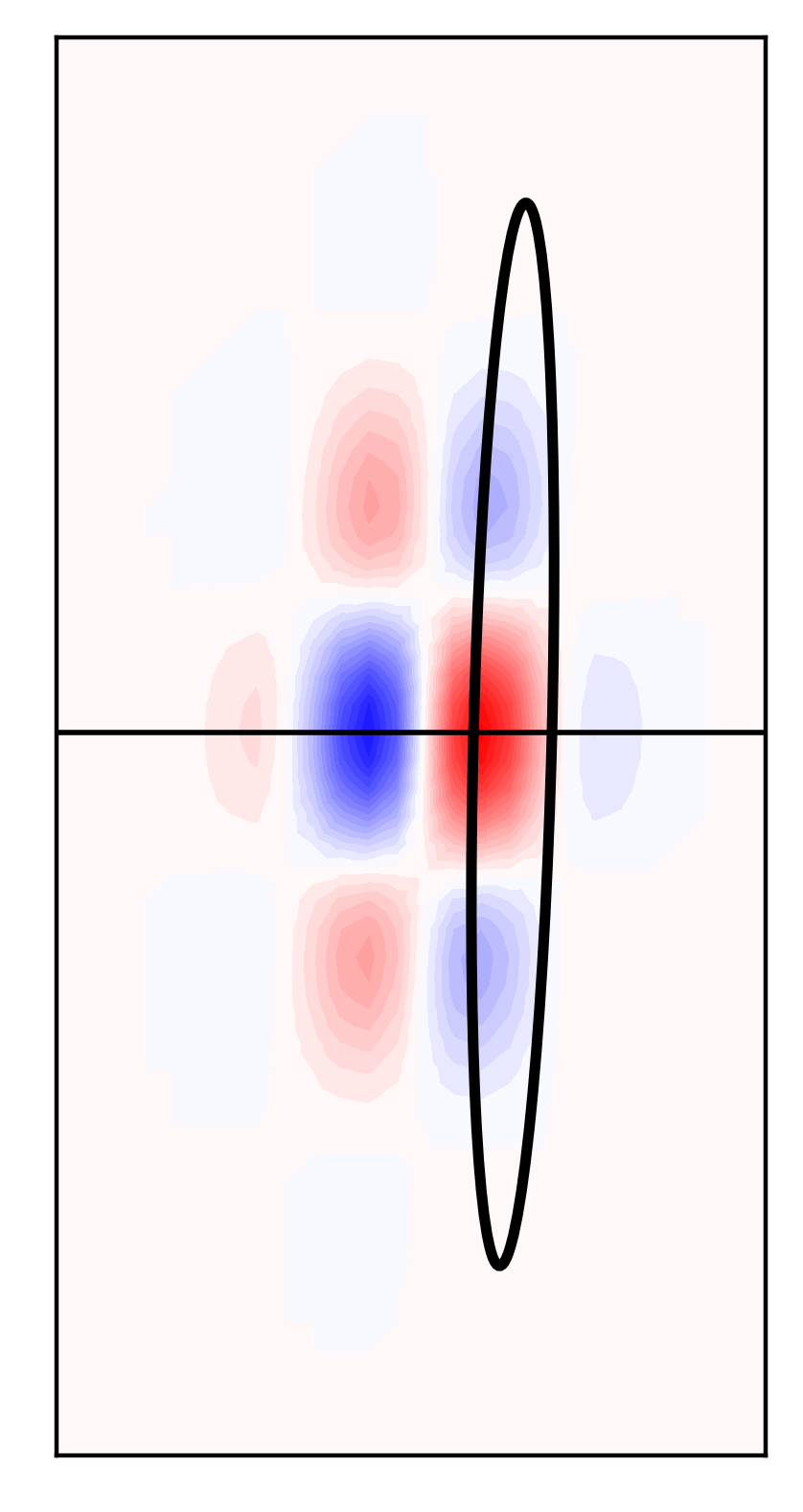} & 
			\includegraphics[scale=0.25, valign=c]{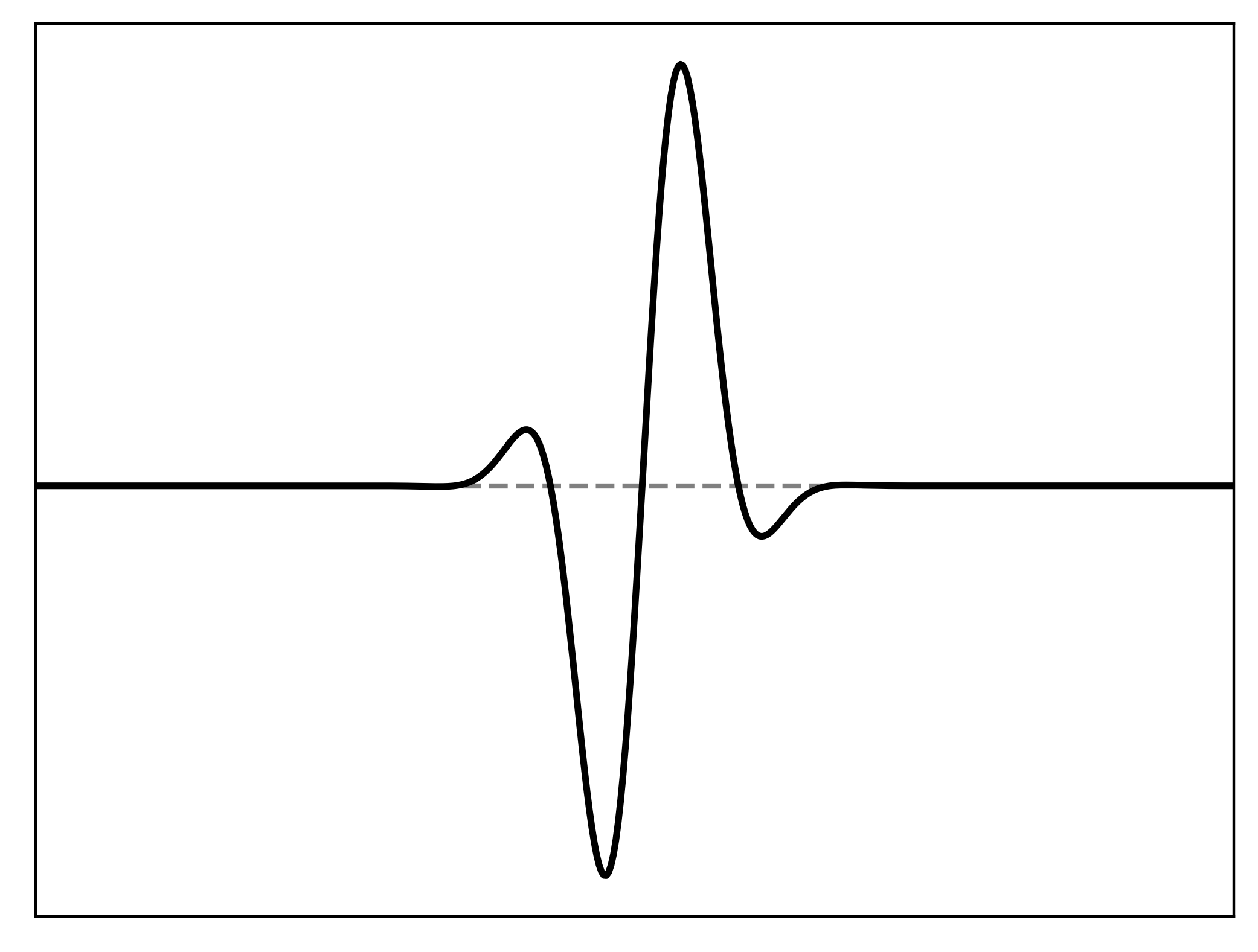} & 
			\includegraphics[scale=0.25, valign=c]{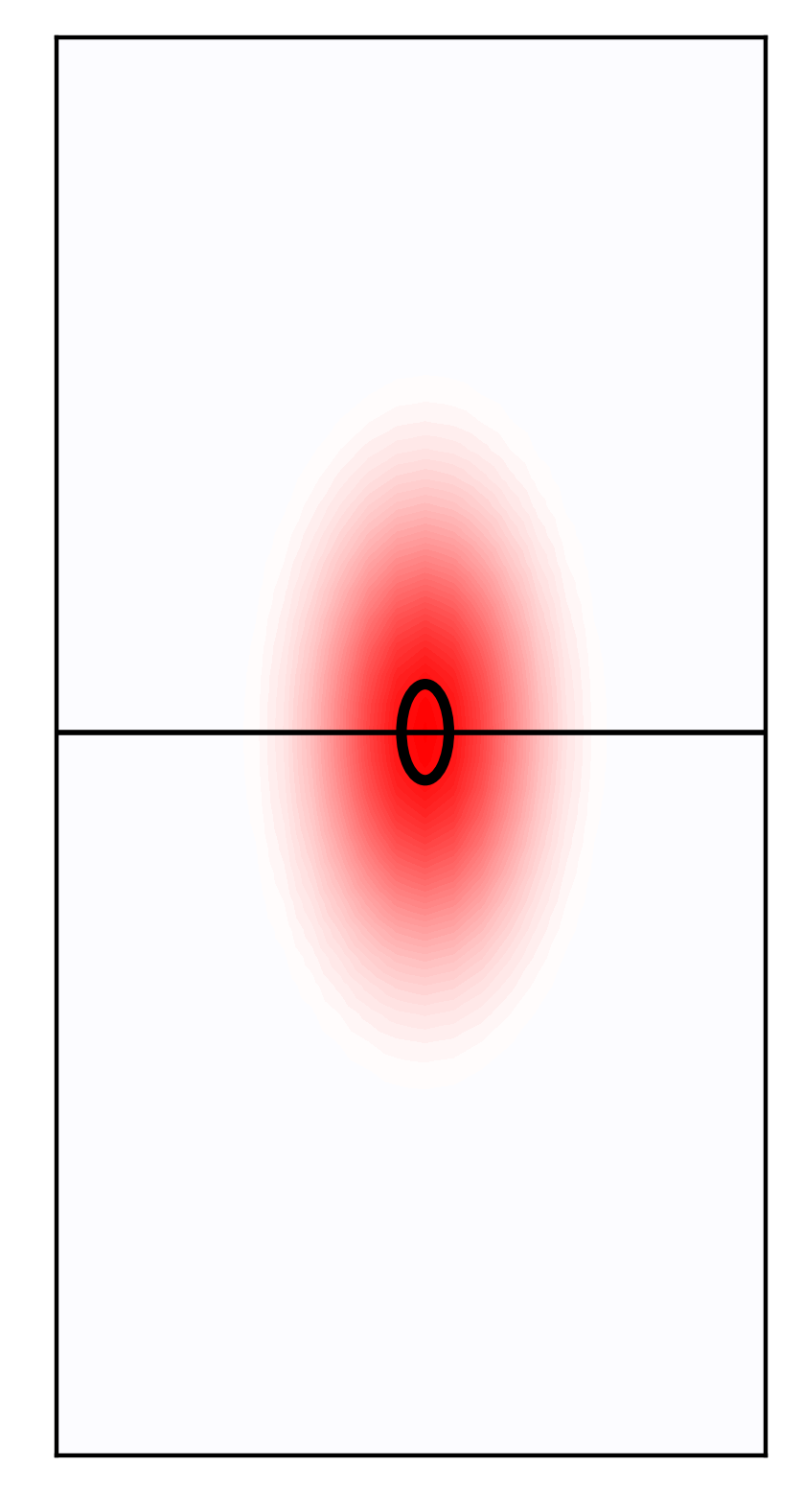} & 
			\includegraphics[scale=0.25, valign=c]{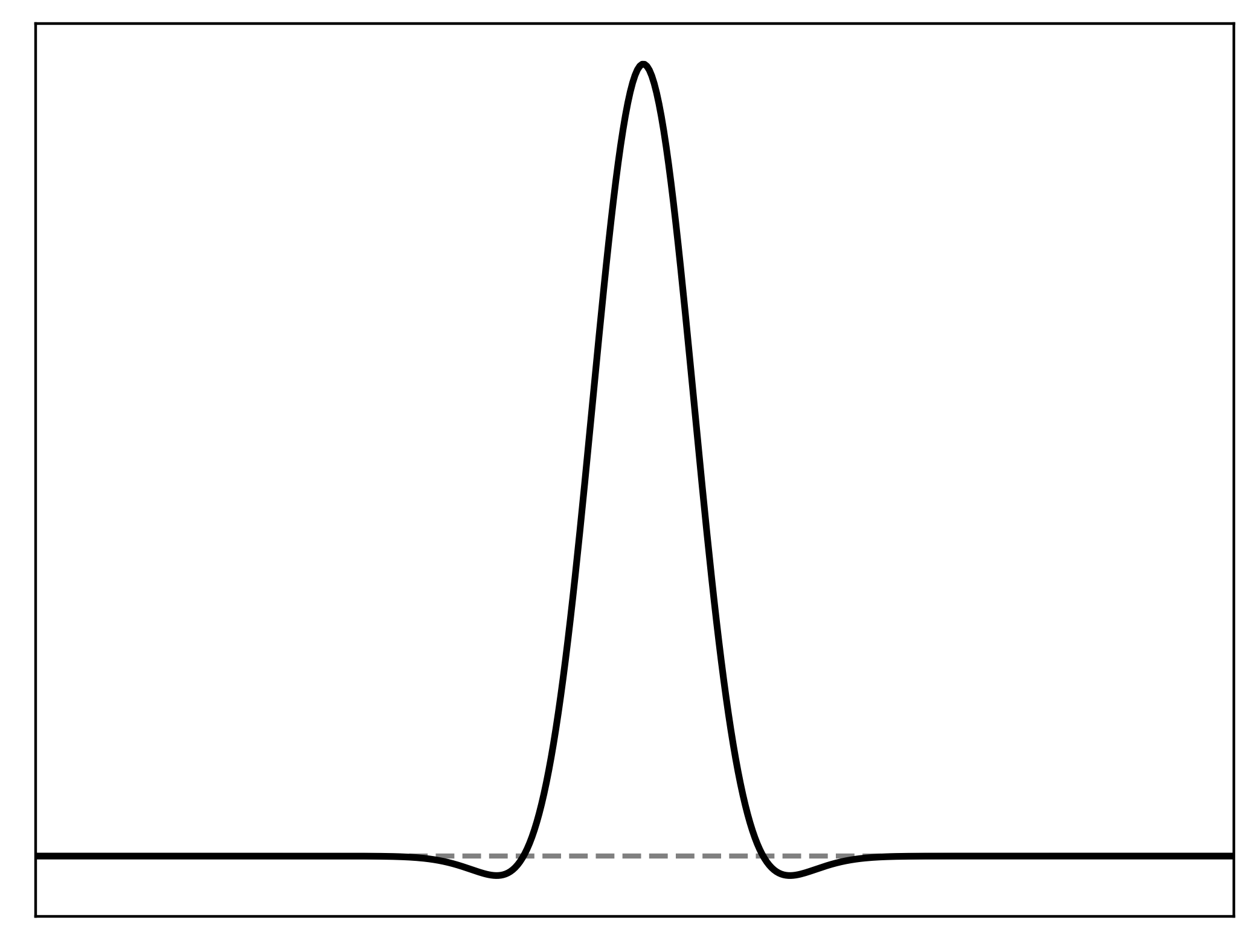} 
		\end{tabular}
	}
	\caption{Illustration of the influence of negative numbers in the integral kernel on the robustness of the ellipsoid estimates for the supports of impulse responses. Left two columns: Blur kernel given in Equation \ref{eq:frog_kernel}. Right two columns: Ricker wavelet-type kernel given by $\Aker(y,x) = (1 - a \gamma)\exp\left(-\gamma/2\right)$, where $\gamma=(y-x)^T \Sigma^{-1} (y-x)$, and $\Sigma=\operatorname{diag}\left(0.0025, 0.01\right)$. Ordered from top to bottom, the results are obtained with $a\in\{1.0, 20.0, 27.0\}$ for the left two columns, and $a\in \{0.0, 0.23, 0.249\}$ for the right two columns. Columns 1 and 3: impulse responses with estimated support ellipsoids indicated by the black ellipses. Red and blue represent positive and negative numbers in the integral kernel, respectively. Columns 2 and 4: one-dimensional slice along the horizontal line indicated in the two-dimensional plots. The dashed gray line is at zero.
	}
	\label{fig:robustness_to_negatives}
\end{figure}

If there are negative numbers in the integral kernel (violating Assumption 3), the ellipsoid estimation procedure may incur errors or fail, leading to poor performance or failure of the PSF-based method. In Figure \ref{fig:robustness_to_negatives} we investigate the robustness of the ellipsoid support estimation procedure to violations of Assumption 3. We study two integral kernel examples, both of which are parameterized by a quantity that controls how negative the kernels are. We make the following observations:
\begin{itemize}
	\item The larger and more numerous the negative numbers are, the more inaccurate the ellipsoid support estimate is.
	\item The further away from the center of the ellipsoid the negative numbers are, the more influence they have on the quality of the ellipsoid support estimate. This is because moment formulas (Equations \ref{eq:define_mean} and \ref{eq:define_cov}) assign more weight to entries in the kernel that are further from the center.
	\item Negative numbers affect the ellipsoid estimation method more if they are isolated, and less if they are balanced by nearby positive numbers.
	\item As kernels become more negative, the ellipsoid estimation performs well up to a certain threshold that depends on the spatial distribution of negative and positive entries. After that threshold is crossed, the estimation rapidly transitions to performing poorly and ultimately failing.
\end{itemize} 
For the kernel in the left two columns of Figure \ref{fig:robustness_to_negatives}, negative numbers are interspersed with positive numbers, allowing us to include a large amount of negative numbers before the ellipsoid estimation fails. For the kernel in the right two columns, the ellipsoid estimate fails with tiny amounts of negative numbers because the negative numbers are far away from the mean and not balanced by positive numbers at similar distances and angles. In the bottom two rows, we see the aforementioned threshold effect, in which the ellipsoid estimation method rapidly transitions from performing reasonably well to performing poorly with only a small change to the integral kernel.

\subsection{Finite element discretization}
\label{sec:finite_element_kernel}

In computations, functions are discretized and replaced by finite-dimensional vectors, and operators mapping between infinite-dimensional spaces are replaced by operators mapping between finite-dimensional spaces. In this paper we discretize the functions that $\mathcal{A}$ and $\mathcal{A}^T$ are applied to using continuous finite elements satisfying the Kronecker property (defined below). With minor modifications, the PSF-based method could be used with more general finite element methods, or other discretization schemes such as finite differences or finite volumes. These restrictions on discretization only apply to functions $u$ that $\mathcal{A}$ and $\mathcal{A}^T$ are applied to. Other functions that arise internally during the process of computing actions of $\mathcal{A}$, such as state variables in a PDE that is solved in a subproblem, may be discretized with any method.

Let $\febasis_1, \febasis_2, \dots, \febasis_\fedim$ be a set of continuous finite element basis functions used to discretize the problem on a mesh with mesh size parameter $h$, let $V_h := \Span\left(\febasis_1, \febasis_2, \dots, \febasis_\fedim\right)$ be the corresponding finite element space under the $L^2$ inner product, and let $\fepoint_i \in \mathbb{R}^\gdim$, $i=1,\dots, \fedim$ be the Lagrange nodes associated with the functions $\febasis_i$. We assume that the finite element basis satisfies the Kronecker property, i.e., $\febasis_i(\fepoint_i)=1$ and $\febasis_i(\fepoint_j)=0$ if $i\neq j$. 
For $u_h \in V_h$ we write $\mathbf{u} \in \mathbb{R}^\nsamplepts_\massmatrix$ to denote the coefficient vector for $u_h$ with respect to the finite element basis, i.e.,
$
u_h(x) = \sum_{i=1}^\fedim \mathbf{u}_i \febasis_i(x).
$
Linear functionals $\rho_h \in V_h'$ have coefficient dual vectors $\boldsymbol{\rho}\in \mathbb{R}^\nsamplepts_{\massmatrix^{-1}}$, with entries $\boldsymbol{\rho}_i = \rho_h(\febasis_i)$ for $i=1,\dots,\nsamplepts$.
Here, $\massmatrix \in \mathbb{R}^{\fedim \times \fedim}$ denotes the sparse finite element mass matrix which has entries $\massmatrix_{ij}=\int_\Omega \febasis_i(x) \febasis_j(x) \mathrm{d}x$ for $i,j=1,\dots,\fedim$.
The space $\mathbb{R}^\fedim_\massmatrix$ is $\mathbb{R}^\fedim$ with the inner product $(\mathbf{u},\mathbf{w})_\massmatrix := \mathbf{u}^T \massmatrix \mathbf{w}$, and $\mathbb{R}^\fedim_{\massmatrix^{-1}}$ is the analogous space with $\massmatrix^{-1}$ replacing $\massmatrix$. Direct calculation shows that $\mathbb{R}^\fedim_\massmatrix$ and $\mathbb{R}^\fedim_{\massmatrix^{-1}}$ are isomorphic to $V_h$ and $V_h'$ as Hilbert spaces, respectively.

After discretization, the operator $\Aop:L^2(\Omega) \rightarrow L^2(\Omega)'$ is replaced by an operator $A_h:V_h \rightarrow V_h'$, which becomes an operator
$
\mathbf{A}:\mathbb{R}^\fedim_\massmatrix \rightarrow \mathbb{R}^\fedim_{\massmatrix^{-1}}
$
under the isomorphism discussed above. The PSF-based method is agnostic to the computational procedure for approximating $\Aop$ with $\mathbf{A}$. What is important is that we do not have direct access to matrix entries $\mathbf{A}_{ij}$. Rather, we have a computational procedure that allows us to compute matrix-vector products $\mathbf{u}\mapsto \mathbf{A}\mathbf{u}$ and $\mathbf{w}\mapsto \mathbf{A}^T\mathbf{w}$, and computing these matrix-vector products is costly. The PSF-based method mitigates this cost by performing as few matrix-vector products as possible.
Of course, matrix entries can be computed via matrix-vector products as $\mathbf{A}_{ij} = \left(\mathbf{A}\mathbf{e}_j\right)_i$,
where $\mathbf{e}_j=(0,\dots,0,1,0,\dots,0)^T$ is the length $\fedim$ unit vector with one in the $j$th coordinate and zeros elsewhere. But computing the matrix-vector product $\mathbf{e}_j \mapsto \mathbf{A}\mathbf{e}_j$ is costly, and therefore wasteful if we do not use other matrix entries in the $j$th column of $\mathbf{A}$. Hence, methods for approximating $\mathbf{A}$ are computationally intractable if they require accessing scattered matrix entries from many different rows and columns of $\mathbf{A}$. 

The operator $A_h:V_h \rightarrow V_h'$ can be written in integral kernel form,~\eqref{eq:kernel_representation}, but with $\Aker$ replaced by a slightly different integral kernel,
$\Aker_h$, which we do not know, and which differs from $\Aker$ due to discretization error. Since the functions in $V_h$ are continuous at $x$, the delta distribution $\delta_x$ is a continuous linear functional on $V_h$, which has a discrete dual vector $\boldsymbol{\delta}_x \in \mathbb{R}^\fedim_{\massmatrix^{-1}}$ with entries $\left(\boldsymbol{\delta}_x\right)_i = \febasis_i(x)$ for $i=1,\dots,\fedim$. Additionally, it is straightforward to verify that the Riesz representation, $\rho_h^* \in V_h$, of a functional $\rho \in V_h'$ has coefficient vector
$
\boldsymbol{\rho}^* = \massmatrix^{-1} \boldsymbol{\rho}.
$
Therefore, the formula for the impulse response from~\eqref{eq:impulse_response_delta_action} becomes 
$
\boldsymbol{\impulseresponse}_x = \left(A_h \delta_x^*\right)^* =  \massmatrix^{-1}\mathbf{A} \massmatrix^{-1} \boldsymbol{\delta}_x,
$
and the $(y,x)$ kernel entry of $\Aker_h$ may be written as
$
	\Aker_h(y,x) = \boldsymbol{\delta}_y^T \boldsymbol{\impulseresponse}_x = \boldsymbol{\delta}_y^T \massmatrix^{-1}\mathbf{A} \massmatrix^{-1} \boldsymbol{\delta}_x.
$
Now define $\mathbf{\Aker} \in \mathbb{R}^{\fedim \times \fedim}$ to be the following dense matrix of kernel entries evaluated at all pairs of Lagrange nodes:
\begin{equation}
\label{eq:Akerpcmat_entries}
\mathbf{\Aker}_{ij} := \Aker_h(\fepoint_i, \fepoint_j).
\end{equation}
Because of the Kronecker property of the finite element basis, we have $\boldsymbol{\delta}_{\fepoint_i} = \mathbf{e}_i$. Thus, we have
$
	\Aker_h(\fepoint_i,\fepoint_j) = \left(\massmatrix^{-1}\mathbf{A} \massmatrix^{-1}\right)_{ij},
$
which implies
\begin{equation}
\label{eq:boldA}
	\mathbf{A} = \massmatrix \mathbf{\Aker} \massmatrix.
\end{equation}
Broadly, we will construct an H-matrix approximation of $\mathbf{A}$ by forming an H-matrix approximation of $\boldsymbol{\Aker}$, then multiplying $\boldsymbol{\Aker}$ by $\massmatrix$ (or a lumped mass version of $\massmatrix$) on the left and right using H-matrix methods. Classical H-matrix construction methods require access to arbitrary matrix entries $\mathbf{\Aker}_{ij}$, but these matrix entries are not easily accessible. The bulk of the PSF-based method is therefore dedicated to forming approximations of these matrix entries that can be evaluated rapidly.

\paragraph{Lumped mass matrix} At the continuum level, $\Aker$ is assumed to be non-negative. However, entries of $\boldsymbol{\Aker}$ involve inverse mass matrices, which typically contain negative numbers. 
We therefore recommend replacing the mass matrix, $\massmatrix$, with a positive diagonal \emph{lumped mass} approximation. 
Here, we use the lumped mass matrix in which the $i$th diagonal entry of the lumped mass matrix is the sum of all entries in the $i$th row of the mass matrix. Other mass lumping techniques may be used.

\section{Key innovations}
\label{sec:prerequisites}

In this section we present two key innovations that the PSF-based method is based on.
First, we define moments of the impulse responses, $\phi_x$, show how these moments can be computed efficiently, and use these moments to form ellipsoid shaped a-priori estimates for the supports of the impulse responses (Section~\ref{sec:intromoments}). 
Second, we describe an improved method to approximate impulse responses from other nearby impulse responses, which we call ``normalized local mean displacement invariance'' (Section~\ref{sec:local_mean_displacement_invariance}).

\subsection{Impulse response moments and ellipsoid support estimate}
\label{sec:intromoments}

The impulse response $\impulseresponse_x$ may be interpreted as a scaled probability distribution because of the non-negative integral kernel property. Let $\spatialvol:\Omega \rightarrow \mathbb{R}$,
\begin{equation}
\label{eq:define_vol}
\spatialvol(x) := \int_\Omega \impulseresponse_x(y) \mathrm{d}y,
\end{equation}
denote the spatially varying scaling factor, and for $i,j=1,\dots,\gdim$ define $\spatialmean:\Omega \rightarrow \mathbb{R}^\gdim$ and $\spatialcov:\Omega \rightarrow \mathbb{R}^{\gdim \times \gdim}$ as follows:
\begin{align}
\spatialmean^i(x) :=& \frac{1}{V(x)}\int_\Omega \impulseresponse_x(y) y^i ~\mathrm{d}y \label{eq:define_mean} \\
\spatialcov^{ij}(x) :=& \frac{1}{V(x)}\int_\Omega \impulseresponse_x(y) \left(y^i - \spatialmean^i(x)\right) \left(y^j - \spatialmean^j(x)\right) ~\mathrm{d}y, \label{eq:define_cov}
\end{align}
where $\spatialmean^i(x)$ and $y^i$ denote the $i^\text{th}$ components of the vectors $\spatialmean(x)$ and $y$, respectively, and $\spatialcov^{ij}(x)$ denotes the $(i,j)$ entry of the matrix $\spatialcov(x)$.
The quantities $\spatialmean(x)\in \mathbb{R}^\gdim$ and $\spatialcov(x) \in \mathbb{R}^{\gdim \times \gdim}$ are the mean and covariance of the normalized version of $\impulseresponse_x$, respectively. 

The direct approach to compute $\spatialvol(x)$, $\spatialmean(x)$, and $\spatialcov(x)$ is to apply $\Aop$ to a point source centered at $x$ to obtain $\impulseresponse_x$, per~\eqref{eq:impulse_response_delta_action}. Then one can post process $\impulseresponse_x$ to determine $\spatialvol(x)$, $\spatialmean(x)$, and $\spatialcov(x)$. However, this direct approach is not feasible because our algorithm for picking sample points (Section \ref{sec:sample_point_selection} and Figure \ref{fig:frog_batches}) needs to know $V(x)$, $\spatialmean(x)$, and $\spatialcov(x)$ before we compute $\phi_x$. Computing $\phi_x$ in order to determine $V(x)$, $\spatialmean(x)$, and $\spatialcov(x)$ would be extremely computationally expensive, and defeat the purpose of the PSF-based method, which is to reduce the computational cost by computing impulse responses in batches. Fortunately, it is possible to compute $\spatialvol(x)$, $\spatialmean(x)$, and $\spatialcov(x)$ indirectly, \emph{for all points $x\in\Omega$ simultaneously}, by applying $\Aop^T$ to one constant function, $\gdim$ linear functions, and $\gdim(\gdim+1)/2$ quadratic functions (e.g., 6 total operator applications in two spatial dimensions and 10 in three spatial dimensions). This may be motivated by analogy to matrices. If $\mathbf{A}\in \mathbb{R}^{\fedim \times \fedim}$ is a matrix with $i^\text{th}$ column $\mathbf{a}_i$ and $\mathbf{w} \in \mathbb{R}^\fedim$, then
\begin{equation*}
\mathbf{A}^T \mathbf{w} = \begin{bmatrix}
\horzbar & \mathbf{a}_1^T & \horzbar \\
& \vdots & \\
\horzbar & \mathbf{a}_\fedim^T & \horzbar
\end{bmatrix}
\mathbf{w} = 
\begin{bmatrix}
\mathbf{a}_1^T \mathbf{w} \\
\vdots \\
\mathbf{a}_\fedim^T \mathbf{w}
\end{bmatrix}.
\end{equation*}
By computing one matrix-vector product of $\mathbf{A}^T$ with $\mathbf{w}$, we compute the inner product of each column of $\mathbf{A}$ with $\mathbf{w}$ simultaneously. The operator case is analogous, with $\phi_x$ taking the place of a matrix column. We have
\begin{equation}
\label{eq:operator_simultaneous_ips}
\left(\Aop^T w\right)^*(x) = \int_\Omega \Aker(y,x) w(y) dy = \left(\phi_x, w\right)_{L^2(\Omega)}.
\end{equation}
By computing one operator application of $\Aop^T$ to $w$, we compute the inner product of each $\phi_x$ with $w$, for all points $x$ simultaneously. 

Let $C$, $L^i$, and $Q^{ij}$ be the following constant, linear, and quadratic functions:
\begin{equation*}
C(x) := 1, \qquad
L^i(x) := x^i, \qquad
Q^{ij}(x) := x^i x^j
\end{equation*}
for $i,j=1,\dots,\gdim$. Using the definition of $\spatialvol$ in~\eqref{eq:define_vol} and using~\eqref{eq:operator_simultaneous_ips}, we have
\begin{equation*}
	\spatialvol(x) = \int_\Omega \phi_x(y) C(y) ~ dy = \left(\phi_x, C\right)_{L^2(\Omega)} = \left(\mathcal{A}^T C\right)^*(x).
\end{equation*}
Hence, we compute $\spatialvol(x)$ for all $x$ simultaneously by applying $\Aop^T$ to $C$. Analogous manipulations show that $\spatialmean(x)$ and $\spatialcov(x)$ may be computed for all points $x$ simultaneously by applying $\Aop^T$ to the functions $L^i$ and $Q^{ij}$, respectively. We have
\begin{subequations}
	\label{eq:vol_mean_var}
	\begin{align}
	\spatialvol =& \left(\Aop^T C\right)^* \\
	\spatialmean^i =& \left(\Aop^T L^i\right)^* / \spatialvol \\
	\spatialcov^{ij} =& \left(\Aop^T Q^{ij}\right)^* / \spatialvol - \spatialmean^i\cdot \spatialmean^j
	\end{align}
\end{subequations}
for $i,j=1,\dots, \gdim$. Here, $u/w$ denotes pointwise division, $\left(u/w\right)(x) = u(x)/w(x)$, and $u\cdot w$ denotes pointwise multiplication, $(u \cdot w)(x) = u(x)w(x)$.

We approximate the support of $\impulseresponse_x$ with the ellipsoid
\begin{equation}
\label{eq:support_ellipsoid}
E_x := \{x' \in \Omega: (x' - \spatialmean(x))^T \spatialcov(x)^{-1} (x' - \spatialmean(x)) \le \tau^2\},
\end{equation}
where $\tau$ is a fixed constant (see Figure \ref{fig:frog_moments_ellipsoid}). The ellipsoid $E_x$ is the set of points within $\tau$ standard deviations of the mean of the Gaussian distribution with mean $\spatialmean(x)$ and covariance $\spatialcov(x)$, i.e., the Gaussian distribution which has the same mean and covariance as the normalized version of $\impulseresponse_x$. The quantity $\tau$ is a parameter that must be chosen appropriately. The larger $\tau$ is, the larger the ellipsoid $E_x$ is, and the more conservative the estimate is for the support of $\impulseresponse_x$. However, in Section~\ref{sec:sample_point_selection} we will see that the cost of the PSF-based method depends on how many non-overlapping ellipsoids $E_x$ we can ``pack'' in the domain $\Omega$ (more ellipsoids is better), and choosing a larger value of $\tau$ means that fewer ellipsoids will fit in $\Omega$. In practice, we find that $\tau=3.0$ yields a reasonable balance between these competing interests, and use $\tau=3.0$ in all numerical results, except for Figure \ref{fig:PSF_tau_convergence}, where we study the effects of varying $\tau$. The fraction of the ``mass'' of $\impulseresponse_x$ residing outside of $E_x$ is less than $1/\tau^2$ by Chebyshev's inequality, though this bound is typically conservative.

\subsection{Local mean displacement invariance}
\label{sec:local_mean_displacement_invariance}

Let $x$ and $x'$ be points in $\Omega$ that are close to each other, and consider the following approximations:
\begin{align}
\phi_x(y) &\approx \phi_{x'}(y) \label{eq:translate1}\\
\phi_x(y) &\approx \phi_{x'}(y-x+x') \label{eq:translate2}\\
\phi_x(y) &\approx \phi_{x'}\left(y-\spatialmean(x)+\spatialmean(x')\right) \label{eq:translate3}\\
\phi_x(y) &\approx \phi_{x'}\left(y-\spatialmean(x)+\spatialmean(x')\right) \spatialvol(x) / \spatialvol(x'). \label{eq:translate4}
\end{align}
These are four different ways to approximate an impulse response by a nearby impulse response, with each successive approximation building upon the previous ones. The PSF-based method uses~\eqref{eq:translate4}, which is the most sophisticated. Approximation~\eqref{eq:translate1} says that $\phi_x$ can be approximated by $\phi_{x'}$ when $x$ and $x'$ are close. Operators satisfying~\eqref{eq:translate1} can be well approximated via low rank CUR approximation.
However, the required rank in the low rank approximation can be large, which makes algorithms based on~\eqref{eq:translate1} expensive. Operators that satisfy~\eqref{eq:translate2} are called ``locally translation invariant'' because integral kernel entries $\Aker(y,x)$ for such operators are approximately invariant under translation of $x$ and $y$ by the same displacement, i.e., $x \rightarrow x+h$ and $y \rightarrow y+h$. It is straightforward to show that if equality holds in~\eqref{eq:translate2}, then $\Aop$ is a convolution operator. Locally translation invariant operators act like convolutions locally, and can therefore be well approximated by PC approximations.

Approximation~\eqref{eq:translate3} improves upon~\eqref{eq:translate1} and~\eqref{eq:translate2}, and generalizes both. On one hand, if~\eqref{eq:translate1} holds, then $\spatialmean(x) \approx \spatialmean(x')$, and so~\eqref{eq:translate3} holds. On the other hand, translating a distribution translates its mean, so if~\eqref{eq:translate2} holds, then $\spatialmean(x')-\spatialmean(x) \approx x' - x$, so again~\eqref{eq:translate3} holds. But approximation~\eqref{eq:translate3} can hold in situations where neither~\eqref{eq:translate1} nor~\eqref{eq:translate2} holds. For example, because the expected value commutes with affine transformations,~\eqref{eq:translate3} will hold when $\Aop$ is locally translation invariant with respect to a translated and rotated frame of reference, while~\eqref{eq:translate2} will not. Additionally,~\eqref{eq:translate3} generalizes to operators $\Aop:L^2(\Omega_1) \rightarrow L^2(\Omega_2)'$ that map between function spaces on different domains $\Omega_1$ and $\Omega_2$, and even operators that map between domains with different spatial dimensions. In contrast,~\eqref{eq:translate2} does not naturally generalize to operators that map between function spaces on different domains, because the formula $y-x+x'$ requires vectors in $\Omega_2$ and $\Omega_1$ to be added together.  We call~\eqref{eq:translate3} ``local mean displacement invariance,'' and illustrate~\eqref{eq:translate3} in Figure~\ref{fig:mean_displacement_invariance}.

We use approximation~\eqref{eq:translate4}, which is the same as~\eqref{eq:translate3}, except for the factor $\spatialvol(x)/\spatialvol(x')$. This factor makes the approximation more accurate if $\spatialvol(x)$ varies widely. Approximation~\eqref{eq:translate4} is equivalent to~\eqref{eq:translate3}, but with $\phi_x$ replaced by its normalized version, $\phi_x/\spatialvol(x)$. We call~\eqref{eq:translate4} \emph{normalized local mean displacement invariance}.

\section{Operator approximation algorithm}
\label{sec:method}

Before presenting the technical details of the algorithm in Sections \ref{sec:sample_point_selection}--\ref{sec:flip_eigs}, we first provide an overview.

We use~\eqref{eq:vol_mean_var} to compute $\spatialvol$, $\spatialmean$, and $\spatialcov$ by applying $\Aop^T$ to polynomial functions. Then we use~\eqref{eq:support_ellipsoid} to form ellipsoid shaped estimates for the support of the $\phi_x$'s, \emph{without} computing them (see Figure \ref{fig:frog_moments_ellipsoid}). This allows us to compute large numbers of $\impulseresponse_{x_i}$ in ``batches,'' $\eta_b$ (see Figures~\ref{fig:batches_intro}~and~\ref{fig:frog_batches}). We compute one batch, denoted $\eta_b$, by applying $\Aop$ to a weighted sum of point sources (Dirac comb) associated with a batch, $S_b$, of points $x_i$ scattered throughout $\Omega$ (Section~\ref{sec:get_impulse_response}). The batch of points, $S_b$, is chosen via a greedy ellipsoid packing algorithm so that, for $x_i,x_j \in S_b$, the support ellipsoid for $\impulseresponse_{x_i}$ and the support ellipsoid for $\impulseresponse_{x_j}$ do not overlap if $i \neq j$ (Section~\ref{sec:sample_point_selection}). Because these supports do not overlap (or do not overlap much), we can post process $\eta_b$ to recover the functions $\impulseresponse_{x_i}$ associated with all points $x_i \in S_b$. With one application of $\Aop$, we recover many $\impulseresponse_{x_i}$ (Section~\ref{sec:get_impulse_response}). The process is repeated until a desired number of batches is reached.

Once the batches $\eta_b$ are computed, we approximate the integral kernel $\Aker(y,x)$ at arbitrary points $(y,x)$ by interpolation of translated and scaled versions of the computed $\impulseresponse_{x_i}$ (Section~\ref{sec:approximate_kernel_entries} and Figure \ref{fig:hmatrix_neighbors_frog}). The key idea behind the interpolation is the normalized local mean displacement invariance assumption discussed in Section~\ref{sec:local_mean_displacement_invariance}. Specifically, we approximate $\Aker(y,x) = \phi_x(y)$ by a weighted linear combination of the values $\frac{\spatialvol(x)}{\spatialvol(x_i)}\phi_{x_i}(y - \spatialmean(x) + \spatialmean(x_i))$ for a small number of sample points $x_i$ near $x$. The weights are determined by radial basis function (RBF) interpolation.

The ability to rapidly evaluate approximate kernel entries $\Aker(y,x)$ allows us to construct an H-matrix approximation, $\boldsymbol{\Aker}_H \approx \mathbf{\Aker}$, using the conventional adaptive cross H-matrix construction method (Section~\ref{sec:h_matrix_construction_short}). In this method, one forms low rank approximations of off-diagonal blocks of the matrix by sampling rows and columns of those blocks. We then convert $\boldsymbol{\Aker}_H$ into an H-matrix approximation $\mathbf{A}_H \approx \mathbf{A}$. 

When $\Aop$ is symmetric positive semi-definite, $\mathbf{A}_H$ may be non-symmetric and indefinite due to errors in the approximation. In this case, one may optionally symmetrize $\mathbf{A}_H$, then modify it via low rank updates to remove erroneous negative eigenvalues (Section \ref{sec:flip_eigs}). The complete algorithm for constructing $\mathbf{A}_H$ is shown in Algorithm~\ref{alg:construct_Atilde}. The computational cost is discussed in Section~\ref{sec:overall_cost}.

\begin{algorithm2e}
	\SetAlgoNoLine
	\SetKwInOut{Input}{Input}
	\SetKwInOut{Output}{Output}
	{	
		\Input{Linear operator $\Aop$, parameter $\nbatch$}
		\Output{H-matrix $\mathbf{A}_H$}
		
		Compute $V, \mu$, and $\Sigma$ (Equations~\eqref{eq:vol_mean_var} in Section~\ref{sec:intromoments})
		
		\For{$k=1,2,\dots,\nbatch$}{
			Choose a batch of sample points, $\pointbatch_k$ (Section~\ref{sec:sample_point_selection})
			
			Compute impulse response batch $\combresponse_k$ by applying $\Aop$ to the Dirac comb for $\pointbatch_k$ (Section~\ref{sec:get_impulse_response})
			
		}

		Form H-matrix approximation $\boldsymbol{\Aker}_H$ of integral kernel (Sections~\ref{sec:approximate_kernel_entries} and \ref{sec:h_matrix_construction_short})

		Form H-matrix approximation $\mathbf{A}_H$ of $\Aop$ (Section~\ref{sec:h_matrix_construction_short})
		
		(optional) Modify $\mathbf{A}_H$ to make it symmetric and remove negative eigenvalues (Section \ref{sec:flip_eigs})
		
	}
	\caption{Construct PSF H-matrix approximation}
	\label{alg:construct_Atilde}
\end{algorithm2e}

\subsection{Sample point selection via greedy ellipsoid packing}
\label{sec:sample_point_selection}

We choose sample points, $x_i$, in batches $\pointbatch_k$. We use a greedy ellipsoid packing algorithm to choose as many points as possible per batch, while ensuring that there is no overlap between the support ellipsoids, $E_{x_i}$, associated with the sample points within a batch.

We start with a finite set of candidate points $\candidatepts$ and build $\pointbatch_k$ incrementally with points selected from $\candidatepts$. For simplicity of explanation, here $\pointbatch_k$ and $\candidatepts$ are mutable sets that we add points to and remove points from. First we initialize $\pointbatch_k$ as an empty set. Then we select the candidate point $\candidatepoint_i \in \candidatepts$ that is the farthest away from all points in previous sample point batches $S_1 \cup S_2 \cup \dots \cup S_{k-1}$. Candidate points for the first batch $S_1$ are chosen randomly from $\candidatepts$.
Once $\candidatepoint_i$ is selected, we remove $\candidatepoint_i$ from $\candidatepts$. Then we perform the following checks:
\begin{enumerate}
	\item We check whether $\candidatepoint_i$ is sufficiently far from all of the previously chosen points in the current batch, in the sense that $E_{\candidatepoint_i} \cap E_{\candidatepoint_j} = \{\}$ for all $\candidatepoint_j \in \pointbatch_k$.
	\item We make sure that $\spatialvol(\candidatepoint_i)$ is not too small, by checking whether $\spatialvol(\candidatepoint_i) > \epsilon_V \spatialvol_\text{max}$. Here, $\spatialvol_\text{max}$ is the largest value of $\spatialvol(\candidatepoint_j)$ over all points $q$ in the initial set of candidate points, and $\epsilon_V$ is a small threshold parameter (we use $\epsilon_V=10^{-5}$).
	\item We make sure that all eigenvalues of $\spatialcov(\candidatepoint_i)$ are positive, and the aspect ratio of $E_{\candidatepoint_i}$ (square root of the ratio of the largest eigenvalue of $\spatialcov(\candidatepoint_i)$ to the smallest) is bounded by a constant $1/\epsilon_\Sigma$ (we use $1/\epsilon_\Sigma=20$). Negative integral kernel entries due to discretization error can cause $\spatialcov(\candidatepoint_i)$ to be indefinite or highly ill-conditioned.
\end{enumerate}
If $\candidatepoint_i$ passes these checks 
then we add $\candidatepoint_i$ to $\pointbatch_k$. Otherwise we discard $\candidatepoint_i$. This process repeats until there are no more points in $\candidatepts$.  
We repeat the point selection process to construct several batches of points $\pointbatch_1, \pointbatch_2, \dots, \pointbatch_{\nbatch}$. For each batch, $\candidatepts$ is initialized as the set of all Lagrange nodes for the finite element basis functions used to discretize the problem, except for points in previous batches.

We check whether $E_{\candidatepoint_i} \cap E_{\candidatepoint_j} = \{\}$ in a two stage process. First, we check whether the axis aligned bounding boxes for the ellipsoids intersect. This quickly rules out intersections of ellipsoids that are far apart. Second, if the bounding boxes intersect, we check if the ellipsoids intersect using the ellipsoid intersection test in \cite{GilitschenskiHanebeck12}.

\subsection{Impulse response batches}
\label{sec:get_impulse_response}

We compute impulse responses, $\impulseresponse_{x_i}$, in batches by applying $\Aop$ to Dirac combs. The Dirac comb, $\diraccomb_k$, associated with a batch of sample points, $\pointbatch_k$, is the following weighted sum of Dirac distributions (point sources) centered at the points $x_i \in \pointbatch_k$:
\begin{equation*}
	\diraccomb_k := \sum_{x_i \in \pointbatch_k} \delta_{x_i} / \spatialvol(x_i).
\end{equation*}
We compute the \emph{impulse response batch}, $\eta_k$, by applying $\Aop$ to the Dirac comb:
\begin{equation}
	\label{eq:dirac_comb_H_action}
	\combresponse_k := \left(\Aop \diraccomb_k^*\right)^*
	=\sum_{x_i \in \pointbatch_k} \impulseresponse_{x_i} / \spatialvol(x_i).
\end{equation}
The last equality in~\eqref{eq:dirac_comb_H_action} follows from linearity and the definition of $\impulseresponse_{x_i}$ in~\eqref{eq:impulse_response_delta_action}. Since the points $x_i$ are chosen so that the ellipsoid $E_{x_i}$ that (approximately) supports $\impulseresponse_i$, and the ellipsoid $E_{x_j}$ that (approximately) supports $\impulseresponse_j$ do not overlap when $i \neq j$, we have (approximately)
\begin{equation}
\label{eq:varphi_eval}
	\impulseresponse_{x_i}(z) =
	\begin{cases}
		\combresponse_k(z) \spatialvol(x_i), & z \in E_{x_i} \\
		0, & \text{otherwise}
	\end{cases}
\end{equation}
for all $x_i \in \pointbatch_k$. By applying the operator once, $\xi_k \mapsto \left(\Aop \diraccomb_k^*\right)^*$, we recover $\impulseresponse_{x_i}$ for every point $x_i \in \pointbatch_k$. 

Each point source, $\delta_{x_i}$, is scaled by $1/\spatialvol(x_i)$ so that the resulting scaled impulse responses within $\eta_k$ are comparable in magnitude. Without this scaling, the portion of $\phi_{x_i}$ outside of $E_{x_i}$, which we neglect, may overwhelm $\phi_{x_j}$ for a nearby point $x_j$ if $\spatialvol(x_i)$ is much larger than $\spatialvol(x_j)$. Note that we are not in danger of dividing by zero, because the ellipsoid packing procedure from Section~\ref{sec:sample_point_selection} excludes $x_i$ from consideration as a sample point if $\spatialvol(x_i)$ is smaller than a predetermined threshold.

\subsection{Approximate integral kernel entries}
\label{sec:approximate_kernel_entries}

Here, we describe how to rapidly evaluate arbitrary entries of an approximation to the integral kernel by performing radial basis function interpolation of translated and scaled versions of nearby known impulse responses. In Section~\ref{sec:h_matrix_construction_short} we use this procedure for rapidly evaluating kernel entries to construct the H-matrix approximation of $\mathbf{A}$.

Given $(y,x)\in \Omega \times \Omega$, let
$
	z_i := y - \spatialmean(x) + \spatialmean(x_i)
$
and define
\begin{equation}
\label{eq:fxyxp}
f_i := \frac{\spatialvol(x)}{\spatialvol(x_i)}\phi_{x_i}\left(z_i\right)
\end{equation}
for $i=1,\dots,\numnbr$, where $\{x_i\}_{i=1}^{\numnbr}$ are the $\numnbr$ nearest sample points to $x$, excluding points $x_i$ for which $z_i \notin \Omega$. Here, $\numnbr$ is a small user-defined parameter, e.g., $\numnbr=10$. We find the $\numnbr$ nearest sample points to $x$ by querying a precomputed kd-tree~\cite{Bentley75} of all sample points.  We check whether $z_i \in \Omega$ by querying a precomputed axis aligned bounding box tree (AABB tree)~\cite{Ericson04} of the mesh cells used to discretize the problem. 
Note that $\phi_{x_i}\left(z_i\right)$ is well-defined because $z_i \in \Omega$, and $\frac{\spatialvol(x)}{\spatialvol(x_i)}$ is well-defined because the sample point choosing procedure in Section~\ref{sec:sample_point_selection} ensures that $\spatialvol(x_i)>0$. Per the discussion in Section~\ref{sec:local_mean_displacement_invariance}, we expect $\Aker(y,x) \approx f_i$ for $i=1,\dots,\numnbr$. The closer $x_i$ is to $x$, the better we expect the approximation to be. We therefore approximate $\Aker(y,x)$ by interpolating the (point,value) pairs
$
\left\{\left(x_i, f_i\right)\right\}_{i=1}^{\numnbr}
$
at the point $x$. 
Interpolation is performed using the following radial basis function~\cite{Wendland04} scheme:
\begin{equation}
\Aker(y,x) \approx \widetilde{\Aker}(y,x) := \sum_{i=1}^{\numnbr} \rbfweight_i~ \rbf\left(\|x-x_i\|\right),
\end{equation}
where $\rbfweight_i$ are weights, and 
$\rbf(r):= \exp\left(-\frac{1}{2}\left(C_\text{RBF}\frac{r}{r_0}\right)^2\right)$
is a Gaussian kernel radial basis function. Here, $r_0:=\diam\left(\{x_i\}_{i=1}^{\numnbr}\right)$ is the diameter of the set of sample points used in the interpolation, and $C_\text{RBF}$ is a user-defined shape parameter that controls the width of the kernel function. The vector of weights, $\rbfweight = (\rbfweight_1, \rbfweight_2, \dots, \rbfweight_{\numnbr})^T$, is found as the solution to the $\numnbr \times \numnbr$ linear system
\begin{equation}
B \rbfweight = f,
\end{equation}
where $B \in \mathbb{R}^{\numnbr \times \numnbr}$,
$
B_{ij} := \rbf\left(\|x_i - x_j\|\right),
$
and $f \in \mathbb{R}^{\numnbr}$ has entries $f_i$ from~\eqref{eq:fxyxp}. 

To evaluate $f_i$, we check whether $z_i \in E_{x_i}$ using~\eqref{eq:support_ellipsoid}. If $z_i \notin E_{x_i}$, then $z_i$ is outside the estimated support of $\impulseresponse_{x_i}$, so we set $f_i=0$. If $z_i \in E_{x_i}$, we look up the batch index $b$ such that $x_i \in S_b$, and evaluate $f_i$ via the formula
$
f_i = \spatialvol(x)\eta_b\left(z_i\right),
$
per~\eqref{eq:varphi_eval}. Note that $z_i$ is typically not a gridpoint of the mesh used to discretize the problem, even if $y$, $x$, and $x_i$ are gridpoints. Hence,
evaluating $\eta_b\left(z_i\right)$ requires determining which mesh cell contains $z_i$, then evaluating finite element basis functions on that mesh cell. Fortunately, the mesh cell containing $z_i$ was determined as a side effect of querying the AABB tree of mesh cells when we checked whether $z_i \in \Omega$. 

The shape parameter, $C_\text{RBF}$, mediates a tradeoff between accuracy and stability. Small $C_\text{RBF}$ is required for RBF interpolation with Gaussian kernels to achieve high accuracy, but small $C_\text{RBF}$ also makes RBF interpolation less robust to errors or nonsmoothness in the function being interpolated. 
For our numerical results involving Hessians in inverse problems governed by PDEs (Sections~\ref{sec:ice}~and~\ref{sec:adv}), high accuracy is not required because the PSF-based method is used to construct a preconditioner. Hence, for these Hessian approximations we use a conservative choice of $C_\text{RBF}=3.0$ to ensure robustness. For our numerical results for the blur problem example (Section~\ref{sec:frog}), we use a smaller value of $C_\text{RBF}=0.5$ so that the RBF interpolation accuracy is not a limiting factor as we study convergence of the PSF-based method.

\subsection{Hierarchical matrix construction}
\label{sec:h_matrix_construction_short}

We form an H-matrix approximation $\mathbf{A}_H \approx \mathbf{A}$ by forming an H-matrix representation $\mathbf{\Aker}_H$ of $\mathbf{\Aker}$ then multiplying $\mathbf{\Aker}$ with mass matrices $\massmatrix$ per~\eqref{eq:boldA} to form
$
\mathbf{A}_H = \massmatrix \mathbf{\Aker}_H \massmatrix.
$
Here, we use a diagonal lumped mass matrix, so these matrix-matrix multiplications are trivial. If a non-diagonal mass matrix is used, one may form an H-matrix representation of the mass matrix, then perform the matrix-matrix multiplications in \eqref{eq:boldA} using H-matrix methods.
We use H1 matrices in the numerical results, but any other H-matrix format could be used instead. For more details on H-matrices, see~\cite{Hackbusch15}. 

We form $\mathbf{\Aker}_H$ using the standard geometrical clustering/adaptive cross method implemented within the HLIBpro software package~\cite{Kriemann08}. For details about the algorithms used for geometrical clustering, H-matrix construction, and H-matrix operations in HLIBpro, we refer the reader to~\cite{BormGrasedyckHackbusch03,GrasedyckKriemannLeBorne08,Kriemann13}. Although $\mathbf{\Aker}$ is a dense $\fedim \times \fedim$ matrix, constructing $\mathbf{\Aker}_H$ only requires evaluation of $O(\hrank \fedim \log \fedim)$ kernel entries $\mathbf{\Aker}_{ij} = \widetilde{\Aker}(\fepoint_i, \fepoint_j)$ (see \cite{BebendorfRjasanow03}), and these entries are computed via the radial basis function interpolation method described in Section~\ref{sec:approximate_kernel_entries}. Here, $\hrank$ is the rank of the highest rank block in the H-matrix. 
We emphasize that the dense matrix $\mathbf{\Aker}$ is never formed.

\subsection{Symmetrizing and flipping negative eigenvalues (optional)}
\label{sec:flip_eigs}

In many applications, one seeks to approximate an operator
$\mathcal{H} = \Aop + \mathcal{R},$
where $\Aop$ is a symmetric positive semi-definite operator that we approximate with the PSF-based method to form an H-matrix $\mathbf{A}_H$, and $\mathcal{R}$ is a symmetric positive definite operator that may be easily converted to an H-matrix $\mathbf{R}_H$ without using the PSF-based method. For example, in inverse problems $\mathcal{H}$ is the Hessian, $\Aop$ is the data misfit term in the Hessian which is dense and available only matrix-free, and $\mathcal{R}$ is the regularization term, which is typically an elliptic differential operator that becomes a sparse matrix after discretization.

The PSF-based approximation $\mathbf{A}_H$, and therefore $\mathbf{A}_H + \mathbf{R}_H$, may be non-symmetric and indefinite because of approximation error. 
This is undesirable because symmetry and positive semi-definiteness are important properties which should be preserved if possible. Also, lacking these properties may prevent one from using highly effective algorithms to perform further operations involving $\mathbf{A}_H + \mathbf{R}_H$, such as using $\mathbf{A}_H + \mathbf{R}_H$ as a preconditioner in the conjugate gradient method.

We modify $\mathbf{A}_H$ to make it symmetric and remove negative eigenvalues via the following procedure. First, we symmetrize $\mathbf{A}_H$ via
$
\mathbf{A}_H^{\text{sym}} := \frac{1}{2}\left(\mathbf{A}_H + \mathbf{A}_H^T\right).
$
Next, we find negative eigenvalues and their corresponding eigenvectors for the generalized eigenvalue problem
$\mathbf{A}_H^{\text{sym}} \mathbf{u} = \lambda \mathbf{R}_H$
using a Cayley shift-and-invert Krylov scheme \cite{LehoucqEtAl98}. We flip the signs of these eigenvalues to be positive instead of negative (i.e., $\lambda \rightarrow |\lambda|$) by performing a low rank update to $\mathbf{A}_H^{\text{sym}}$. We observe that the eigenvectors associated with large erroneous negative eigenvalues tend to be directions that are nevertheless ``important'' to $\mathcal{A}$, so flipping the eigenvalues instead of setting them to zero tends to yield better approximations. The primary computational task in the Cayley shift-and-invert scheme is the solution of shifted linear systems of the form
$\left(\mathbf{A}_H^{\text{sym}} + \mu_i \mathbf{R}_H\right) \mathbf{x} = \mathbf{b},$
for a small number of positive shifts $\mu_i$. We solve these linear systems by factorizing the matrices $\mathbf{A}_H^{\text{sym}} + \mu_i \mathbf{R}_H$ using fast H-matrix methods. We compute and flip all eigenvalues $\lambda < \epsilon_\text{flip}$ which are less than some threshold $\epsilon_\text{flip} \in (-1,0]$. By choosing $\epsilon_\text{flip} > -1$, we ensure that the modified version of $\mathbf{A}_H^{\text{sym}} + \mathbf{R}_H$ is positive definite. Choosing $\epsilon_\text{flip}=0$ would remove all erroneous negative eigenvalues. However, this is computationally infeasible if $\Aop$ has a large or infinite cluster of eigenvalues near zero, a common situation for Hessians in ill-posed inverse problems. We therefore recommend choosing $\epsilon_\text{flip} < 0$. In our numerical results, we use $\epsilon_\text{flip}=-0.1$.

\section{Computational cost}
\label{sec:overall_cost}

The computational cost of the PSF-based method may be divided into the costs to perform the following tasks: (1) Computing impulse response moments and batches (Lines 1 and 4 in Algorithm \ref{alg:construct_Atilde}); (2) Building the H-matrix (Lines 5 and 6 in Algorithm \ref{alg:construct_Atilde}); (3) Performing linear algebra operations with the H-matrix. This may optionally include the symmetric positive semi-definite modifications described in Section \ref{sec:flip_eigs}. In target applications, (1) is the dominant cost because applying $\Aop$ to a vector requires an expensive computational procedure such as solving a PDE, and (1) is the only step that requires applying $\Aop$ to vectors.
All operations that do not require applications of $\Aop$ to vectors are polylog linear (i.e., $O(N \log(N)^b)$ for some $b$), and therefore scalable, in the size of the problem, $N$.
We now describe these costs in detail. For convenience, Table~\ref{tab:vars} lists variable symbols and their approximate sizes.

\begin{table}
	\begin{tabular}{lll}
		Symbol & Typical size & Variable name \\
		\hline
		$\fedim$ & $10^3$--$10^9$ & Number of finite element degrees of freedom \\
		$\nbatch$ & $1$--$25$ & Number of batches \\
		$\hrank$ & $5$--$50$ & H-matrix rank \\
		$\numnbr$ & $5$--$15$ & Number of nearest neighbors for RBF interpolation \\
		$\gdim$ & $1$--$3$ & Spatial dimension \\
		$\nsamplepts$ & $10^1$--$10^4$ & Total number of sample points (all batches) \\
		$|\pointbatch_i|$ & $1$--$500$ & Number of sample points in the $i$th batch
	\end{tabular}
	\caption{Symbols used for variables in computational cost estimates, and approximate ranges for their sizes in practice.}
	\label{tab:vars}
\end{table}

\paragraph{(1) Computing impulse response moments and batches} Computing $\spatialvol$, $\spatialmean$, and $\spatialcov$ requires applying $\Aop^T$ to $1$, $\gdim$, and $\gdim(\gdim+1)/2$ vectors, respectively. This works out to $3$ applications of $\Aop^T$ in one spatial dimension, $6$ in two dimensions, and $10$ in three dimensions. 
Computing each $\eta_i$ requires applying $\Aop$ to one vector, so computing $\{\eta_i\}_{i=1}^{\nbatch}$ requires $\nbatch$ operator applications. In total, computing all impulse response moments and batches therefore requires
\begin{equation*}
	1 + \gdim + \gdim(\gdim+1)/2 + \nbatch \qquad \text{operator applications.}
\end{equation*}
In a typical application one might have $\gdim=2$ and $\nbatch=5$, in which case a modest $11$ operator applications are required.

Computing the impulse response batches also requires choosing sample point batches via the greedy ellipsoid packing algorithm described in Section~\ref{sec:sample_point_selection}. Choosing the $i$th batch of sample points may require performing up to $\fedim |\pointbatch_i|$ ellipsoid intersection tests, where $|\pointbatch_i|$ is the number of sample points in the $i$th batch. Choosing all of the sample points therefore requires performing at most
\begin{equation*}
	\fedim \nsamplepts \qquad \text{ellipsoid intersection tests,}
\end{equation*}
where $\nsamplepts$ is the total number of sample points in all batches. The multiplicative dependence of $\fedim$ with $\nsamplepts$ is undesirable since $\nsamplepts$ may be large, and reducing this cost is possible with more involved computational geometry methods. 
However, from a practical perspective, the cost of choosing sample points is small compared to other parts of the algorithm, and hence such improvements are not pursued here.

\paragraph{(2) Building the H-matrix} 
Classical H-matrix construction techniques
require evaluating $O(\hrank \fedim \log \fedim)$ matrix entries of the approximation \cite{BebendorfRjasanow03}, where $\hrank$ is the H-matrix rank, i.e,  the maximum rank among the blocks of the H-matrix. To evaluate one matrix entry, first one must find the $\numnbr$ nearest sample points to a given point, where $\numnbr$ is the number of impulse responses used in the RBF interpolation. This is done using a precomputed kd-tree of sample points, and requires $O(\numnbr \log \nsamplepts)$ floating point and logical elementary operations. Second, one must find the mesh cells that the points $\{z_i\}_{i=1}^{\numnbr}$ reside in. This is done using an AABB tree of mesh cells, and requires $O(\numnbr \log \fedim)$ elementary operations. Third, one must evaluate finite element basis functions on those cells, which requires $O(\numnbr)$ elementary operations. Finally, the radial basis function interpolation requires solving a $\numnbr \times \numnbr$ linear system, which requires $O(\numnbr^3)$ elementary operations. Therefore, building the H-matrix requires
\begin{equation*}
	O\left(\left(\hrank \fedim \log \fedim\right)\left(\numnbr \log \fedim + \numnbr^3\right)\right) \qquad \text{elementary operations}.
\end{equation*}

\paragraph{(3) Performing linear algebra operations with the H-matrix} It is well known that H-matrix methods for matrix-vector products, matrix-matrix addition, matrix-matrix multiplication, matrix factorization, matrix inversion, and low rank updates require performing $O\left(\hrank^a \fedim \log(\fedim)^b\right)$
elementary operations,
where $a,b \in \{0,1,2,3\}$ are constants which depend on the type of H-matrix used and the operation being performed~\cite{GrasedyckHackbusch03}\cite[Section 2.1]{Kriemann13}.
For our numerical results involving Hessians (Sections~\ref{sec:ice} and~\ref{sec:adv}), we use one matrix-matrix addition to add the H-matrix approximation of the data misfit term in the Hessian to the regularization term in the Hessian.
Symmetrizing $\mathbf{A}_H$ requires one matrix-matrix addition. Flipping negative eigenvalues to be positive requires a handful (typically around 5) of matrix-matrix additions and matrix factorizations to factor the required shifted linear systems, and a number of factorized solves that is proportional to the number of erroneous negative eigenvalues.

In summary, computing all the necessary ingredients to evaluate kernel entries of the PSF-based approximation requires a handful of operator applications (e.g., $6+\nbatch$ operator applications in two dimensions, or $10+\nbatch$ operator applications in three dimensions, with $\nbatch$ typically in the range 1--25), plus comparatively cheap additional overhead costs, most notably performing ellipsoid intersection tests while choosing sample point batches. Once these ingredients are computed, no more operator applications (and thus PDE solves) are required, and approximate kernel entries can be evaluated rapidly. Constructing the H-matrix from kernel entries requires a number of elementary operations that scales polylog linearly in $N$. Using the H-matrix to perform linear algebra operations also scales polylog linearly in $N$, though the details of these costs depend heavily on the type of H-matrix and operation being performed.

\section{Numerical results}
\label{sec:numerical_results}

We use the PSF-based method to approximate the Newton (or Gauss-Newton) Hessians
in inverse problems governed by PDEs which model steady state ice
sheet flow~\cite{PetraEtAl12} (Section~\ref{sec:ice}) and advective-diffusive transport of a contaminant~\cite{PetraStadler11}
(Section~\ref{sec:adv}), and to approximate the integral kernel in a blur problem that is not based on PDEs (Section~\ref{sec:frog}). These problems are described in detail in their respective sections. 

In both PDE-based inverse problems (Sections~\ref{sec:ice} and~\ref{sec:adv}), to reconstruct the unknown parameter fields, denoted $q$, the inverse problems are formulated as nonlinear least squares
optimization problems, whose objective functions consist of a data misfit term
(between the observations and model output) and a bi-Laplacian
regularization term following~\cite{VillaPetraGhattas21}. The
regularization is centered at a constant function $q_0(x)$. To
mitigate boundary effects
%For stokes q_0(x) = 10.5$%
we use a constant coefficient Robin boundary condition as
in~\cite{Roininen14}. The parameters for the bi-Laplacian operator are chosen
so that the Green's function of the Hessian of the regularization has
a characteristic length of $0.25$ of the domain radius. For the
specific setup, we refer the reader to~\cite[Section
  2.2]{VillaPetraGhattas21}. In all numerical results we choose the
regularization parameter (which controls the overall strength of the regularization) using the Morozov discrepancy
principle~\cite{Vogel02}. 

We solve the ice sheet inverse problem with an inexact Newton preconditioned
conjugate gradient (PCG) scheme and a globalizing Armijo line search~\cite{NocedalWright99}. The Newton search directions,
$\mathbf{\searchdir}$, are obtained by solving 
\begin{equation}
	\label{eq:newton_system}
	\mathbf{H} \mathbf{\searchdir} = - \mathbf{g} \qquad \text{or} \qquad \mathbf{H}_\text{gn} \mathbf{\searchdir} = - \mathbf{g},
\end{equation}
wherein we choose the initial guess as the discretization of the constant
function $q_0$. Here, $\mathbf{g}$, $\mathbf{H}$ and $\mathbf{H}_\text{gn}$ are
the discretized gradient, Hessian, and Gauss-Newton Hessian of the
inverse problem objective function, respectively, evaluated at the current Newton iterate.
To ensure positive definiteness of the Hessian we use
$\mathbf{H}_\text{gn}$ for the first five iterations, and $\mathbf{H}$ for all
subsequent iterations. The Newton iterations are terminated when
$\|\mathbf{g}\| < 10^{-6} \|\mathbf{g}_0\|$, where $\mathbf{g}_0$ is the gradient evaluated at the initial guess. Systems~\eqref{eq:newton_system} are solved inexactly using an inner PCG
iteration, which is terminated early based on the Eisenstat-Walker \cite{EisenstatWalker96} and Steihaug \cite{Steihaug83} conditions. The inverse problem governed by the advection-diffusion PDE is linear, hence Newton's method converges in one iteration. In this case the Newton linear system,~\eqref{eq:newton_system}, is solved using PCG, using termination tolerances described in Section \ref{sec:adv}.

We use the framework described in this paper to generate Hessian preconditioners. We build H-matrix approximations, $\mathbf{A}_H$, of the data misfit Gauss-Newton Hessian (the term in $\mathbf{H}_\text{gn}$ that arises from the data misfit). The 
approximations are indicated by ``PSF ($\nbatch$)'', where $\nbatch$ is the number of impulse response batches used to build the approximation. The Hessian of the regularization term is a combination of
stiffness and mass matrices, which are sparse. Therefore, we form
H-matrix representations of these matrices and combine them into an H-matrix approximation of the regularization term in the Hessian, $\mathbf{R}_H$, using
standard sparse H-matrix techniques. Then, H-matrix approximations of the Gauss-Newton Hessian, 
$\mathbf{H}_\text{gn}\approx \widetilde{\mathbf{H}}:=\mathbf{A}_H + \mathbf{R}_H$,
are formed by adding $\mathbf{A}_H$ to $\mathbf{R}_H$ using fast H-matrix arithmetic. We modify $\widetilde{\mathbf{H}}$ to be (approximately) symmetric positive semi-definite via the procedure described in Section~\ref{sec:flip_eigs}. We factor $\mathbf{\preconditioner}$
using fast H-matrix methods, then use the factorization as a
preconditioner. We approximate $\mathbf{H}_\text{gn}$ rather than $\mathbf{H}$ because $\mathbf{H}$ more often has negative values in its integral kernel. The numerical results show that~$\widetilde{\mathbf{H}}$ is a good preconditioner for both~$\mathbf{H}_\text{gn}$ and~$\mathbf{H}$.

\subsection{Example 1: Inversion for the basal friction coefficient in an ice sheet flow problem}
\label{sec:ice}

\begin{figure}
	\begin{center} 
    \begin{subfigure}{0.99\textwidth}
    \begin{center} 
	\includegraphics[scale=0.2]{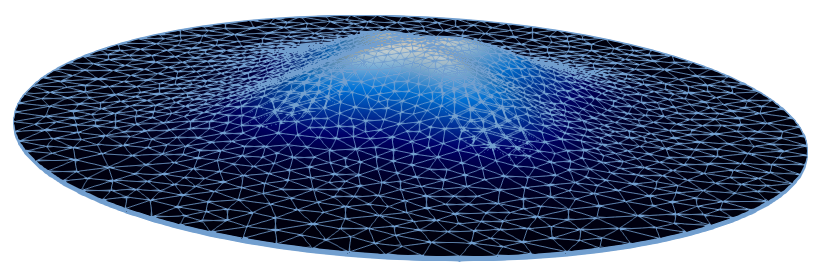}
	\end{center}
	\caption{Ice sheet model geometry}
	\label{fig:ice_mountain_mesh}
    \end{subfigure} \\
	\begin{subfigure}{0.49\textwidth}
		\includegraphics[scale=0.2]{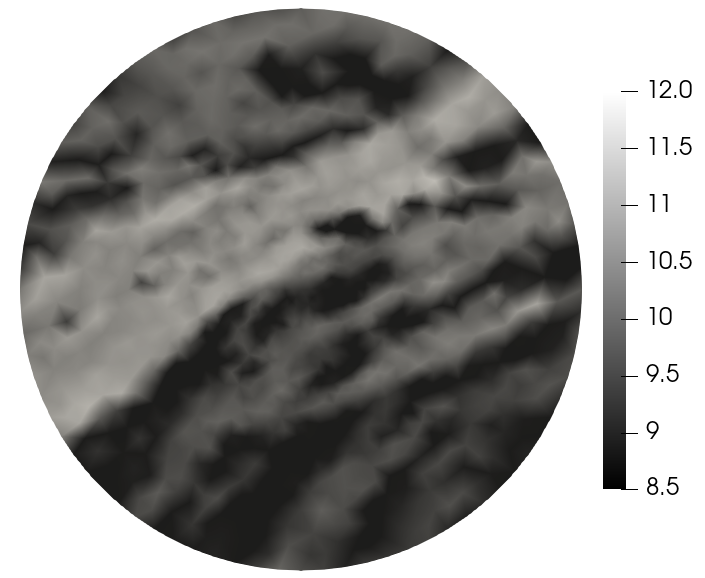}
		\caption{$\basalfriction_\text{true}$}
		\label{fig:true_beta}
	\end{subfigure}
	\begin{subfigure}{0.49\textwidth}
		\includegraphics[scale=0.2]{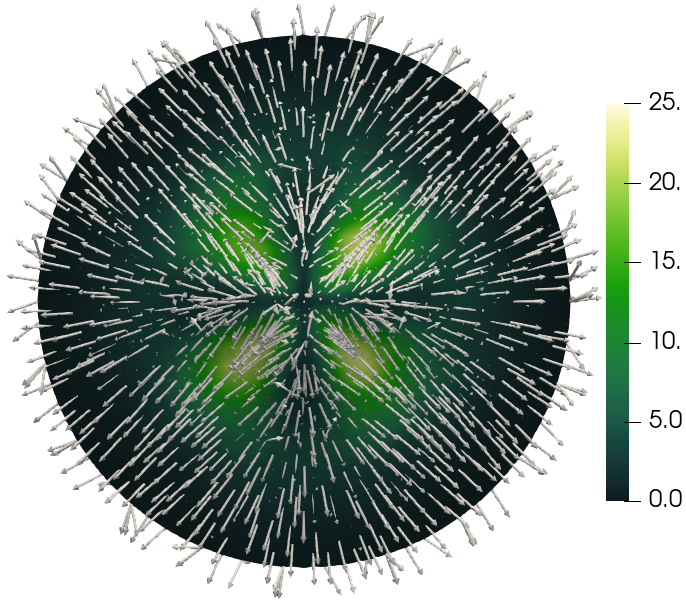}
		\caption{$\velocity_\text{true}$}
		\label{fig:stokes_velocity}
	\end{subfigure}
    \end{center}
	\caption{(Ice sheet) 
		(\ref{fig:ice_mountain_mesh}) Bird's eye view of the ice sheet discretized by a mesh of tetrahedra. Color indicates the height of the base of the ice sheet (i.e., the mountain topography). The radius of the domain is $10^4$ meters, the maximum height of the mountain is $2.1 \times 10^3$ meters, and the average thickness of the ice sheet is 250 meters.		
		(\ref{fig:true_beta}) True parameter, $\basalfriction_\text{true}$. (\ref{fig:stokes_velocity}) True velocity, $\velocity_\text{true}$. Arrows indicate the direction of $\velocity_\text{true}$ and color indicates the magnitude of $\velocity_\text{true}$.}
	\label{fig:true_beta_u}
\end{figure}

For this example, we consider a sheet of ice flowing down a mountain (see Figure~\ref{fig:ice_mountain_mesh}). Given observations of the tangential component of the ice velocity on the top surface of the ice, we invert for the logarithm of the unknown spatially varying basal friction Robin coefficient field, which governs the resistance to sliding along the base of the ice sheet. 
The setup, which we briefly summarize, follows~\cite{IsaacEtAl15,PetraEtAl12}.
The region of ice is denoted by
$\icedomain \subset\mathbb{R}^{3}$. The basal, lateral and top parts of the boundary $\partial \icedomain$ are denoted by $\Gamma_{b}$, $\Gamma_{l}$, and $\Gamma_{t}$, respectively.
The governing equations are the linear incompressible Stokes
equations,
\begin{subequations}\label{Stokeseqn}
	\begin{align}
          -\nabla \cdot \stress(v, p)=\bodyforce \text{ and }
          \nabla\cdot \velocity &=0 \quad \text{ in }\icedomain, \label{Stokeseqn:1}\\
	\stress(v, p)\normalvec
	&=0\,\,\,\,\,\,\text{ on }\Gamma_{t}, \label{Stokeseqn:3}\\
	\velocity \cdot\normalvec =0 \text{ and } \tangentop\left(\stress(v, p)\normalvec+\exp\left(\basalfriction\right)\velocity\right)
	&=0\,\,\,\,\,\,\text{ on }\Gamma_{b},\label{Stokeseqn:4} \\
	\stress(v, p)\normalvec+\stokesrobincoeff \velocity &=0\,\,\,\,\,\,\text{ on }\Gamma_{l}.\label{Stokeseqn:5}	
	\end{align}
\end{subequations}
The solution to these equations is the pair $(\velocity, \pressure)$, where $\velocity$ is the ice flow velocity field\footnote{We do not use bold to denote vector or tensor fields to avoid confusion with vectors that arise from finite element discretizations, which are already denoted with bold.} and $\pressure$ is the pressure field. Here, $\basalfriction$ is the unknown logarithmic
basal friction field (large $\basalfriction$ corresponds to large resistance to sliding) defined on the surface $\Gamma_b$. The quantity $\bodyforce$ is the body force density due to gravity, $\stokesrobincoeff=10^6$ is a Robin boundary condition constant,  $\normalvec$ is the outward unit normal and $\tangentop$ is the tangential projection operator that restricts a vector field to its tangential component along the boundary. 
We employ a Newtonian constitutive law, $\stress(v, p)= 2\eta
\dot{\strain}(v) -I\pressure$,
where $\stress$ is the stress tensor and
$\dot{\strain}(v)= \frac{1}{2}\left(
\nabla\velocity+\nabla\velocity^{\top} \right)$ is the strain rate tensor~\cite{IsaacEtAl15}. Here, $\eta$ is the viscosity and $I$ is the identity operator. 
Note that while the PDE is linear, the parameter-to-solution map, $\basalfriction \mapsto (\velocity, \pressure)$, is nonlinear.

\begin{figure}
	\begin{subfigure}{0.32\textwidth}
		\includegraphics[scale=0.17]{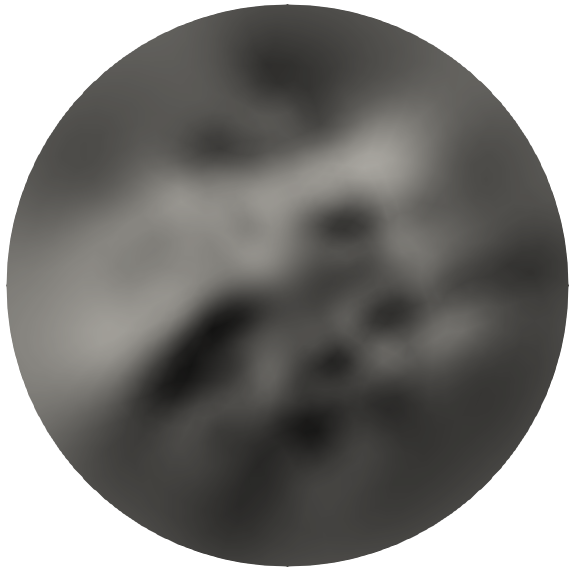}
	\end{subfigure}
	\begin{subfigure}{0.32\textwidth}
		\includegraphics[scale=0.17]{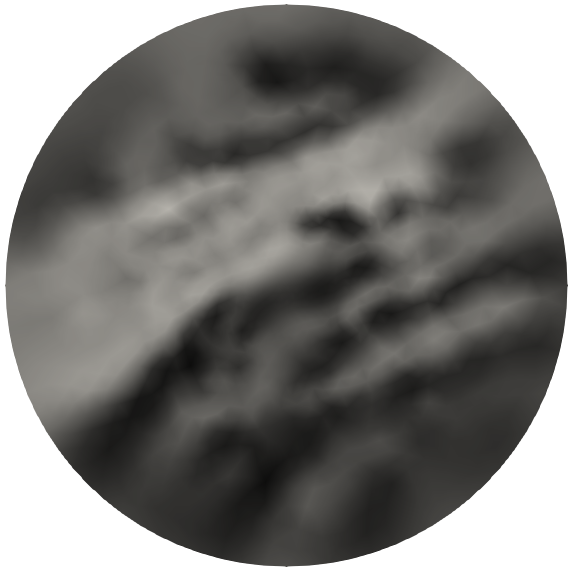}
	\end{subfigure}
	\begin{subfigure}{0.32\textwidth}
	\includegraphics[scale=0.17]{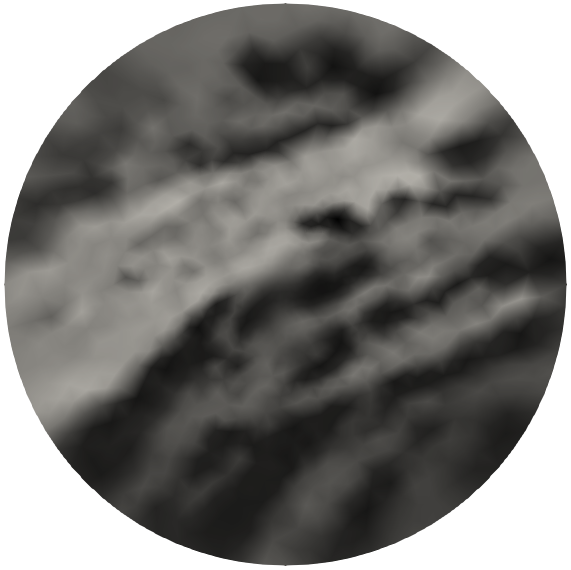}
    \end{subfigure}
	\caption{(Ice sheet) The log basal friction parameter, with color scale as in Figure~\ref{fig:true_beta}, computed from the PDE constrained optimization problem with noise levels: $25\%$ (left), $5.0\%$ (middle), and $1.0\%$ (right).}
	\label{fig:stokes_reconstructions} 
\end{figure}

The pressure, $\pressure$, is discretized with first order scalar continuous Galerkin finite elements defined on a mesh of tetrahedra. The velocity, $\velocity$, is discretized with second order continuous Galerkin finite elements on the same mesh. The parameter $\basalfriction$ is discretized with first order scalar continuous Galerkin finite elements on the mesh of triangles that results from restricting the tetrahedral mesh to the basal boundary, $\Gamma_b$. Note that $\Gamma_b$ is a two-dimensional surface embedded in three dimensions due to the mountain topography. The PSF-based method involves translating impulse responses. Hence it requires either a flat domain, or a notion of local parallel transport. We therefore generate a flattened version of $\Gamma_b$, denoted by $\Omega \subset \mathbb{R}^2$, by ignoring the height coordinate. The parameter $\basalfriction$ is viewed as a function on $\Gamma_b$ for the purpose of solving the Stokes equations, and as a function on $\Omega$ for the purpose of building Hessian approximations and defining the regularization.
The observations are generated by adding multiplicative
Gaussian noise to the tangential component of the velocity field
restricted to the top surface of the geometry. We use $5\%$ noise in all cases, except for Figure~\ref{fig:stokes_reconstructions} and Table \ref{tab:condition_number} where the noise is varied from $1\%$ to $25\%$ and the regularization is determined by the Morozov discrepancy principle for each noise level. The true basal friction
coefficient and resulting velocity fields, which are obtained by
solving~\eqref{Stokeseqn}, are shown in Figure~\ref{fig:true_beta_u}.

\begin{table}
	{\footnotesize
		\begin{center}
			\begingroup
			\setlength{\tabcolsep}{3pt}
			\renewcommand{\arraystretch}{1.1}
			\begin{tabular}{c| c c c | c c c | c c c}
				&  \multicolumn{3}{|c|}{PSF (5)} & \multicolumn{3}{|c|}{REG} & \multicolumn{3}{|c}{NONE} \\
				\hline
				Iter & 
				\#CG & \#Stokes & $\|\mathbf{g}\|$ & 
				\#CG & \#Stokes & $\|\mathbf{g}\|$ & 
				\#CG & \#Stokes & $\|\mathbf{g}\|$ \\
				0 &
				1 & 4 & 1.9e+7 &
				3 & 8 & 1.9e+7 &
				1 & 4 & 1.9e+7 \\
				1 &
				2 & 6  & 6.1e+6 &
				8 & 18 & 8.4e+6 &
				2 & 6  & 6.1e+6 \\
				2 &
				4 & 10 & 2.6e+6 &
				16 & 34 & 4.1e+6 &
				4 & 10 & 2.6e+6 \\
				3 &
				2 & 6+22 & 6.9e+5 &
				34 & 70 & 1.8e+6 &
				14 & 30 & 6.9e+5 \\
				4 &
				3 & 8 & 4.4e+4 &
				52 & 106 & 5.6e+5 &
				29 & 60 & 1.3e+5 \\
				5 &
				5 & 12 & 2.2e+3 &
				79 & 160 & 9.4e+4 &
				38 & 78 & 1.0e+4 \\
				6 &
				0 & 2 & 1.1e+1 &
				102 & 206 & 6.5e+3 &
				58 & 118 & 1.8e+2 \\
				7 &
				--- & --- & --- &
				151 & 304 & 1.2e+2 &
				0 & 2 & 5.5e-1 \\
				8 & 
				--- & --- & --- &
				0 & 2 & 2.9e-1 &
				--- & --- & --- \\
				\hline
				Total & 
				17 & 70 & --- &
				445 & 908 & --- &
				146 & 308 & --- \\
			\end{tabular}
			\endgroup
		\end{center}
	}
	\caption{(Ice sheet) Convergence history for solving the Stokes inverse problem using inexact Newton PCG to tolerance $10^{-6}$. 
	Preconditioners shown are the PSF-based method with five batches (PSF (5)) constructed at the third iteration, regularization preconditioning (REG), and no preconditioning (NONE). Columns \#CG show the number of PCG iterations used to solve the Newton system for $\mathbf{\widehat{\basalfriction}}$. Columns $\|\mathbf{g}\|$ show the $l^2$ norm of the gradient at $\mathbf{\basalfriction}$. Columns \#Stokes show the total number of Stokes PDE solves performed in each Newton iteration. 
	Under PSF (5) and in row Iter 3, we write $6+22$ to indicate that $6$ Stokes solves were used during the standard course of the iteration, and $22$ Stokes solves were used to build the PSF (5) preconditioner.}
	\label{tab:newton_convergence_table}
\end{table}

Table~\ref{tab:newton_convergence_table} shows the performance of the preconditioner for accelerating the solution of the optimization problem to reconstruct $\mathbf{q}$ from observations with $5\%$ noise. We build the PSF (5) preconditioner in the third Gauss-Newton iteration, and reuse it for all subsequent Gauss-Newton and Newton iterations. No preconditioning is used in the iterations before the PSF (5) preconditioner is built. We compare the PSF-based method with the most commonly used existing preconditioners: no preconditioning (NONE), and preconditioning by the regularization term in the Hessian (REG). The results show that using PSF (5) reduces the total number of Stokes PDE solves to 70, as compared to 908 for regularization preconditioning and 308 for no preconditioning, a reduction in cost of roughly $5\times$--$10\times$. For problems with a larger physical domain and correspondingly more observations, such as continental scale ice sheet inversion, the speedup will be even greater. This is because the rank of the data misfit Hessian will increase, while the locality of the impulse responses will remain the same. In Figure~\ref{fig:stokes_reconstructions} we show reconstructions for $1\%$, $5\%$, and $25\%$ noise.
\begin{figure}
	\begin{subfigure}{0.52\textwidth}
		\centering
		\includegraphics[scale=0.9]{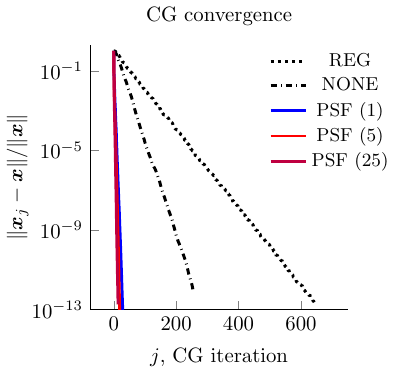}
	\end{subfigure}
	\begin{subfigure}{0.47\textwidth}
	\begin{center}
		\includegraphics[scale=0.9]{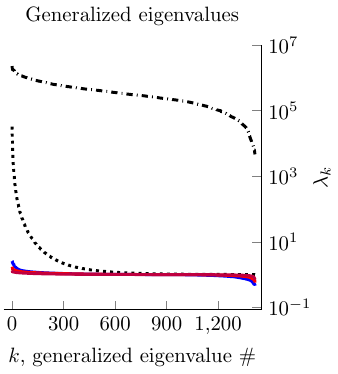}
	\end{center}
    \end{subfigure}
	\caption{(Ice sheet) Left: Convergence history for solving $\mathbf{H} \mathbf{x} = \mathbf{b}$ using PCG, where $\mathbf{b}$ has i.i.d.\ random entries drawn from the standard Gaussian distribution and $\mathbf{H}$ is evaluated at the solution of the inverse problem. 
	Results in these figures are shown for the PSF-based preconditioners with 1, 5, and 25 batches (PSF (1), PSF (5), and PSF (25), respectively), regularization preconditioning (REG), and no preconditioning (NONE). The preconditioner is constructed using $\mathbf{H}_\text{gn}$.
	Right: Generalized eigenvalues for generalized eigenvalue problem $\mathbf{H} \mathbf{u}_{k}=\lambda_{k} \mathbf{\preconditioner} \mathbf{u}_{k}$. Here, $\mathbf{H}$ is the Hessian and the matrices $\mathbf{\preconditioner}$ are the same Hessian approximations used in the left sub-figure, with NONE corresponding to the identity matrix.}
	\label{fig:krylov_convergence}
\end{figure} 

Next, we build PSF~(1), PSF~(5), and PSF~(25) preconditioners based on the Gauss-Newton Hessian evaluated at the converged solution $\mathbf{\basalfriction}$ (note: $k$ in PSF~($k$) refers to the number of batches; this is not to be confused with the noise levels which range over the same numerical values). We use PCG to solve a linear system with the Hessian as the coefficient operator and a right hand side vector with random independent and identically distributed (i.i.d.) entries drawn from the standard Gaussian distribution. In Figure~\ref{fig:krylov_convergence} (left) we compare the convergence of PCG for solving this linear system using the PSF (1), PSF (5), PSF (25), REG, and NONE preconditioners. PCG converges fastest with the PSF-based preconditioners, with PSF (25) converging fastest, followed by PSF (5), followed by PSF (1), as expected. 
In Figure~\ref{fig:krylov_convergence} (right) we show the generalized eigenvalues for the generalized eigenvalue problem
$
\mathbf{H} \mathbf{u} = \lambda \mathbf{\preconditioner} \mathbf{u}.
$ 
The matrix $\mathbf{\preconditioner}$ is one of the PSF (1), PSF (5), or PSF (25) Gauss-Newton Hessian approximations, the regularization Hessian (REG), or the identity matrix (NONE). With the PSF-based preconditioners, the generalized eigenvalues cluster near one, with more batches yielding better clustering. 

In Table~\ref{tab:condition_number}, we show the condition number of the preconditioned Hessian for noise levels ranging from $1\%$ to $25\%$. Note that the condition number using PSF-based preconditioners is extremely small (ranging between 1 and 10) and relatively stable over this range of noise levels. As expected, PSF (25) outperforms PSF (5), which outperforms PSF (1). All PSF-based preconditioners outperform regularization and no preconditioning by several orders of magnitude for all noise levels. 

\begin{table}
	\begin{center}
			{\small
			\begingroup
			\setlength{\tabcolsep}{4pt}
			\renewcommand{\arraystretch}{1.25}
			\begin{tabular}{c| c c c c c}
				noise    & \multicolumn{5}{c}{COND$(\mathbf{\preconditioner}^{-1} \mathbf{H})$ } \\ \cline{2-6}
				level    & REG     &	NONE  & PSF $(1)$ & PSF $(5)$ & PSF $(25)$ \\ \hline 
				$25\%$   & 1.01e+3 & 2.96e+3  & 1.34e+0   & 1.30e+0   & 1.18e+0    \\ 
				$11\%$   & 7.40e+3 & 1.05e+3  & 2.27e+0   & 1.55e+0   & 1.31e+0    \\   
				$5.0\%$  & 3.29e+4 & 4.96e+2  & 5.61e+0   & 3.06e+0   & 1.92e+0    \\ 
				$2.2\%$  & 1.66e+5 & 8.89e+2  & 1.58e+1   & 8.07e+0   & 4.03e+0    \\  
				$1.0\%$  & 5.36e+5 & 1.61e+3  & 7.17e+1   & 1.93e+1   & 9.19e+0    \\   
			\end{tabular}
			\endgroup
			}
		\end{center}
	\caption{(Ice sheet) Condition number for $\mathbf{\preconditioner}^{-1} \mathbf{H}$ for the PSF-based preconditioners with 1, 5, and 25 batches (PSF (1), PSF (5), and PSF (25), respectively), no preconditioner (NONE) and regularization preconditioning (REG). All operators are evaluated at the soutions of the inverse problems for 
	their respective noise levels.}
	\label{tab:condition_number}
\end{table}

\subsection{Example 2: Inversion for the initial condition in an advective-diffusive transport problem}
\label{sec:adv}

\begin{figure}
	\begin{center}
	\begin{subfigure}{0.24\textwidth}
		\begin{center}
		\includegraphics[scale=0.12]{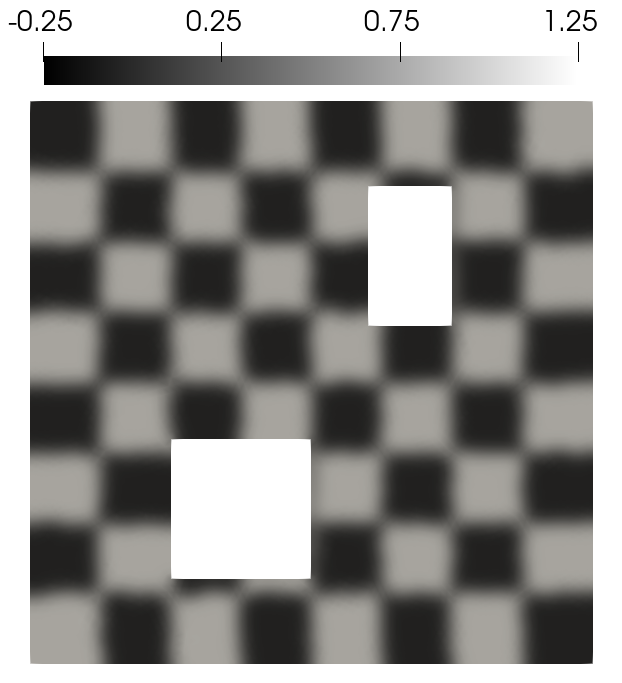}
	\end{center}
				\caption{True $q$}
				\label{fig:adv_initial_condition}
	\end{subfigure}
	\begin{subfigure}{0.24\textwidth}
		\begin{center}
		\includegraphics[scale=0.12]{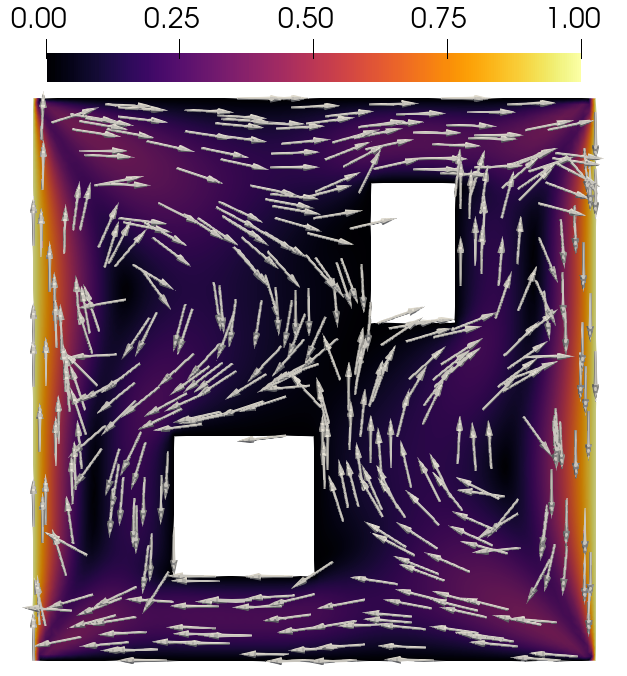}
		\end{center}
				\caption{$v$}
				\label{fig:adv_velocity}
	\end{subfigure}
	\begin{subfigure}{0.24\textwidth}
		\begin{center}
		\includegraphics[scale=0.12]{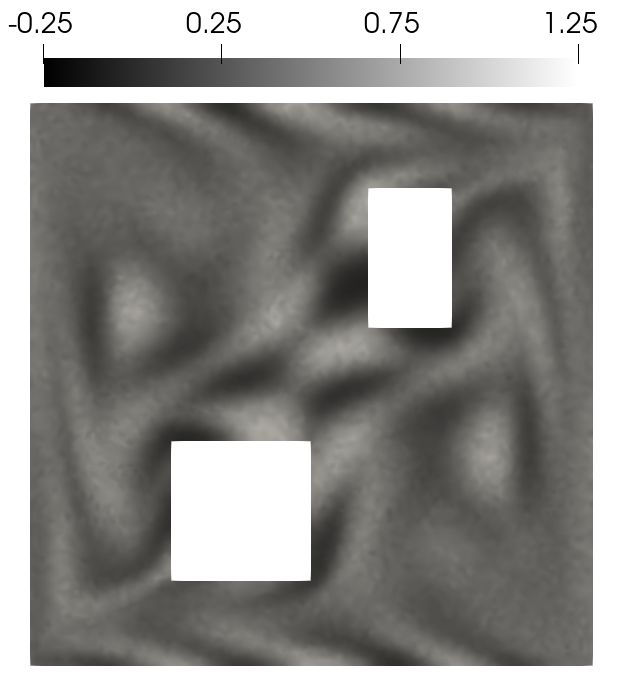}
		\end{center}
				\caption{Noisy observations}
				\label{fig:adv_observations}
	\end{subfigure}
	\begin{subfigure}{0.24\textwidth}
		\begin{center}
		\includegraphics[scale=0.12]{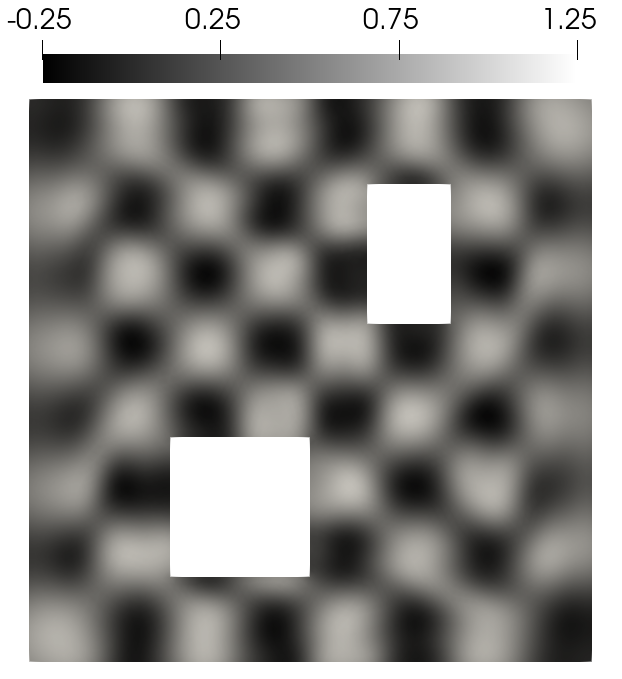}
		\end{center}
				\caption{Reconstructed $q$}
				\label{fig:adv_reconstruction}
	\end{subfigure}
	\end{center}
	\caption{(Advective-diffusive transport) (\ref{fig:adv_initial_condition}) True initial condition. (\ref{fig:adv_velocity}) Velocity field. Color indicates magnitude of velocity vector. (\ref{fig:adv_observations}) Noisy observations of concentration at the final time. (\ref{fig:adv_reconstruction}) Reconstructed initial condition.}
	\label{fig:adv_inversion1} 
\end{figure}

\begin{figure}
	\centering
	{
		\begin{tabular}{ccc}
			& $T=0.5$ 
			& $T=2.0$ \\
			$\kappa=1.0$e${-4}$  & \includegraphics[scale=0.12, valign=c]{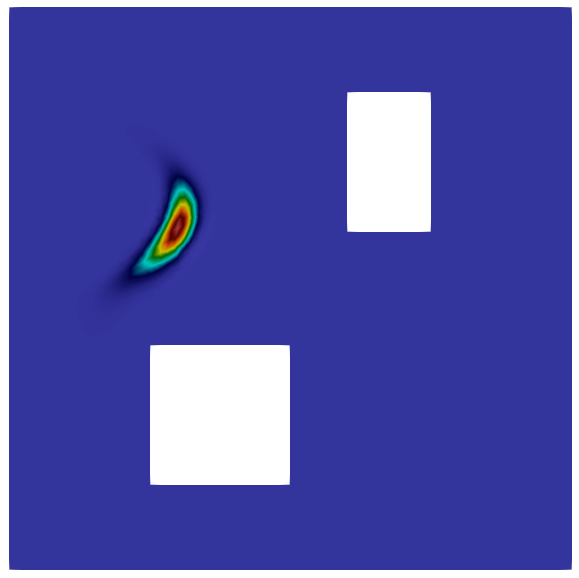} 
			& \includegraphics[scale=0.12, valign=c]{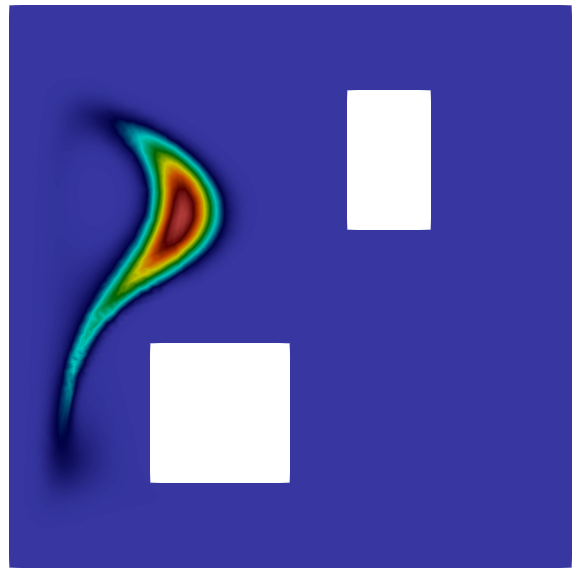} \\
			$\kappa=1.0$e${-3}$ & \includegraphics[scale=0.12, valign=c]{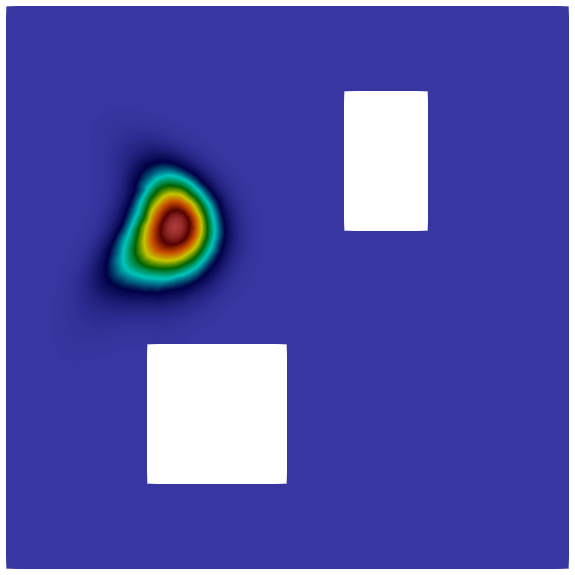} 
			& \includegraphics[scale=0.12, valign=c]{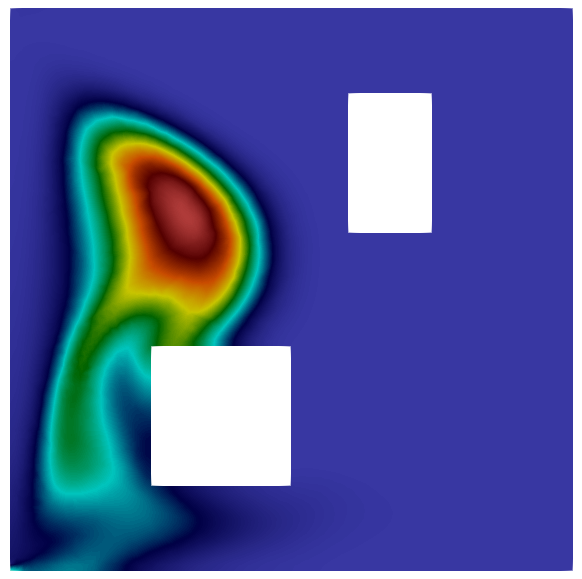}
		\end{tabular}
	}
	\caption{(Advective-diffusive transport) Impulse responses for small and large diffusion parameters $\kappa$ and simulation times $T$.}
	\label{tbl:adv_vary_kappa_tf}
\end{figure}

Here, we consider a time-dependent advection-diffusion equation in
which we seek to infer the unknown spatially varying initial condition, $\ipar$, from noisy observation of the full state  
at a final time, $T$. This PDE
models advective-diffusive transport in a domain $\advdomain \subset \R^d$, which is
depicted in Figure~\ref{fig:adv_inversion1}. In this case, the state, $\concentration(x,t)$, could be interpreted as the concentration of a contaminant. The problem description
below closely follows~\cite{PetraStadler11,VillaPetraGhattas21}. The domain
boundaries $\partial \advdomain$ include the outer boundaries as well as the
internal boundaries of the rectangles, which represent buildings.  The
parameter-to-observable map $\iFF$ in this case maps an initial
condition $\ipar \in L^2(\advdomain)$ to the concentration field at a final time, $\concentration(x, T)$,  through solution of the advection-diffusion
equation given by
\begin{equation}\label{eq:ad}
  \begin{aligned}
    \concentration_t - \kappa\Delta \concentration + v\cdot\nabla \concentration &= 0 & \quad&\text{in
    }\advdomain\times (0,T), \\
    \concentration(\cdot, 0) &= \ipar  &&\text{in } \advdomain , \\
    \kappa\nabla \concentration\cdot \nu &= 0 &&\text{on } \partial \advdomain \times (0,T).
  \end{aligned}
\end{equation}
Here, $\kappa > 0$ is a diffusivity coefficient, $\nu$ is the boundary unit normal vector, and $T > 0$ is the
final time.  The velocity field, $v:\advdomain \rightarrow \mathbb{R}^d$, is computed by solving the
steady-state Navier-Stokes equations for a two dimensional flow with
Reynolds number 50, with boundary conditions $v(x)=(0,1)$ on the left boundary, $v(x)=(0,-1)$ on the right boundary, and $v(x)=(0,0)$ on the top and bottom boundaries, as in~\cite[Section 3]{PetraStadler11}. We use a checkerboard image for the initial condition (Figure~\ref{fig:adv_initial_condition}) and add $5\%$ multiplicative noise to generate a synthetic observation at the final time, $T$. The initial condition, velocity field, noisy observations, and reconstructed initial condition are shown in Figure~\ref{fig:adv_inversion1}. We use $\kappa=3.2$e${-1}$ and $T=1.0$ for all results, except for Table \ref{tab:adv_convergence} and Figure~\ref{tbl:adv_vary_kappa_tf} where we vary $\kappa$ and $T$.

In Table~\ref{tab:adv_convergence} we show the number of PCG iterations, $j$, required to solve the Newton linear system to a relative error tolerance of $\|\widehat{\mathbf{q}} - \widehat{\mathbf{q}}_j\| < 10^{-6} \|\widehat{\mathbf{q}}\|$. The solution of the Newton system to which we compare, $\widehat{\mathbf{q}}$, is found via another PCG iteration with a relative residual tolerance of $10^{-11}$. We show results for $T$ ranging from $0.5$ to $2.0$ and $\kappa$ ranging from $10^{-4}$ to $10^{-3}$ using the PSF-based preconditioners with $1$, $5$, and $25$ batches, regularization preconditioning, and no preconditioning. The results show that PSF-based preconditioning outperforms regularization preconditioning and no preconditioning in all cases except one. The exception is $T=2.0$ and $\kappa=1.0$e${-3}$, in which PSF (1) performs slightly worse than regularization preconditioning but better than no preconditioning. Adding more batches yields better results, and the impact of adding more batches is more pronounced here than in the ice sheet example. For example, in the mid-range values $T=1.0$ and $\kappa=3.2$e${-4}$, PSF (1), PSF (5), and PSF (25) require $1.3\times$, $2.3\times$, and $5.1\times$ fewer PCG iterations, respectively, as compared to no preconditioning, and exhibit greater improvements as compared to regularization preconditioning. The PSF preconditioners perform best in the high rank regime where $T$ is small, which makes sense given that short simulation times yield more localized impulse responses (see Figure~\ref{tab:adv_convergence}). For example, for $\kappa=3.2$e${-4}$ using the PSF (5) preconditioner yields 140, 202, and 379 iterations for $T=0.5$, $1.0$, and $2.0$, respectively. The performance of the PSF preconditioners as a function of $\kappa$ does not have as clear of a trend. Reducing $\kappa$ makes the impulse responses thinner and hence easier to fit in batches, but also increases the complexity of the impulse response shapes, which may reduce the accuracy of the RBF interpolation. The greatest improvements are seen for $T=0.5$ and $\kappa=1$e${-3}$, for which PSF (25) requires roughly $10\times$ and $20\times$ fewer PCG iterations than no preconditioning and regularization preconditioning, respectively.

\begin{figure}
	\begin{subfigure}{0.52\textwidth}
		\centering
		\includegraphics[scale=0.9]{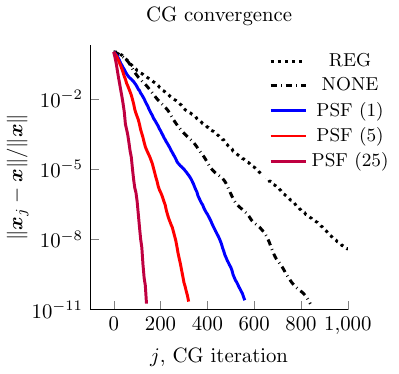}
	\end{subfigure}
	\begin{subfigure}{0.47\textwidth}
		\centering
		\includegraphics[scale=0.9]{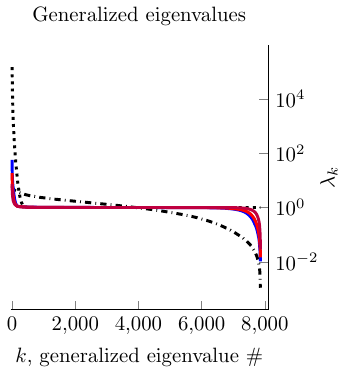}
	\end{subfigure}
	\caption{(Advective-diffusive transport) Left: Convergence history for solving $\mathbf{H} \mathbf{x} = \mathbf{b}$ using PCG, where $\mathbf{b}$ has i.i.d.\ random entries drawn from the standard Gaussian distribution. Right: Generalized eigenvalues for generalized eigenvalue problem $\mathbf{H} \mathbf{u}_{k}=\lambda_{k} \mathbf{\preconditioner} \mathbf{u}_{k}$. Here, $\mathbf{H}$ is the Hessian and the preconditioner, $\mathbf{\preconditioner}$, is the PSF-based approximation for 1, 5, or 25 batches (PSF (1), PSF (5), and PSF (25), respectively), the regularization Hessian (REG), or the identity matrix (NONE).}
	\label{fig:adv_krylov_convergence}
\end{figure}

In Figure~\ref{fig:adv_krylov_convergence} (left) we show the convergence of PCG for solving $\mathbf{H}\mathbf{x}=\mathbf{b}$, where $\mathbf{b}$ has i.i.d.\ random entries drawn from a standard Gaussian distribution. The preconditioners used, $\widetilde{\mathbf{H}}$, are the PSF-based preconditioners with $1$, $5$, or $25$ batches, the regularization Hessian, and the identity matrix (i.e., no preconditioning). The results show that PCG converges fastest with the PSF-based preconditioners, with more batches yielding faster convergence. In Figure~\ref{fig:adv_krylov_convergence} (right), we show the eigenvalues for the generalized eigenvalue problem $\mathbf{H}\mathbf{u}=\lambda \widetilde{\mathbf{H}}\mathbf{u}$, where the $\widetilde{\mathbf{H}}$ are the preconditioners stated above. With the PSF-based preconditioners the eigenvalues cluster near one, and more batches yields better clustering. With the regularization preconditioner the trailing eigenvalues cluster near one, while the leading eigenvalues are amplified.

\begin{table}
	\begin{center}
		{\small
		\begingroup
		\setlength{\tabcolsep}{8pt}
		\renewcommand{\arraystretch}{1.0}
		\begin{tabular}{c| c | c c c c c}
			& $\kappa$    & REG     &	NONE   & PSF $(1)$ & PSF $(5)$ & PSF $(25)$ \\ 
			\hline 
			& 1.0e-4    & 584 & 317  & 311   & 151   & 56    \\ 
			$T=0.5$ & 3.2e-4    & 685 & 311  & 233   & 140   & 44    \\ 
			& 1.0e-3    & 702 & 324  & 122   & 71   & 33    \\ 
			\hline
			& 1.0e-4    & 634 & 449  & 539   & 288   & 100    \\ 
			$T=1.0$ & 3.2e-4    & 681 & 459  & 350   & 202   & 90    \\ 
			& 1.0e-3    & 574 & 520  & 266   & 260   & 208    \\ 
			\hline
			& 1.0e-4    & 609 & 591  & 548   & 520   & 165    \\ 
			$T=2.0$ & 3.2e-4    & 524 & 645  & 318   & 379   & 170    \\ 
			& 1.0e-3    & 349 & 786  & 381   & 262   & 158    \\ 
		\end{tabular}
		\endgroup
		}
	\end{center}
	\caption{(Advective-diffusive transport) Number of PCG iterations required to solve the Newton linear system to tolerance $||\widehat{\mathbf{q}}_j - \widehat{\mathbf{q}}|| < 10^{-6}||\widehat{\mathbf{q}}||$, where $\widehat{\mathbf{q}}_j$ is the $j$th iterate, and $\widehat{\mathbf{q}}$ is the solution of the Newton linear system. Iteration counts are shown for a variety of different diffusion parameters $\kappa$, simulation times $T$, and preconditioners.}
	\label{tab:adv_convergence}
\end{table}

\subsection{Example 3: Spatially varying blurring problem}
\label{sec:frog}

Here, we define a PDE-free spatially varying blur problem, in which the impulse response, $\phi_x$,
is a bumpy blob that is centered
near $x$, and is rotated and scaled in a manner that depends on $x$
(see Figure \ref{fig:frog_kernel_impulse_responses}). This blur problem is used in
Sections~\ref{sec:intro} and~\ref{sec:prelims} to visually illustrate
various stages and aspects of the PSF-based method, and robustness (or lack thereof) of the PSF-based method to violations of the non-negative kernel assumption (Section~\ref{sec:conditions_2}). In this section, the blur problem is used to compare the PSF-based method to hierarchical off-diagonal low rank (HODLR) and global low rank (GLR) methods, to study convergence of the PSF-based method, and to investigate the effect of the ellipsoid size parameter $\tau$.
The closed form
expression for the integral kernel is given by
\begin{align}
\label{eq:frog_kernel}
\Phi(y,x) = (1 - a f(y,x)) g(x) \exp \left(-\frac{1}{2}(h(y,x)^T
C^{-1} h(y,x) \right),
\end{align}
where $f(y,x) = \cos\left(h^1(y,x) / \sqrt{c_1 /
	2}\right)\sin\left(h^2(y,x) / \sqrt{c_2 / 2}\right)$, with
$h^i(y,x)$ the $i^{th}$ component of $R(\theta(x)) (y-x)$, with
$R(\theta(x))$ a two-dimensional rotation matrix by angle $\theta(x) =
(x^1 + x^2) \pi / 2$, $g(x) = x^1(1-x^1)x^2(1-x^2)$, and with $C = L^2
\operatorname{diag}(c_1, c_2)$. The constant $L$ controls the width of
the blob, and $c_1/c_2$ controls its aspect ratio. The constant $a$
represents deviation from a Gaussian. When $a=0$, $\phi_x$ is a
Gaussian, and as $a$ increases, $\phi_x$ becomes non-Gaussian. When $a
> 1$, the integral kernel contains negative values, which allows us to
study the robustness of the PSF-based method to violations of Assumption 3
(Section~\ref{sec:conditions_2}).

\begin{table}
	\begin{center}
		{\small
			\begingroup
			\setlength{\tabcolsep}{8pt}
			\renewcommand{\arraystretch}{1.0}
			\begin{tabular}{c| c | c c c c c}
				& Error tol.    &	\#applies PSF   & \#applies HODLR 	& \#applies RSVD 	\\
				\hline 
				& 20\%    		& 11  				& 592				& 354				\\
				$L=1$ 			& 10\% 			& 16  				& 772				& 520				\\
				& 5\% 			& 22   				& 924				& 674				\\
				\hline
				& 20\%    		& 8  				& 852 				& 1316				\\
				$L=1/2$ 		& 10\%    		& 9  				& 1144				& 1916				\\
				& 5\% 	  		& 12   				& 1404				& 2456				\\
				\hline
				& 20\%    	 	& 7					& 932				& 2624				\\
				$L=1/3$ 		& 10\%    		& 8  				& 1264 				& 3734				\\
				& 5\% 			& 8   				& 1520 				& 4660				\\
			\end{tabular}
			\endgroup
		}
	\end{center}
	\caption{(Blur) Comparison of cost to
		approximate the blur kernel from Equation
		\eqref{eq:frog_kernel} using the PSF-based method, the randomized
		HODLR (hierarchical off diagonal low rank) method, and GLR
		(global low rank) approximation using randomized SVD. The
		quantity $L$ scales the width of the impulse responses,
		hence it influences the rank of the operator. Large $L$
		means low rank, and small $L$ means high rank. The second
		column (``Error tol'') is the relative error in the
		approximation of the kernel measured in the Frobenius norm,
		$||\boldsymbol{\Phi} -
		\widetilde{\boldsymbol{\Phi}}||_\text{Fro} /
		||\boldsymbol{\Phi}||_\text{Fro}$. The remaining three
		columns show the number of operator applies required to
		achieve the given error tolerances, using the PSF, HODLR,
		and GLR methods.}
	\label{tab:frog_psf_vs_hodlr_vs_glr}
\end{table}

In Table~\ref{tab:frog_psf_vs_hodlr_vs_glr} we compare the
cost to approximate the blur kernel from
Equation~\eqref{eq:frog_kernel} using the PSF-based method, the
randomized HODLR (hierarchical off diagonal low rank)
method~\cite{Martinsson16,Hartland_2023} with 8 levels, and GLR (global low rank)
approximation using double-pass randomized SVD~\cite{HalkoMartinssonTropp11}. For
these results we vary the quantity $L$ to scale the width of the
impulse responses and hence the rank of the operator. For each case,
we calculate the relative error in the approximation of the kernel
measured in the Frobenius norm, $||\boldsymbol{\Phi} -
\widetilde{\boldsymbol{\Phi}}||_\text{Fro} /
||\boldsymbol{\Phi}||_\text{Fro}$ and show the number of operator
applications required to achieve 20\%, 10\%, and 5\% relative error by
each method. The results reveal superior performance of the PSF
method as compared to the HODLR and GLR methods for all cases. We note
that as we increase the rank and decrease the error tolerance, the
performance of the HODLR and GLR methods deteriorates.

\begin{figure}
	\begin{subfigure}{0.52\textwidth}
		\begin{center}
			\includegraphics[scale=1.0]{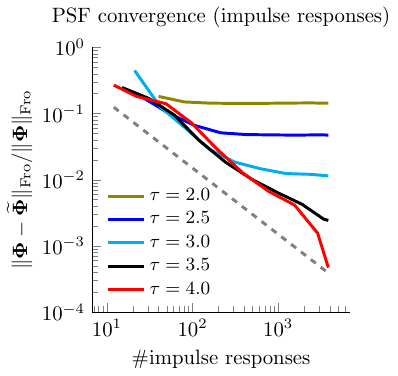}
		\end{center}
	\end{subfigure}
	\begin{subfigure}{0.47\textwidth}
		\begin{center}
			\includegraphics[scale=1.0]{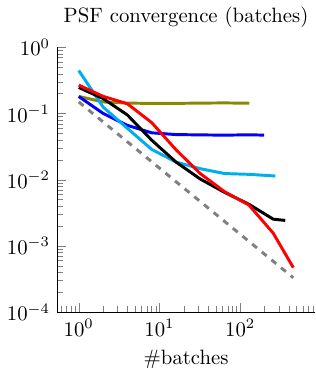}
		\end{center}
	\end{subfigure}
	\caption{(Blur) Relative error 
		for different ellipsoid size parameters, $\tau$, vs. the total number of impulse responses (left), and the number of batches (right). 
		The dashed gray lines show linear convergence rates, i.e., $\text{constant} \times \left(\# \text{impulse responses}\right)^{-1}$ on the left, and  $\text{constant} \times \# \text{batches}^{-1}$ on the right.}
	\label{fig:PSF_tau_convergence}
\end{figure} 

In Figure \ref{fig:PSF_tau_convergence} we show the convergence of the PSF-based method on the blur kernel as a function of the total number of impulse responses (left), and the number of batches (right). We show convergence for several ellipsoid size parameters $\tau$, ranging from $2.0$--$4.0$. The results in Figure \ref{fig:PSF_tau_convergence} (left) show that the relative error decreases as $\text{constant} \times \left(\# \text{impulse responses}\right)^{-1}$
suggesting linear convergence. The linear convergence stalls at a limit that depends on $\tau$. Increasing $\tau$ lowers this limit, allowing the PSF-based method to achieve higher accuracy. 
In Figure \ref{fig:PSF_tau_convergence} (right), the results show that before this limit is reached, the convergence is faster for smaller $\tau$ in terms of the number of batches. This is expected because smaller $\tau$ results in more impulse responses per batch. Larger $\tau$ causes the PSF-based method to converge more slowly than smaller $\tau$, but with larger $\tau$ the PSF-based method stalls at a lower level of error than it does with smaller $\tau$ (see, e.g., $\tau=4.0$ vs. $\tau=2.5$).

\section{Conclusions}
\label{sec:conclusions}

We presented an efficient matrix-free PSF-based method for approximating operators with locally supported non-negative integral kernels. The PSF-based method requires access to the operator only via application of the operator to a small number of vectors. The idea of the PSF-based method is to compute batches of impulse responses by applying the operator to Dirac combs of scattered point sources, then interpolate these impulse responses to approximate entries of the operator's integral kernel. The interpolation is based on a new principle we call ``local mean displacement invariance,'' which generalizes classical local translation invariance. The ability to quickly approximate arbitrary integral kernel entries permits us to form an H-matrix approximation of the operator. Fast H-matrix arithmetic is then used to perform further linear algebra operations that cannot be performed easily with the original operator, such as matrix factorization and inversion. The supports of the impulse responses are estimated to be contained in ellipsoids, which are determined a-priori via a moment method that involves applying the operator to a small number of polynomial functions. Point source locations for the impulse response batches are chosen using a greedy ellipsoid packing procedure, in which we choose as many impulse responses per batch as possible, while ensuring that the corresponding ellipsoids do not overlap. We applied the PSF-based method to approximate the Gauss-Newton Hessians in an ice sheet flow inverse problem governed by a linear Stokes PDE, and an advective-diffusive transport inverse problem governed by an advection-diffusion PDE. We saw that preconditioners based on the PSF-based approximation cluster the eigenvalues of the preconditioned Hessian near one, and allow us to solve the inverse problems using roughly $5\times$--$10\times$ fewer PDE solves. 
For larger domains with more observations, the rank of the data misfit Hessian will increase, while the locality of impulse responses will remain the same. Hence, we expect the speedup will be even greater for such problems. Although the PSF-based method is not applicable to all Hessians, it is applicable to many Hessians of practical interest. For these Hessians, the PSF-based method offers a \emph{data scalable} alternative to conventional low rank approximation, due to the ability to form high rank approximations of an operator using a small number of operator applications, and thus PDE solves.

\section*{Acknowledgments}
We thank J.J. Alger, Longfei Gao, Mathew Hu, and Rami Nammour for helpful discussions. We thank Nicole Aretz, Trevor Heise, and the anonymous reviewers for editing suggestions. We thank Georg Stadler for help with the domain setup for the ice sheet problem.

\FloatBarrier

\bibliographystyle{siamplain}
\bibliography{localpsf}

\end{document}